\documentclass[12pt,oneside,english]{amsart}
\usepackage[latin9]{inputenc}
\usepackage[letterpaper]{geometry}
\usepackage{babel}
\usepackage{textcomp}
\usepackage{mathrsfs}
\usepackage{amstext}
\usepackage{amsthm}
\usepackage{amssymb}
\usepackage[unicode=true,pdfusetitle,
 bookmarks=true,bookmarksnumbered=false,bookmarksopen=false,
 breaklinks=false,pdfborder={0 0 0},pdfborderstyle={},backref=false,colorlinks=false]
 {hyperref}

\makeatletter
\numberwithin{equation}{section}
\numberwithin{figure}{section}
\theoremstyle{plain}
\newtheorem*{thm*}{\protect\theoremname}
\theoremstyle{plain}
\newtheorem{thm}{\protect\theoremname}[section]
\theoremstyle{plain}
\newtheorem{prop}[thm]{\protect\propositionname}
\theoremstyle{plain}
\newtheorem{lem}[thm]{\protect\lemmaname}
\theoremstyle{definition}
\newtheorem{defn}[thm]{\protect\definitionname}
\theoremstyle{plain}
\newtheorem{cor}[thm]{\protect\corollaryname}
\theoremstyle{remark}
\newtheorem{claim}[thm]{\protect\claimname}
\theoremstyle{remark}
\newtheorem{rem}[thm]{\protect\remarkname}
\theoremstyle{remark}
\newtheorem*{claim*}{\protect\claimname}

\usepackage{tikz}
\usepackage{tikz-cd}
\usepackage{pgfplots}

\usepackage{amsmath}
\usepackage{amssymb}

\usepackage{enumitem}
\setlist{leftmargin=0.5cm}

\usepackage[cmtip,all]{xy}
\usepackage{quiver}

\input xypic

\newcommand{\ootimes}{
  \mathbin{
    \mathchoice
      {\buildcircleotimes{\displaystyle}}
      {\buildcircleotimes{\textstyle}}
      {\buildcircleotimes{\scriptstyle}}
      {\buildcircleotimes{\scriptscriptstyle}}
  }
}

\newcommand\buildcircleotimes[1]{%
  \begin{tikzpicture}[baseline=(X.base), inner sep=0, outer sep=0]
    \node[draw,circle] (X)  {$#1\otimes$};
  \end{tikzpicture}%
}

\usepackage{mathabx,graphicx}

\makeatother

\providecommand{\claimname}{Claim}
\providecommand{\corollaryname}{Corollary}
\providecommand{\definitionname}{Definition}
\providecommand{\lemmaname}{Lemma}
\providecommand{\propositionname}{Proposition}
\providecommand{\remarkname}{Remark}
\providecommand{\theoremname}{Theorem}

\begin{document}
\global\long\def\r{\mathbb{\mathbb{R}}}%

\global\long\def\n{\mathbb{\mathbb{N}}}%

\global\long\def\c{\mathbb{\mathbb{C}}}%

\global\long\def\b{\mathbb{\mathbb{B}}}%

\global\long\def\z{\mathbb{\mathbb{\mathbb{Z}}}}%

\global\long\def\q{\mathbb{\mathbb{\mathbb{Q}}}}%

\global\long\def\k{\textrm{\textbf{k}}}%

\global\long\def\h{\textrm{\textbf{H}}}%

\global\long\def\l{\ell}%

\global\long\def\o{\mathcal{O}}%

\global\long\def\nil{\emptyset}%

\global\long\def\longsquiggly{\xymatrix{{}\ar@{~>}[r]  &  {}}
 }%

\global\long\def\Hom{\mathrm{Hom}}%

\global\long\def\End{\mathrm{End}}%

\global\long\def\iso{\xrightarrow{{\sim}}}%

\global\long\def\idop{\textrm{Id}}%

\global\long\def\lra{\longrightarrow}%

\global\long\def\acts{\;\circlearrowright\;}%

\global\long\def\acted{\;\circlearrowleft\;}%

\global\long\def\dg{\mathsf{dg}}%

\global\long\def\op{\mathsf{op}}%

\global\long\def\vep{\varepsilon}%

\global\long\def\t{\tau}%

\global\long\def\xx{\tilde{x}}%

\global\long\def\tt{\tilde{\tau}}%

\global\long\def\ee{\tilde{E}}%

\global\long\def\veps{\varepsilon}%

\global\long\def\ci{\circ}%

\title{Tensor $2$-Product for $\mathfrak{sl}_{2}$: Extensions to the Negative
Half}
\author{Matthew McMillan}
\begin{abstract}
In a recent paper, the author defined an operation of tensor product
for a large class of $2$-representations of $\mathcal{U}^{+}$, the
positive half of the $2$-category associated to $\mathfrak{sl}_{2}$.
In this paper, we prove that the operation extends to give an operation
of tensor product for $2$-representations of the full $2$-category
$\mathcal{U}$: when the inputs are $2$-representations of the full
$\mathcal{U}$, the $2$-product is also a $2$-representation of
the full $\mathcal{U}$. As in the previous paper, the $2$-product
is given for a simple $2$-representation $\mathcal{L}(1)$ and an
abelian $2$-representation $\mathcal{V}$ taken from the $2$-category
of algebras.

This is the first construction of an operation of tensor product for
higher representations of a full Lie algebra in the abelian setting.
\end{abstract}

\maketitle
\tableofcontents{}

\section{Introduction}

\subsection{Background and motivation}

This paper is the second part in a series by the author, starting
with \cite{mcmillanTensor2product2representations2022}, about an
abelian tensor $2$-product operation for $2$-representations of
Lie algebras. This $2$-product is designed with a view to the program
by Crane and Frenkel \cite{craneFourDimensionalTopological1994} seeking
a higher representation theory in order to upgrade known $3d$ topological
invariants, such as the TQFT of Witten-Reshetikhin-Turaev \cite{wittenQuantumFieldTheory1989,reshetikhinInvariants3manifoldsLink1991},
to $4d$ invariants.

Prior work in this program involved building categories with Grothendieck
groups equal to various representations, including specific tensor
products, and these categories have been used to define homological
link invariants. This includes early work by Bernstein-Frenkel-Khovanov
\cite{bernsteinCategorificationTemperleyLiebAlgebra1999} and later
Stroppel and others \cite{stroppelCategorificationTemperleyLiebCategory2005,frenkelCategorificationFinitedimensionalIrreducible2007,mazorchukCombinatorialApproachFunctorial2009,sartoriCategorificationTensorProduct2015,sussanCategoryMathfrakslLink2007}
using category $\mathcal{O}$ of $\mathfrak{gl}_{n}$ for tensor products
of simples in type $A$, and work by Webster \cite{websterKnotInvariantsHigher2017,websterTensorProductAlgebras2016}
using diagrammatic methods for tensor products of simples in other
types. We expect these categories to be equivalent (in an appropriate
sense) to the tensor $2$-products of $2$-representations produced
by the operation studied in this paper when the factors are simple
$2$-representations.

The received notion of $2$-representation was introduced and developed
in \cite{chuangDerivedEquivalencesSymmetric2008a,laudaCategorificationQuantumMathfraksl2010,rouquier2KacMoodyAlgebras2008,khovanovDiagrammaticApproachCategorification2009,khovanovDiagrammaticApproachCategorification2011}.
A very general definition of tensor $2$-product operation for $2$-representations
of Kac-Moody algebras in the setting of $\mathcal{A}_{\infty}$-algebras
is expected from Rouquier in \cite{rouquierHopfCategories}. This
general definition does not come with explicit constructions. 

In \cite{mcmillanTensor2product2representations2022}, the present
author gave an explicit abelian model for the tensor $2$-product
$\mathcal{L}(1)\ootimes\mathcal{V}$ in the case of $\mathfrak{sl}_{2}^{+}$.
This is the construction of an algebra, bimodule, and bimodule maps
producing a $2$-action of $\mathcal{U}^{+}$, the positive half of
the $2$-category corresponding to the enveloping algebra of $\mathfrak{sl}_{2}$.
Here $\mathcal{L}(1)$ is a certain simple $2$-representation and
$\mathcal{V}$ is a given abelian $2$-representation taken from the
$2$-category of algebras and satisfying two additional hypotheses.

A related tensor $2$-product for the case of $\mathfrak{gl}(1|1)^{+}$,
which does not involve homotopical complications that are present
for $\mathfrak{sl}_{2}^{+}$ (due to the absence of endomorphisms
$x\in\End(E)$ in the relevant Hecke algebra), was applied by Manion-Rouquier
in \cite{manionHigherRepresentationsCornered2020} to describe Heegaard-Floer
theory for surfaces. Their construction has not been extended to the
full $\mathfrak{gl}(1|1)$.

It is not clear whether a $2$-representation theory for $\mathfrak{sl}_{2}^{+}$
could suffice to build a TQFT, and it is natural to ask whether the
construction in \cite{mcmillanTensor2product2representations2022}
can be extended to $\mathfrak{sl}_{2}$. The main result of this paper
is a proof that it can indeed be extended. It gives, then, the first
case of a $2$-product operation in the abelian setting for a \emph{full}
Lie algebra or super-Lie algebra, while \cite{mcmillanTensor2product2representations2022}
gave an operation for a \emph{half} Lie algebra (in an abelian setting),
whereas \cite{manionHigherRepresentationsCornered2020} used an operation
for a \emph{half super}-Lie algebra (in a $\dg$-setting).

\subsection{Result}

Let $\mathcal{U}$ be the $2$-category associated with the enveloping
algebra of $\mathfrak{sl}_{2}$, as given in Rouquier \cite[§4.1.3]{rouquier2KacMoodyAlgebras2008}
or Vera \cite[§3.2]{veraFaithfulnessSimple2representations2020}.
Let $\mathcal{U}^{+}$ be the monoidal category associated to the
positive half of the enveloping algebra of $\mathfrak{sl}_{2}$. As
in \cite[§1.2]{mcmillanTensor2product2representations2022}, we work
with $2$-representations in the abelian $2$-category of algebras,
bimodules, and bimodule maps.

Let $A$ be a $k$-algebra for a field $k$, let $E$ be an $(A,A)$-bimodule,
and let $x\in\End(E)$, $\tau\in\End(E^{2})$ be bimodule endomorphisms
satisfying the nil affine Hecke relations: 
\begin{gather}
\tau^{2}=0,\nonumber \\
\tau E\circ E\tau\circ\tau E=E\tau\circ\tau E\circ E\tau,\label{eq:nil-aff-Hecke-rels}\\
\tau\circ Ex=xE\circ\tau+1,\;Ex\circ\tau=\tau\circ xE+1.\nonumber 
\end{gather}
(The notation $xE$ indicates the endomorphism $x\otimes\idop_{E}$
in $\End(E^{2})$, etc.) The data $(A,E,x,\tau)$ determines a $2$-representation
of $\mathcal{U}^{+}$.

Now assume that $(A,E,x,\tau)$ has a weight decomposition $A=\prod_{\lambda\in\z}A_{\lambda}$
(cf.~\cite[§4.3.1]{mcmillanTensor2product2representations2022}).
The data $(A,E,x,\tau)$ extends to determine a $2$-represen\-tation
of the full $2$-category $\mathcal{U}$ when the functor $E\otimes_{A}-$
admits a right adjoint functor $F$ (with unit $\eta$ and counit
$\veps$) such that the \textquotedblleft commutator'' maps $\rho_{\lambda}$
(determined by $x$, $\tau$, $\eta$, $\veps$; see §\ref{subsec:Commutator-defs}
below) are isomorphisms in each weight $\lambda\in\z$.

A simple $2$-representation $\mathcal{L}(1)$ of $\mathcal{U}$ that
categorifies the fundamental representation $L(1)$ of $\mathfrak{sl}_{2}$
may be given by the following data. Let the $k$-algebra be $k[y]_{+1}\times k[y]_{-1}$
(decomposed into weight algebras), and the triple be $(k[y],y,0)$.
Let $x$ act by multiplication by $y$. Let $y\in k[y]_{-1}$ act
on $k[y]$ on the right by multiplication, and $y\in k[y]_{+1}$ act
by zero; swap them for the left action.

Let $P_{n}=k[x_{1},\dots,x_{n}]$ be the polynomial algebra. Then
$P_{n}$ acts on $E^{n}$ with $x_{i}\in P_{n}$ acting by the endomorphism
$E^{n-i}xE^{i-1}$.
\begin{thm*}[Main result]
 Suppose $(A,E,x,\tau)$ gives the data of a $2$-representation
$\mathcal{V}$ of $\mathcal{U}^{+}$ such that $\mathcal{V}$ has
a weight decomposition. Define the left-dual $(A,A)$-bimodule $F=\Hom_{A}(_{A}E,A)$.
Suppose $E$ has the following properties: 
\begin{itemize}
\item $_{A}E$ is finitely generated and projective, so $(E\otimes_{A}-,F\otimes_{A}-)$
is an adjunction where the unit $\eta$ and counit $\veps$ arise
from the duality pairing, 
\item $E^{n}$ is free as a $P_{n}$-module, 
\item $E$ and $F$ are locally nilpotent, 
\item The maps $\rho_{\lambda}$ determined by the given data are isomorphisms
for each $\lambda\in\z$.
\end{itemize}
\raggedright These properties imply that $(A,E,F,x,\tau,\eta,\veps)$
determines an integrable $2$-representation of $\mathcal{U}$.

\medskip{}
Now let $C$ be the $k$-algebra, $\tilde{E}$ the $(C,C)$-bimodule,
and $\tilde{x}$ and $\tilde{\tau}$ the bimodule endomorphisms constructed
in \cite{mcmillanTensor2product2representations2022}. Let $\tilde{F}=\Hom_{C}(_{C}\tilde{E},C)$.
Then: 
\begin{itemize}
\item $_{A}\tilde{E}$ is finitely generated and projective, so $(\tilde{E}\otimes_{C}-,\tilde{F}\otimes_{C}-)$
is an adjunction with unit $\tilde{\eta}$ and counit $\tilde{\veps}$
arising from the duality pairing, 
\item $\tilde{E}$ and $\tilde{F}$ are locally nilpotent, 
\item The maps $\tilde{\rho}_{\lambda}$ determined by the given data are
isomorphisms, so: 
\item $(C,\tilde{E},\tilde{F},\tilde{x},\tilde{\tau},\tilde{\eta},\tilde{\veps})$
determines an integrable $2$-representation of $\mathcal{U}$.
\end{itemize}
\end{thm*}
The data $(C,\tilde{E},\tilde{x},\tilde{\tau})$ determines a $2$-representation
of $\mathcal{U}^{+}$ that we interpreted in \cite{mcmillanTensor2product2representations2022}
as the result $\mathcal{L}(1)\ootimes\mathcal{V}$ of a $2$-product
operation (with the factors considered as $2$-representations of
$\mathcal{U}^{+}$). One reason to interpret the structure in this
way was that it results from a categorification of the Hopf coproduct
formula. Another reason was that it recovers the expected structures
in some known cases. For details, see \cite{mcmillanTensor2product2representations2022}
in §1.3, §1.4, as well as in Remark 3.4 about the effect of $E'$
and thus $\tilde{E}$ on the Grothendieck group. Since the additional
components $\tilde{F}$, $\tilde{\eta}$, and $\tilde{\veps}$ are
fully determined by $(C,\tilde{E},\tilde{x},\tilde{\tau})$, in this
article we interpret the $2$-representation determined by the combined
data $(C,\tilde{E},\tilde{F},\tilde{x},\tilde{\tau},\tilde{\eta},\tilde{\veps})$
as the result $\mathcal{L}(1)\ootimes\mathcal{V}$ of a $2$-product
operation with the factors considered as $2$-representations of $\mathcal{U}$.

We emphasize that for an integrable $2$-representation of $\mathcal{U}^{+}$
given by the data $(A,E,x,\tau)$, the fact that the data determines
a $2$-representation of the full $2$-category $\mathcal{U}$ is
equivalent to the data having a property: namely that $_{A}E$ is
f.g.~projective, and the commutator maps $\rho_{\lambda}$ determined
by the data are isomorphisms. When this holds, then (according to
the theorem) the maps $\tilde{\rho}_{\lambda}$ of the product are
also isomorphisms. So the new data $(C,\tilde{E},\tilde{x},\tilde{\tau})$
inherits the property of determining an action of the full $\mathcal{U}$.

\subsection{Outline summary}

The paper is organized as follows:
\begin{itemize}
\item In §2 we introduce the relevant background theory for extensions of
$2$-representations of $\mathcal{U}^{+}$ to the full $2$-category
$\mathcal{U}$. This section builds on the background theory and definitions
of \cite{mcmillanTensor2product2representations2022}. We include
a discussion of the adjunction, the commutator maps $\rho_{\lambda}$,
and the condition of integrability as it relates to our product construction.
\item In §3 we define and study important bimodules, giving concrete algebraic
models for them in the manner of §3.2 of \cite{mcmillanTensor2product2representations2022}.
\item In §4 we consider the left dual to $\tilde{E}$, namely $\tilde{F}=\Hom_{C}(_{C}\tilde{E},C)$,
and we show how to describe it concretely by using the $B$ side of
the equivalence described in §3.3.2 of \cite{mcmillanTensor2product2representations2022}.
\item In §5.1 we study the tensor products $\tilde{E}\otimes_{C}\tilde{E}$
and $\tilde{E}\otimes_{C}\tilde{F}$ and $\tilde{F}\otimes_{C}\tilde{E}$,
and describe their structure as $(A[y],A[y])$-bimodules. In §5.2
we compute explicit formulas for $\tilde{\rho}_{\lambda}$ in terms
of the structures found in §5.1. In §5.3 we use the formulas from
§5.2 to show that each $\tilde{\rho}_{\lambda}$ is an isomorphism.
\end{itemize}

\subsection{Acknowledgments}

I thank Raphaël Rouquier for advice and encouragement during this
project. This work was supported by the NSF through grant DMS-1702305.

\section{Background: extending $\mathcal{U}^{+}$ actions to $\mathcal{U}$
actions}

\subsection{$2$-Representations of $\mathcal{U}$}

We begin with a description of a $2$-represen\-tation of the full
$2$-category $\mathcal{U}$ associated to the Lie algebra $\mathfrak{sl}_{2}$.
The $2$-category $\mathcal{U}$ that we mean is defined in §4.1.3
of \cite{rouquier2KacMoodyAlgebras2008}, but with $\tau$ replaced
by $-\tau$ in the Hecke relations. We do not repeat that definition
here since we work with the concrete data of $2$-representations
and not with the $2$-category $\mathcal{U}$ itself.

In \cite{mcmillanTensor2product2representations2022}, a $2$-representation
was defined as a strict monoidal functor from $\mathcal{U}^{+}$ to
a monoidal category of the form $\mathsf{Bim}_{k}(A)$, which is defined
for a $k$-algebra $A$ as follows: the objects of $\mathsf{Bim}_{k}(A)$
are $(A,A$)-bimodules, and the morphisms of $\mathsf{Bim}_{k}(A)$
are bimodule maps. The monoidal structure on $\mathsf{Bim}_{k}(A)$
is given by tensor product of bimodules over $A$. This monoidal category
$\mathsf{Bim}_{k}(A)$ may also be interpreted as a $2$-category
with a single object $A$, where the $1$-morphisms are given by tensor
product with $(A,A)$-bimodules, and the $2$-morphisms are bimodule
maps.

A $2$-representation of the full $\mathcal{U}$ is defined in terms
of weights (see Def.~4.25 of \cite{mcmillanTensor2product2representations2022}).
When $A$ is provided with a weight decomposition $A=\prod_{\lambda\in\z}A_{\lambda}$,
then the $2$-category $\mathsf{Bim}_{k}(A)$ with single object $A$
may be expanded to a $2$-category with objects given by the weight
algebras $A_{\lambda}$, morphisms given by $(A_{\mu},A_{\lambda})$-bimodules,
and $2$-morphisms given by bimodule maps. With this interpretation,
a $2$-representation of $\mathcal{U}$ may be described as a strict
$2$-functor $\mathcal{U}\to\mathsf{Bim}_{k}(A)$ given on objects
by $\mathbf{1}_{\lambda}\mapsto A_{\lambda}$.

According to Prop.~2.4 of \cite{mcmillanTensor2product2representations2022},
a $2$-representation of $\mathcal{U}^{+}$ in $\mathsf{Bim}_{k}(A)$
is equivalent to the data of a $k$-algebra $A$ together with a bimodule
$_{A}E_{A}$ and bimodule maps $x\in\End(E)$, $\tau\in\End(E^{2})$
satisfying relations (\ref{eq:nil-aff-Hecke-rels}). This paper will
rely on the following analogue of that proposition: 
\begin{prop}
\label{prop:def-2rep-full-U} The data of a $2$-representation $\mathcal{U}\to\mathsf{Bim}_{k}(A)$
for a $k$-algebra $A=\prod_{\lambda\in\z}A_{\lambda}$ consists of
bimodules $_{A}E_{A}$, $_{A}F_{A}$ (having weights $+2$ and $-2$),
the unit $\eta$ and counit $\veps$ of an adjunction $(E,F)$, and
bimodule maps $x\in\End(E)$, $\tau\in\End(E^{2})$ that satisfy relations
(\ref{eq:nil-aff-Hecke-rels}), all such that $\rho_{\lambda}$ (defined
below in terms of $x$, $\tau$, $\eta$, $\veps$) is an isomorphism
for each $\lambda$.
\end{prop}

(Bimodules $E$, $F$ are said to have weight $+2$ and $-2$, respectively,
when $e_{j}Ee_{i}=\delta_{i+2,j}\cdot e_{i+2}Ee_{i}$ and $e_{j}Fe_{i}=\delta_{i-2,j}\cdot e_{i-2}Fe_{i}$.)

In this paper, a symbol $\mathcal{V}$ is used sometimes to denote
a $2$-representation of $\mathcal{U}^{+}$, and sometimes to denote
the extension of the former to a $2$-representation of $\mathcal{U}$.
This is an abuse of notation because the first $\mathcal{V}$ is a
monoidal category, and the second $\mathcal{V}$ is a $2$-category.
This abuse is justifiable when both types of category are determined
by the same data.

\subsection{\label{subsec:Commutator-defs}Commutator morphisms}

Here we define the commutator morphisms. Assume we are given the data
of a $k$-algebra $A$, bimodules $_{A}E_{A}$, $_{A}F_{A}$ which
determine endofunctors of $A\text{-mod}$ by tensor product on the
left, the unit $\eta$ and counit $\veps$ of an adjunction $(E,F)$,
and endofunctors $x$ and $\tau$ satisfying (\ref{eq:nil-aff-Hecke-rels}).
Assume that $A$ has a weight decomposition $A=\prod_{\lambda\in\z}A_{\lambda}$,
and $E$ and $F$ have weights $+2$ and $-2$. Let us use the notation
$E_{\lambda}=E\cdot A_{\lambda}$ and $_{\mu}E_{\lambda}=A_{\mu}\cdot E\cdot A_{\lambda}$,
so $E=\bigoplus_{\mu,\lambda}\phantom{}_{\mu}E_{\lambda}.$ In this
paper we also use a convention that \textquoteleft $\oplus$\textquoteright{}
and \textquoteleft $\sum$\textquoteright{} denote the components
of a map to and from a direct sum, respectively.

We define $\sigma:EF\to FE$ by: 
\[
\sigma=(FE\varepsilon)\circ(F\tau F)\circ(\eta EF):EF\to FE.
\]
For $\lambda\in\z_{\geq0}$ we define: 

\begin{equation}
\rho_{\lambda}=\sigma\oplus\bigoplus_{i=0}^{\lambda-1}\varepsilon\circ x^{i}F:EF_{\lambda}\to FE_{\lambda}\oplus A_{\lambda}^{\oplus\lambda},\label{eq:rho_+}
\end{equation}
 and for $\lambda\in\z_{\leq0}$: 
\begin{equation}
\rho_{\lambda}=\left(\sigma,\sum_{i=0}^{-\lambda-1}Fx^{i}\circ\eta\right):EF_{\lambda}\oplus A_{\lambda}^{\oplus-\lambda}\to FE_{\lambda}.\label{eq:rho_-}
\end{equation}
(The summation terms are neglected when $\lambda=0$.)

\subsection{Conventions}

We adopt the conventions of \cite{mcmillanTensor2product2representations2022},
so the reader may consult §2.3 of that text for additional details.
Assume we are given data $(A,E,x,\tau)$ determining a $2$-representation
$\mathcal{V}$ of $\mathcal{U}^{+}$. Assume that $_{A}E$ is f.g.~projective
and that $E^{n}$ is free as a $P_{n}$-module.

The construction of the product $\mathcal{L}(1)\ootimes\mathcal{V}$
in \cite{mcmillanTensor2product2representations2022} makes use of
the $(A[y],A[y])$-bimodule $E[y]$, and the endomorphism $x-y\in\End(E[y])$.
Write $E_{y}$ for the quotient $E[y]\big/(x-y)E[y]$, and $\pi:E[y]\to E_{y}$
for its projection.

Concatenation of the symbols for bimodules indicates tensor product
over some algebra that is determined by context. Sometimes this algebra
could be either $A$ or $A[y]$, so we stipulate that if the expression
for a bimodule contains \textquoteleft $y$\textquoteright , it will
be understood as an $(A[y],A[y])$-bimodule, and if the expression
lacks \textquoteleft $y$\textquoteright , it will be understood as
an $A$-module. We suppress isomorphisms such as: 
\[
E[y]E_{y}=E[y]\otimes_{A[y]}E_{y}\iso E\otimes_{A}E_{y}=EE_{y}.
\]

Extend $x$ to an element of $\End(E[y])$ by $x:ey^{n}\mapsto x(e)y^{n}$
and $\tau$ to $\End(E^{2}[y])$ by $\tau:eey^{n}\mapsto\tau(ee)y^{n}$.
When writing formulas for morphisms we often write an arbitrary element
of $E[y]$ with the single letter \textquoteleft $e$\textquoteright{}
and an arbitrary element of $E^{2}[y]$ with the doubled symbol \textquoteleft $ee$\textquoteright{}
(which is not assumed to be a simple tensor).

We make use of the notation $y_{i}=x_{i}-y$. Here $y_{i}$ indicates
$\bigl(E^{j}xE^{i-1}-y\bigr)$ for some $j$, and context will determine
the value of $j$.

As in §2.3 of \cite{mcmillanTensor2product2representations2022},
let $s\in\mathrm{End}(E^{2})$ be the bimodule map given by $s=\tau\ci(x_{1}-x_{2})-\idop$,
and extended to $E^{2}[y]$ as $x$ and $\tau$ are extended. Note
that $s$ descends to define maps of $(A[y],A[y])$-bimodules $s:E_{y}E\to EE_{y}$
and $s:EE_{y}\to E_{y}E$ such that $s^{2}$ descends to $\idop$.

\subsection{Adding a dual}

Every bimodule $_{A}E_{A}$ has left- and right-dual bimodules, 
\begin{align*}
\phantom{}^{\vee\negmedspace}E & =\Hom_{A}(_{A}E,A),\\
E^{\vee} & =\Hom_{A}(E_{A},A),
\end{align*}
 respectively. 

Now, when $_{A}E$ is f.g.~projective, the canonical morphism $\phantom{}^{\vee\negmedspace}E\otimes_{A}E\to\Hom_{A}(_{A}E,E)$
is an isomorphism of bimodules. More generally, the canonical morphism
of functors $\phantom{}^{\vee\negmedspace}E\otimes_{A}-\to\Hom_{A}(_{A}E,-)$
is an isomorphism. In this situation, the endofunctor $\phantom{}^{\vee\negmedspace}E\otimes_{A}-$
of the category $A\text{-mod}$ is right adjoint to the endofunctor
$E\otimes_{A}-$ of the same category. The triple $(\phantom{}^{\vee\negmedspace}E,\eta,\veps)$
gives the right-dual object for $E$ in the monoidal category $\mathsf{Bim}_{k}(A)$.
Here $\veps:E\otimes_{A}\phantom{}^{\vee\negmedspace}E\to A$ is given
by evaluation, and $\eta:A\to\phantom{}^{\vee\negmedspace}E\otimes_{A}E$
is given via the isomorphism $\phantom{}^{\vee\negmedspace}E\otimes_{A}E\iso\Hom_{A}(_{A}E,E)$
by the right $A$-action whereby $\eta(a):e\mapsto e.a$. (Note that
we say $\phantom{}^{\vee\negmedspace}E$ is the \emph{left}-dual bimodule,
even though it gives the \emph{right}-dual object.)

Conversely, assume that $(E\otimes_{A}-,\phantom{}^{\vee\negmedspace}E\otimes_{A}-)$
is an adjoint pair for some bimodule $_{A}E_{A}$. The adjunction
gives equivalences of functors: 
\[
\Hom_{A}(_{A}E,-)\cong\Hom_{A}(_{A}A,\phantom{}^{\vee\negmedspace}E\otimes_{A}-)\cong\phantom{}^{\vee\negmedspace}E\otimes_{A}-,
\]
so all three are both right- and left-exact functors. So $_{A}E$
is projective. Furthermore, these functors commute with infinite direct
sums, so $_{A}E$ is finitely generated as well.

In this paper we consider $2$-representations for which the image
of $F$ in $\mathsf{Bim}_{k}(A)$, i.e.~the bimodule $_{A}F_{A}$,
is identically the left-dual bimodule $\phantom{}^{\vee\negmedspace}E$.
There is no loss of generality because any $2$-representation of
$\mathcal{U}$ in $\mathsf{Bim}_{k}(A)$ is equivalent to one of these.
(For any $2$-representation in $\mathsf{Bim}_{k}(A)$, the endofunctor
$_{A}F\otimes_{A}-$ of $A\text{-mod}$ is right adjoint to $_{A}E\otimes_{A}-$,
and is therefore unique up to unique isomorphism.) A $2$-representation
of $\mathcal{U}$ given by the data $(A,E,F,x,\tau,\eta,\veps)$ in
$\mathsf{Bim}_{k}(A)$ is said to \emph{extend} a $2$-representation
$(A,E,x,\tau)$ of $\mathcal{U}^{+}$ when $F=\phantom{}^{\vee\negmedspace}E$
and $\eta$, $\veps$ arise from the duality.

It was a hypothesis of the main theorem of \cite{mcmillanTensor2product2representations2022}
that $_{A}E$ is f.g.~projective. This condition was needed in order
to show that $E'X$ was a perfect complex (for example). In light
of the above, we see that the existence of an extension of the $2$-representation
of $\mathcal{U}^{+}$ to a $2$-representation of $\mathcal{U}$ in
$\mathsf{Bim}_{k}(A)$ also necessitates that hypothesis.

The following lemma is a consequence of the foregoing discussion.
\begin{lem}
Suppose the data $(A,E,x,\tau)$ determines a $2$-representation
of $\mathcal{U}^{+}$ in $\mathsf{Bim}_{k}(A)$ having a weight decomposition.
This data extends to determine a $2$-representation of $\mathcal{U}$,
with $F=\phantom{}^{\vee\negmedspace}E$, if and only if $_{A}E$
is f.g.~projective and the commutator morphisms $\rho_{\lambda}$
(determined by $x$, $\tau$, $\eta$, $\veps$) are isomorphisms.
\end{lem}

In \cite{mcmillanTensor2product2representations2022} the author defined
the data $(C,\tilde{E},\tilde{x},\tilde{\tau})$ of a product $2$-represen\-tation
of $\mathcal{U}^{+}$ in terms of given data $(A,E,x,\tau)$ satisfying
some conditions. In that paper it was seen that $_{C}\tilde{E}$ is
f.g.~projective, and it follows that $\tilde{F}=\phantom{}^{\vee\negmedspace}\tilde{E}$
is right adjoint to $\tilde{E}$. In this paper we aim to show that
$(C,\tilde{E},\tilde{F},\tilde{x},\tilde{\tau},\tilde{\eta},\tilde{\veps})$
determines a $2$-representation of $\mathcal{U}$. Our argument uses
the above Lemma: it will suffice to show that the commutator morphisms
$\tilde{\rho}_{\lambda}$ are isomorphisms.

\subsection{Integrability}

In the literature, a $2$-representation is typically defined in terms
of weight categories $\mathcal{C}_{\lambda}$ and functors $E$ and
$F$ between them, whereas we have framed our results entirely in
terms of bimodules $E$ and $F$. One reason for this is that a certain
pair of bimodules may determine several functors (by the operation
of tensoring on the left) that act on several reasonable categories
of modules. The most important ones are $A\text{-mod}$ and $A\text{-proj}$.

The distinction between $A\text{-mod}$ and $A\text{-proj}$ interacts
with our results and the hypothesis of integrability in an interesting
way. This interaction is mediated by the property of \textquotedblleft second
adjunction\textquotedblright{} that a $2$-representation of $\mathcal{U}$
may possess. We explain this next. Note that some authors include
the second adjunction in their definition of a $2$-representation,
and for them, this discussion will be of minor significance. It may
be interesting for them to observe, though, that in our construction
of tensor product, the hypothesis of integrability passes from the
factors to the product quite easily, while it is not clear that a
second adjunction alone passes from the factors to the product at
all.

Every $2$-representation of $\mathcal{U}$ given with functors $E$
and $F$ comes with one adjunction $(E,F)$, and with the data of
a \textquotedblleft candidate\textquotedblright{} unit and counit
pair for a second adjunction $(F,E)$. When the $2$-representation
acts on a category $A\text{-mod}$ and $E$ and $F$ are given by
tensoring with bimodules, the first adjunction implies that $_{A}E$
is f.g.~projective. In this case, the upper half $\mathcal{U}^{+}$
also acts on the smaller category $A\text{-proj}$. If the $2$-representation
is assumed to be integrable, and the full $\mathcal{U}$ acts, i.e.~the
$\rho_{\lambda}$ are isomorphisms, then by Theorem 5.16 of \cite{rouquier2KacMoodyAlgebras2008}
the given candidates do provide a second adjunction $(F,E)$. This
adjunction implies that $_{A}F$ is also f.g.~projective, and now
the full $\mathcal{U}$ action may be restricted to $A\text{-proj}$.

Given only the first adjunction with an action of $\mathcal{U}^{+}$,
so $_{A}E$ is f.g.~projective, together with the hypothesis that
$E^{n}$ is free over $P_{n}$, we can form the $2$-representation
of $\mathcal{U}^{+}$ called $\mathcal{L}(1)\ootimes\mathcal{V}$
in \cite{mcmillanTensor2product2representations2022}. In that paper
it was shown that $_{C}\tilde{E}$ is f.g.~projective, so it may
be interpreted either in an action on $C\text{-mod}$ or else in an
action restricted to $C\text{-proj}$. Given also a second adjunction
$(^{\vee\negmedspace}E,E)$ determining an action of the full $\mathcal{U}$,
we know that $\mathcal{U}$ acts on $A\text{-proj}$ through $E$
and $^{\vee\negmedspace}E$ in the $2$-representation $\mathcal{V}$,
but we are not (currently) able to show from this alone that $\mathcal{U}$
acts on $C\text{-proj}$ through $\tilde{E}$ and $^{\vee\negmedspace}\tilde{E}$,
since we do not know that $^{\vee\negmedspace}\tilde{E}$ is f.g.~projective.

Given the first adjunction $(E,\phantom{}^{\vee\negmedspace}E)$ and
also the hypothesis of integrability of an action of the full $\mathcal{U}$,
we know that there is a second adjunction $(^{\vee\negmedspace}E,E)$.
Now the hypothesis of integrability itself passes to the product bimodule
$\tilde{E}$. (Prop.~4.24 of \cite{mcmillanTensor2product2representations2022}.)
Given that we can also show that the product maps $\tilde{\rho}_{\lambda}$
are isomorphisms (the main effort of this paper), so we have an action
of the full $\mathcal{U}$ on $C\text{-mod}$, it follows from integrability
that there is a second adjunction $(^{\vee\negmedspace}\tilde{E},\tilde{E})$
for the product. This implies, in turn, that $_{C}{}^{\vee\negmedspace}\tilde{E}$
is f.g.~projective and that the full $\mathcal{U}$ action may be
restricted to the category $C\text{-proj}$.

To summarize, second adjunctions enable restriction of the full $\mathcal{U}$
action to the subcategories $A\text{-proj}$ and $C\text{-proj}$.
The existence of a second adjunction $(^{\vee\negmedspace}E,E)$ in
$\mathcal{V}$ is \emph{not enough} (with the arguments below) to
guarantee a second adjunction $(^{\vee\negmedspace}\tilde{E},\tilde{E})$
in $\mathcal{L}(1)\ootimes\mathcal{V}$. But \emph{integrability}
of $\mathcal{V}$ is enough to guarantee \emph{integrability} of $\mathcal{L}(1)\ootimes\mathcal{V}$,
as well as to give \emph{both} second adjunctions $(^{\vee\negmedspace}E,E)$
and $(^{\vee\negmedspace}\tilde{E},\tilde{E})$.

\subsection{Background: $2$-product for $\mathcal{U}^{+}$}

We recall some definitions and results from \cite{mcmillanTensor2product2representations2022}.
The reader is encouraged to review that paper and to consult it for
additional details and conventions.
\begin{defn}[Def.~3.1 of \cite{mcmillanTensor2product2representations2022}]
 Let $B$ be the $k$-algebra:
\[
B=\begin{pmatrix}A[y] & E_{y}\\
0 & A[y]
\end{pmatrix}.
\]
The algebra structure of $B$ is given by matrix multiplication, where
products involving generators in $[B]_{12}$ are defined using the
bimodule structure of $E_{y}$.
\end{defn}

Modules over $B$ are naturally written in terms of components. A
left $B$-module is given by a pair $\left(\begin{smallmatrix}M_{1}\\
M_{2}
\end{smallmatrix}\right)$ of left $A[y]$-modules, together with a morphism $\alpha:E_{y}\otimes_{A[y]}M_{2}\to M_{1}$
of left $A[y]$-modules specifying the action of $E_{y}$ generators;
analogously a right $B$-module is given by a pair $\left(\begin{smallmatrix}N_{1} & N_{2}\end{smallmatrix}\right)$
and morphism $\beta:N_{1}\otimes_{A[y]}E_{y}\to N_{2}$.

A bimodule consists of a $2\times2$ matrix with additional data.
The direct sum of coefficients in the top row of such a matrix gives
the top component of the pair corresponding to the left-module structure,
and the bottom row gives the bottom component of the pair; similarly
the columns give the components of the right-module structure. The
additional data consists of $\alpha$ determining \textquoteleft vertical\textquoteright{}
maps and $\beta$ giving \textquoteleft horizontal\textquoteright{}
maps. A matrix together with maps $\alpha$ and $\beta$ determines
a bimodule only if the left and right actions of $E_{y}$ specified
by $\alpha$ and $\beta$ commute. (In this situation the vertical
and horizontal maps respect the decompositions into horizontal and
vertical components, respectively.)

A complex of left $B$-modules is equivalent to a pair of complexes
of $A[y]$-modules and a map of complexes $\alpha$; analogously for
right $B$-modules and for bimodules. Complexes in this paper have
a cohomological grading.
\begin{defn}[Def.~3.2 of \cite{mcmillanTensor2product2representations2022}]
 \label{def:E'} Let $E'$ be the complex of $(B,B)$-bimodules that
is nonzero in degrees $0$ and $1$, where it is given by: 
\[
E'_{0}=\begin{pmatrix}E[y] & E[y]E_{y}\\
0 & E[y]
\end{pmatrix},\;E'_{1}=\begin{pmatrix}E_{y} & E_{y}E_{y}\\
A[y] & E_{y}
\end{pmatrix}.
\]

The left action of a generator in $E_{y}\subset B$ is specified on
vertical columns of $E'_{0}$ by the maps $0:E_{y}\otimes0\to E[y]$
and $s:E_{y}E[y]\to E[y]E_{y}$. The left action on $E'_{1}$ is specified
by the identity map on vertical columns. The right action on $E'_{0}$
is specified by the identity on the top row, and $0$ on the bottom
row. The right action on $E'_{1}$ is specified by the identity on
both rows. The differential $E'_{0}\to E'_{1}$ is given \emph{componentwise}
by $\left(\begin{smallmatrix}\pi & \pi E_{y}\\
0 & \pi
\end{smallmatrix}\right)$.
\end{defn}

\begin{lem}[Lemma 3.3 of \cite{mcmillanTensor2product2representations2022}]
 Let $M=\left(\left(\begin{smallmatrix}M_{1}\\
M_{2}
\end{smallmatrix}\right),\alpha\right)$ be a complex of left $B$-modules (written as a pair of complexes),
where $\alpha:E_{y}\otimes_{A[y]}M_{2}\to M_{1}$ specifies the action
for generators in $E_{y}$. The functor $E'\otimes_{B}-$ on $M$
may be given by: 
\[
\left(\begin{pmatrix}M_{1}\\
M_{2}
\end{pmatrix},\alpha\right)\overset{E'}{\longmapsto}\left(\begin{pmatrix}E[y]M_{1}\overset{\overset{\pi M_{1}}{\curvearrowright}}{\oplus}E_{y}M_{1}[-1]\\
E[y]M_{2}\overset{\overset{\alpha\circ\pi M_{2}}{\curvearrowright}}{\oplus}M_{1}[-1]
\end{pmatrix},\begin{pmatrix}E[y]\alpha\circ sM_{2} & 0\\
0 & \idop_{E_{y}M_{1}}
\end{pmatrix}\right).
\]
Here the top and bottom rows express cocones of the maps $\pi M_{1}$
and $\alpha\circ\pi M_{2}$.
\end{lem}

\begin{defn}[Def.~3.5 of \cite{mcmillanTensor2product2representations2022}]
 Let $X$ be the following complex of $B$-modules:
\begin{align*}
X & =X_{1}\oplus X_{2};\quad X_{1}=\begin{pmatrix}A[y]\\
0
\end{pmatrix},\quad X_{2}=E'X_{1}=\begin{pmatrix}E[y] & \overset{\pi}{\lra} & E_{y}\\
0 & \lra & A[y]
\end{pmatrix},
\end{align*}
where $X_{1}$ lies in degree $0$ and $X_{2}$ in degrees $0$ and
$1$. The $E_{y}$ action on $X_{2}$ is given by $E_{y}\otimes_{A[y]}A[y]\xrightarrow{\sim}E_{y}$,
$e\otimes1\mapsto e$.
\end{defn}

\begin{prop}[Prop.~3.6 of \cite{mcmillanTensor2product2representations2022}]
 The complex $X$ is strictly perfect and generates $\text{per }B$,
the full subcategory of $D^{b}(B)$ of complexes quasi-isomorphic
to strictly perfect complexes.
\end{prop}

Next we recall an important series of bimodules introduced in \cite{mcmillanTensor2product2representations2022}: 
\begin{defn}[Def.~3.16 of \cite{mcmillanTensor2product2representations2022}]
 Let $G_{n}$ denote $\Hom_{K^{b}(B)}(X_{2},E'^{n}X_{1})$.
\end{defn}

Note that the quotient projection to the derived category is an isomorphism
$G_{n}\iso\mathrm{Hom}_{D^{b}(B)}(X_{2},E'^{n}X_{1})$ because $X_{2}$
is strictly perfect. Note also that $G_{1}$ has an algebra structure
given by composition of endomorphisms.
\begin{prop}[Props.~3.18, 3.20, and 3.22 together with 3.27 and 3.28 of \cite{mcmillanTensor2product2representations2022}]
 \label{prop:G-bar_i} There are isomorphisms of $(A[y],A[y])$-bimodules
$\bar{G}_{1}\iso G_{1}$, $\bar{G}_{2}\iso G_{2}$, $\bar{G}_{3}\iso G_{3}$,
where: 
\begin{gather*}
\bar{G}_{1}=\biggl\langle(\theta,\varphi)\in A^{\op}[y]\oplus\End_{A}(_{A}E)[y]\biggr|\\
\begin{split}\varphi & =\_.\theta+y_{1}\varphi_{1}\\
 & \text{ for some }\varphi_{1}\in\End_{A}(_{A}E)[y]\biggr\rangle,
\end{split}
\end{gather*}
\begin{gather*}
\bar{G}_{2}=\biggl\langle(e_{1},e_{2},\xi)\in E[y]^{\oplus2}\oplus\Hom_{A}(_{A}E,E^{2})[y]\biggr|\\
\begin{split}\qquad e_{1}-e_{2} & =y_{1}e'\\
\xi & =\_\otimes e_{1}+y_{2}\xi_{1}\\
\xi_{1} & =\tau(\_\otimes e_{2})+y_{1}\xi'\\
 & \text{ for some }e'\in E[y]\text{ and }\xi'\in\Hom_{A}(_{A}E,E^{2})[y]\biggr\rangle,
\end{split}
\end{gather*}
 
\begin{gather*}
\bar{G}_{3}=\biggl\langle(ee_{1},ee_{2},ee_{3},\chi)\in E^{2}[y]^{\oplus3}\oplus\Hom_{A}(_{A}E,E^{3})[y]\biggr|\\
\begin{split}\qquad ee_{1}-ee_{2} & =y_{2}ee'\\
ee_{3}-ee_{2} & =y_{1}ee''\\
\tau y_{1}(ee_{3})-ee_{1} & =y_{1}ee''',\\
\chi & =\_\otimes ee_{1}+y_{3}\chi_{1}\\
\chi_{1} & =\tau E(\_\otimes ee_{2})+y_{2}\chi_{1}'\\
\chi_{1}' & =E\tau\circ\tau E(\_\otimes ee_{3})+y_{1}\chi''\\
 & \text{ for some }ee^{(k)}\in E^{2}[y]\text{ and }\chi''\in\Hom_{A}(_{A}E,E^{3})[y]\biggr\rangle.
\end{split}
\end{gather*}

Here $\_.\theta\in\mathrm{End}_{A}(_{A}E)[y]$ is the morphism sending
$e$ to $e.\theta$, and $\_\otimes e_{1}\in\mathrm{Hom}_{A}(_{A}E,E^{2})[y]$
sends $e$ to $e\otimes e_{1}$, and $\_\otimes ee_{1}$ sends $e$
to $e\otimes ee_{1}$. Note that $e'$, $\xi_{1}$, and $\xi'$ are
uniquely determined by $(e_{1},e_{2},\xi)$, and $ee'$, $ee''$,
$ee'''$ and $\chi_{1}$, $\chi_{1}'$, $\chi''$ are uniquely determined
by $(ee_{1},ee_{2},ee_{3},\chi)$.
\end{prop}

We rely on this proposition in what follows and do not distinguish
between $G_{i}$ and $\bar{G}_{i}$ for $i=1,2,3$. For example, we
may write $(\theta,\varphi)\in G_{1}$. (We also write $(\theta,\varphi)$
for the element of $G_{1}^{\op}$.) The interpretations of elements
$(\theta,\varphi)$ etc.~as explicit homomorphisms of complexes are
given by Props.~3.18, 3.20, and 3.22 of \cite{mcmillanTensor2product2representations2022}.

Also recall the two complexes of $B$-modules: 
\begin{defn}
Let $R,X_{2}'\in B\text{-cplx}$ be given by: 
\[
R=\begin{pmatrix}E^{2}[y]\xrightarrow{\left(\begin{smallmatrix}\pi E\\
\pi E\ci\tau
\end{smallmatrix}\right)}E_{y}E\oplus E_{y}E\\
0\to E[y]\oplus E[y]
\end{pmatrix},
\]
 
\[
X_{2}'=\begin{pmatrix}\tau y_{1}E^{2}[y] & \overset{\pi E}{\lra} & E_{y}E\\
0 & \lra & E[y]
\end{pmatrix},
\]
 both lying in degrees $0$ and $1$ (cohomological grading), and
the $E_{y}$ action on $R$ specified by $0$ and the canonical map
\[
E_{y}\otimes(E[y]\oplus E[y])\to E_{y}E\oplus E_{y}E,
\]
 and on $X_{2}'$ specified by $0$ and $\idop_{E_{y}E[y]}$.
\end{defn}

Now recall the following three lemmas: 
\begin{lem}[Lem.~3.11 of \cite{mcmillanTensor2product2representations2022}]
We have that $X_{2}'$ is a finite direct sum of summands of $X$.
\end{lem}

The nil-affine Hecke algebra has the structure of an $n!\times n!$
matrix algebra over the symmetric polynomials $P_{n}^{S_{n}}$ (cf.~Prop.~3.4
of \cite{rouquier2KacMoodyAlgebras2008}). Among other things, this
gives a decomposition of $E^{n}$ into submodules called \textquoteleft divided
powers\textquoteright : 
\[
E^{n}\iso\overset{n\text{ copies}}{\overbrace{E^{(n)}\oplus\dots\oplus E^{(n)}}}.
\]
 We will make use of this for $n=2$, where the isomorphism is given
(by extension to left $A[y]$-modules) explicitly as follows: 
\begin{equation}
E^{2}[y]\xrightarrow[\left(\begin{smallmatrix}\tau y_{1}\\
\tau
\end{smallmatrix}\right)]{\sim}\tau y_{1}E^{2}[y]\oplus\tau y_{1}E^{2}[y].\label{eq:divided-squares-iso}
\end{equation}
 The inverse of this map is $(\iota,-y_{2})$, where $\iota:\tau y_{1}E^{2}[y]\hookrightarrow E^{2}[y]$
is the inclusion. The elements $\tau y_{1}$ and $-y_{2}\tau$ are
orthogonal idempotents summing to $\idop$, and $\tau$ gives an isomorphism
from $-y_{2}\tau E^{2}[y]$ to $\tau y_{1}E^{2}[y]$. (To check this:
$\tau(-y_{2}\tau)ee=\tau ee=\tau y_{1}(\tau ee)\in\tau y_{1}E^{2}[y]$,
$\tau(-y_{2}\tau)y_{1}ee=\tau y_{1}ee$, and if $\tau(-y_{2}\tau)ee=0$
then $\tau ee=0$ so $-y_{2}\tau ee=0$.)
\begin{lem}[Lem.~3.12 of \cite{mcmillanTensor2product2representations2022}]
 \label{lem:R-strictly-perfect} There is an isomorphism $R\iso X_{2}'\oplus X_{2}'$
in $B\text{-cplx}$ given by the above isomorphism on the degree $0$
term of the top row, and the identity on all other terms. So $R$
is a finite direct sum of summands of $X$. In particular, $R$ is
strictly perfect.
\end{lem}

\begin{lem}[Lem.~3.13 of \cite{mcmillanTensor2product2representations2022}]
\label{lem:R-quasi-iso-E'X} There is a quasi-isomorphism $R\xrightarrow{q.i.}E'X_{2}$
determined by $\idop_{E^{2}[y]}$ on the degree $0$ term of the top
row and $\left(\begin{smallmatrix}1 & 0\\
1 & -y_{1}
\end{smallmatrix}\right)$ on the degree $1$ term of the bottom row.
\end{lem}

We recall, finally, the main construction from \cite{mcmillanTensor2product2representations2022}.
That construction was given using an equivalence: 
\[
\text{per }B\xrightarrow[\mathscr{H}om_{B}(X,-)]{\sim}\text{per }C,
\]
 where $C=\End_{K^{b}(B)}(X_{1}\oplus X_{2})^{\op}$. The algebra
$C$ can be presented using a matrix of $(A[y],A[y])$-bimodules (see
Prop.~3.31 and Lem.~3.34 of \cite{mcmillanTensor2product2representations2022}):
\[
[C]\iso\begin{pmatrix}\mathrm{End}(X_{1})^{\op} & \mathrm{Hom}(X_{1},X_{2})\\
\mathrm{Hom}(X_{2},X_{1}) & \mathrm{End}(X_{2})^{\op}
\end{pmatrix}\iso\begin{pmatrix}A[y] & y_{1}E[y]\\
F[y] & G_{1}^{\op}
\end{pmatrix}.
\]
 The functor $E'\otimes-$ on $\text{per }B$ translates to the functor
$\mathscr{E}\otimes-$ on $\text{per }C$, with $\mathscr{E}=\mathscr{H}om_{B}(X,E'X)$,
and there is a quasi-isomorphism $\tilde{E}\xrightarrow{q.i.}\mathscr{E}$
with $\tilde{E}=\Hom_{K^{b}(B)}(X,E'X)$. This $\tilde{E}$ is a $(C,C)$-bimodule.
It can be presented as a matrix of $(A[y],A[y])$-bimodules (see §3.4.2
of \cite{mcmillanTensor2product2representations2022}): 
\begin{equation}
[\tilde{E}]\iso\begin{pmatrix}y_{1}E[y] & y_{1}y_{2}E^{2}[y]\\
G_{1} & G_{2}
\end{pmatrix}.\label{eq:Matrix-tilde-E}
\end{equation}
 Using the derived equivalence we also have an isomorphism $\tilde{E}^{2}=\tilde{E}\otimes_{C}\tilde{E}\iso\Hom_{K^{b}(B)}(X,E'^{2}X)$,
which yields a matrix presentation: 
\begin{equation}
[\tilde{E}^{2}]\iso\begin{pmatrix}y_{1}y_{2}E^{2}[y] & y_{1}y_{2}y_{3}E^{3}[y]\\
G_{2} & G_{3}
\end{pmatrix}.\label{eq:Matrix-tilde-E^2}
\end{equation}

Lastly, in §4 of \cite{mcmillanTensor2product2representations2022},
the author defined $(C,C)$-bimodule endomorphisms $\tilde{x}$ and
$\tilde{\tau}$. They are given componentwise by: 
\begin{equation}
[\tilde{x}]\acts[\tilde{E}]\quad\text{ by: }\quad\begin{pmatrix}x & xE\\
\substack{(\theta,\varphi)\mapsto\\
(y\theta,x\circ\varphi)
}
 & \substack{(e_{1},e_{2},\xi)\mapsto\\
(ye_{1},xe_{2},xE\ci\xi)
}
\end{pmatrix},\label{eq:Def-tilde-x}
\end{equation}
 
\begin{equation}
[\tilde{\tau}]\acts[\tilde{E}^{2}]\quad\text{ by: }\quad\begin{pmatrix}\tau & \tau E\\
\substack{(e_{1},e_{2},\xi)\mapsto\\
(e',e',\tau\circ\xi)
}
 & \substack{(ee_{1},ee_{2},ee_{3},\chi)\mapsto\\
(ee',ee',\tau(ee_{3}),\tau E\circ\chi)
}
\end{pmatrix}.\label{eq:Def-tilde-tau}
\end{equation}
 In the last row, $e'$ is determined by $e_{1}-e_{2}=y_{1}e'$, and
$ee'$ is determined by $ee_{1}-ee_{2}=y_{2}ee'$. (See Prop.~\ref{prop:G-bar_i}
above.) In \cite{mcmillanTensor2product2representations2022} it was
established that these endomorphisms satisfy the nil-affine Hecke
relations (\ref{eq:nil-aff-Hecke-rels}).

\section{More bimodules}

We add a new series of bimodules for this paper: 
\begin{defn}
Let $L_{n}$ denote $\Hom_{D^{b}(B)}(E'^{n}X_{1},X_{2})$.
\end{defn}

Note that $L_{1}=G_{1}$. We will only need $L_{1}$ and $L_{2}$
in what follows. Observe that $L_{n}$ has a right $G_{1}^{\op}$-module
structure given by post-composition. We now study $L_{2}$ and provide
it with the structure of a $(G_{1}^{\op},G_{1}^{\op})$-bimodule.

We need an additional feature of the complex $R$: 
\begin{lem}
\label{lem:G-action-R} The complex $R$ carries a right action of
the algebra $G_{1}^{\op}$, where $(\theta,\varphi)\in G_{1}^{\op}$
acts by post-composing with $E\varphi\in\End(E^{2}[y])$ on the top
row of $R_{0}$, namely $E^{2}[y]$, and by the matrix 
\[
\Phi=\begin{pmatrix}\varphi & 0\\
\varphi_{1} & \theta
\end{pmatrix}
\]
 on the bottom row of $R_{1}$, namely $E[y]^{\oplus2}$, and by $E_{y}\Phi$
on the top row of $R_{1}$, namely $E_{y}E^{\oplus2}$. Through the
quasi-isomorphism of Lemma \ref{lem:R-quasi-iso-E'X}, this action
induces the canonical action of $G_{1}^{\op}=\End_{K^{b}(B)}(X_{2})^{\op}$
on $E'X_{2}$ given by functoriality of $E'$.
\end{lem}

\begin{proof}
First we check that the right action of $(\theta,\varphi)$ described
in the lemma gives a morphism of complexes of left $B$-modules. The
action is clearly $A[y]$-linear in the top and bottom rows, and it
is clearly linear over the off-diagonal generators in $E_{y}\subset B$.
The action commutes with the differential on the bottom row. We check
the top row: 
\begin{align*}
\begin{pmatrix}E_{y}\varphi & 0\\
E_{y}\varphi_{1} & E_{y}\theta
\end{pmatrix}\cdot\begin{pmatrix}\pi E\\
\pi E\circ\tau
\end{pmatrix} & =\begin{pmatrix}E_{y}\varphi\circ\pi E\\
E_{y}\varphi_{1}\circ\pi E+E_{y}\theta\circ\pi E\circ\tau
\end{pmatrix}\\
 & =\begin{pmatrix}\pi E\circ E\varphi\\
\pi E\circ E\varphi_{1}+\pi E\circ\tau\circ E\theta
\end{pmatrix}\\
 & =\begin{pmatrix}\pi E\\
\pi E\circ\tau
\end{pmatrix}\circ E\Phi.
\end{align*}

Next we check that the action commutes with multiplication in the
algebra. In $G_{1}^{\op}$ we have $(\theta,\varphi)\cdot(\theta',\varphi')=(\theta\theta',\varphi'\ci\varphi)$,
while the action of the product on the bottom row of $R_{1}$ is given
by: 
\[
\begin{pmatrix}\varphi' & 0\\
\varphi_{1}' & \theta'
\end{pmatrix}\cdot\begin{pmatrix}\varphi & 0\\
\varphi_{1} & \theta
\end{pmatrix}=\begin{pmatrix}\varphi'\ci\varphi & 0\\
\varphi_{1}'\ci\varphi+(\_.\theta')\ci\varphi_{1} & \theta\theta'
\end{pmatrix}.
\]
 Note that 
\[
\varphi'\ci\varphi-\_.\theta\theta'=y_{1}\bigl((\_.\theta')\ci\varphi_{1}+\varphi_{1}'\ci\varphi\bigr),
\]
 so the composition of the actions agrees with the action of the product
on that term. The other terms are trivial to check.

Lastly we check that through the quasi-isomorphism of Lemma \ref{lem:R-quasi-iso-E'X},
this action is compatible with the canonical action on $E'X_{2}$.
Start with the bottom row of $R_{1}$: 
\begin{align*}
\begin{pmatrix}1 & 0\\
1 & -y_{1}
\end{pmatrix}\cdot\begin{pmatrix}\varphi & 0\\
\varphi_{1} & \theta
\end{pmatrix} & =\begin{pmatrix}\varphi & 0\\
\varphi-y_{1}\varphi_{1} & -y_{1}\theta
\end{pmatrix},\\
\begin{pmatrix}\varphi & 0\\
0 & \theta
\end{pmatrix}\cdot\begin{pmatrix}1 & 0\\
1 & -y_{1}
\end{pmatrix} & =\begin{pmatrix}\varphi & 0\\
\theta & -\theta y_{1}
\end{pmatrix}.
\end{align*}
 These agree because $\varphi-y_{1}\varphi_{1}=\theta$. The other
terms are trivial to check.
\end{proof}
Now we compute a model for $L_{2}$ using the strictly perfect $R$
as a replacement for $E'X_{2}$.
\begin{defn}
\label{def:L_2-bar} Define the following $(A[y],A[y])$-sub-bimodule
of $F[y]^{\oplus2}\oplus\Hom_{A}(_{A}E^{2},E)[y]$: 
\begin{gather*}
\bar{L}_{2}=\biggl\langle(f',f,\rho)\in F[y]^{\oplus2}\oplus\Hom_{A}(_{A}E^{2},E)[y]\biggr|\\
\begin{split}\rho & =Ef+Ef'\ci\tau+y_{1}\circ\rho'\\
 & \text{ for some }\rho'\in\Hom_{A}(_{A}E^{2},E)[y]\biggr\rangle.
\end{split}
\end{gather*}
 One easily checks that the set $\bar{L}_{2}$ is closed under the
bimodule operations.
\end{defn}

\begin{prop}
\label{prop:L_2-from-R} There is an isomorphism of $(A[y],A[y])$-bimodules
$\bar{L}_{2}\iso\Hom_{K^{b}(B)}(R,X_{2})$ determined by equivariance
over $E_{y}\subset B$ with the following data: 
\[
(f',f,\rho)\mapsto\left(\begin{pmatrix}(ee,\left(\begin{smallmatrix}0\\
0
\end{smallmatrix}\right))\\
(0,\left(\begin{smallmatrix}e\\
e'
\end{smallmatrix}\right))
\end{pmatrix}\mapsto\begin{pmatrix}(\rho(ee),0)\\
(0,f(e)+f'(e'))
\end{pmatrix}\right).
\]
\end{prop}

\begin{proof}
The proof is seen by directly computing $Z^{0}\mathscr{H}om_{B}(R,X_{2})$.
It is easy to check that the morphism given as the image of $(f',f,\rho)$
is a morphism of complexes of left $B$-modules. The condition $\rho=Ef+Ef'\ci\tau+y_{1}\circ\rho'$
is equivalent to the statement that this morphism has zero differential.
\end{proof}
(Recall the notation from \cite{mcmillanTensor2product2representations2022}:
$ee$ is an arbitrary element of $E^{2}[y]$, not a simple tensor.
It is unrelated to $e$ and $e'$, which are arbitrary in $E[y]$.)
\begin{cor}
\label{cor:L_2-bar=00003DL_2} The isomorphism above, followed by
the canonical isomorphism of functors $\Hom_{K^{b}(B)}(R,-)\iso\Hom_{D^{b}(B)}(R,-)$
applied to $X_{2}$, gives an isomorphism $\bar{L}_{2}\iso L_{2}$
of $(A[y],A[y])$-bimodules.
\end{cor}

\begin{prop}
\raggedright \label{prop:F2-from-R} There is an isomorphism of $(A[y],A[y])$-bimodules
$F^{2}[y]\iso\Hom_{K^{b}(B)}(R,X_{1})$ determined by equivariance
over $E_{y}\subset B$ with the following data: 
\[
F^{2}[y]\ni ff\mapsto\left(\begin{pmatrix}(ee,\left(\begin{smallmatrix}0\\
0
\end{smallmatrix}\right))\\
(0,\left(\begin{smallmatrix}e\\
e'
\end{smallmatrix}\right))
\end{pmatrix}\mapsto\begin{pmatrix}ff(ee)\\
0
\end{pmatrix}\right).
\]
\end{prop}

\begin{proof}
The proof is seen by directly computing $Z^{0}\mathscr{H}om_{B}(R,X_{1})$.
\end{proof}
It is useful to give a model of $G_{2}$ that is compatible with this
model of $L_{2}$ by using the replacement $R$ for $E'X_{2}$.
\begin{defn}
\label{def:G_2'-bar} Define the following $(A[y],A[y])$-sub-bimodule
of $E[y]^{\oplus2}\oplus\Hom_{A}(_{A}E,E^{2})[y]$: 
\begin{gather*}
\bar{G}_{2}'=\biggl\langle(e',e,\xi)\in E[y]^{\oplus2}\oplus\Hom_{A}(_{A}E,E^{2})[y]\biggr|\\
\begin{split}\xi & =\_\otimes e+y_{2}\tau\left(\_\otimes(e-y_{1}e')\right)+y_{1}y_{2}\xi'\\
 & \text{ for some }\xi'\in\Hom_{A}(_{A}E,E^{2})[y]\biggr\rangle.
\end{split}
\end{gather*}
 One quickly checks that the condition is closed under the bimodule
operations. It is sometimes convenient to rewrite the condition as
\[
\xi=\tau y_{1}(\_\otimes e)-y_{2}\tau y_{1}(\_\otimes e')+y_{1}y_{2}\xi'.
\]
\end{defn}

\begin{prop}
\label{prop:G_2'-homs} There is an isomorphism of $(A[y],A[y])$-bimodules
$\bar{G}_{2}'\iso\Hom_{K^{b}(B)}(X_{2},R)$ determined by equivariance
over $E_{y}\subset B$ with the following data: 
\[
(e',e,\xi)\mapsto\left(\begin{pmatrix}(e,0)\\
(0,1)
\end{pmatrix}\mapsto\begin{pmatrix}(\xi(e),0)\\
(0,\left(\begin{smallmatrix}e\\
e'
\end{smallmatrix}\right))
\end{pmatrix}\right).
\]
\end{prop}

\begin{proof}
The proof is seen by directly computing $Z^{0}\mathscr{H}om_{B}(X_{2},R)$.
\end{proof}
The quasi-isomorphism $R\xrightarrow{q.i.}E'X_{2}$ determines an
isomorphism $\bar{G}_{2}'\iso\bar{G}_{2}$, since $X_{2}$ is strictly
perfect, given by $(e',e,\xi)\mapsto(e,e-y_{1}e',\xi)$, with inverse
given by $(e_{1},e_{2},\xi)\mapsto(y_{1}^{-1}(e_{1}-e_{2}),e_{1},\xi)$.
In most of this paper we will use $\bar{G}_{2}'$ instead of $\bar{G}_{2}$
as a model for $G_{2}$.
\begin{defn}
\raggedright Let $U$ denote $\Hom_{K^{b}(B)}(R,R)$. It is canonically
isomorphic to $\Hom_{D^{b}(B)}(E'X_{2},E'X_{2})$.
\end{defn}

Now we describe a model for $U$. For $U$ and $L_{2}$ later in this
paper, as for $G_{n}$, we frequently assume the terms of the models
to denote morphisms of complexes, passing without mention through
the isomorphisms $\bar{U}\iso U$ and $\bar{L}_{2}\iso L_{2}$.
\begin{defn}
\raggedright \label{def:U-bar} Define the following $(A[y],A[y])$-sub-bimodule
of $FE[y]^{\oplus4}\oplus\Hom_{A}(_{A}E^{2},E^{2})[y]$: 

\begin{gather*}
\bar{U}=\biggl\langle(\Phi_{11},\Phi_{21},\Phi_{12},\Phi_{22},\Lambda)\in FE[y]^{\oplus4}\oplus\Hom_{A}(_{A}E^{2},E^{2})[y]\biggr|\\
\begin{split}\Lambda & =\tau y_{1}(E\Phi_{11}+E\Phi_{12}\circ\tau)-y_{2}\tau y_{1}(E\Phi_{21}+E\Phi_{22}\circ\tau)+y_{1}y_{2}\Lambda^{\circ}\\
 & \text{ for some }\Lambda^{\circ}\in\Hom_{A}(_{A}E^{2},E^{2})[y]\biggr\rangle.
\end{split}
\end{gather*}
 Here $\Phi_{ij}$ give the components of the matrix $[\Phi]$ of
a map $\Phi\in\End_{A}(_{A}E[y]\oplus E[y])$. Note that because $y_{1}y_{2}$
is injective, $\Lambda^{\ci}$ is uniquely determined by $(\Phi,\Lambda)$.
The condition on $\Lambda$ is clearly closed under the bimodule operations.
\end{defn}

\begin{prop}
There is an isomorphism of $(A[y],A[y])$-bimodules $\bar{U}\iso U$
determined by equivariance over $E_{y}\subset B$ with the following
data: 
\[
(\Phi,\Lambda)\mapsto\left(\begin{pmatrix}\bigl(ee,\left(\begin{smallmatrix}0\\
0
\end{smallmatrix}\right)\bigr)\\
\bigl(0,\left(\begin{smallmatrix}e\\
e'
\end{smallmatrix}\right)\bigr)
\end{pmatrix}\mapsto\begin{pmatrix}\bigl(\Lambda(ee),\left(\begin{smallmatrix}0\\
0
\end{smallmatrix}\right)\bigr)\\
\bigl(0,[\Phi]\cdot\left(\begin{smallmatrix}e\\
e'
\end{smallmatrix}\right)\bigr)
\end{pmatrix}\right).
\]
\end{prop}

\begin{proof}
The proof is seen by directly computing $U=Z^{0}\mathscr{H}om_{B}(R,R)$.
We must show that the condition on $\Lambda$ is equivalent to the
statement that the image of $(\Phi,\Lambda)$ has zero differential.
One computes directly that the morphism given as this image has zero
differential if and only if the following pair of equations holds:
\[
\left\{ \begin{aligned}\pi E\ci\Lambda & =E_{y}\Phi_{11}\ci\pi E+E_{y}\Phi_{12}\ci\pi E\ci\tau\\
\pi E\ci\tau\Lambda & =E_{y}\Phi_{21}\ci\pi E+E_{y}\Phi_{22}\ci\pi E\ci\tau.
\end{aligned}
\right.
\]
 These are morphisms from $E^{2}[y]$ in the top row of $R_{0}$.
On the left side they are given by applying the image of $(\Phi,\Lambda)$
first (namely $\Lambda$ on $E^{2}[y]$) and then $d$. On the right
side, $E_{y}\Phi$ is induced on the top row of $R_{1}$ by $\Phi$
on the bottom row of $R_{1}$ together with equivariance over $E_{y}\subset B$.

That pair of equations is equivalent to the condition: 
\begin{gather}
\left\{ \begin{aligned}\Lambda & =E\Phi_{11}+E\Phi_{12}\ci\tau+y_{2}\Lambda'\\
\tau\Lambda & =E\Phi_{21}+E\Phi_{22}\ci\tau+y_{2}\Lambda''
\end{aligned}
\right.\label{eq:Lambda-obvious-condition}\\
\text{ for some }\Lambda',\Lambda''\in\Hom_{A}(_{A}E^{2},E^{2})[y].\nonumber 
\end{gather}
 For example, the first equation of the pair is equivalent to $\pi E\ci(\Lambda-E\Phi_{11}-E\Phi_{12}\ci\tau)=0$
because $\pi E$ commutes with $E_{y}\Phi_{ij}$. This identity implies
the first equation of (\ref{eq:Lambda-obvious-condition}) by Lemma
3.7 of \cite{mcmillanTensor2product2representations2022}; cf.~also
the proof of Prop.~3.26 in that paper.
\begin{claim}
Suppose $(\Phi,\Lambda)$ is given such that (\ref{eq:Lambda-obvious-condition})
holds for some $\Lambda'$, $\Lambda''$. Then there is $\Lambda^{\circ}\in\Hom_{A}(_{A}E^{2},E^{2})[y]$
such that 
\begin{equation}
\Lambda=\tau y_{1}(E\Phi_{11}+E\Phi_{12}\circ\tau)-y_{2}\tau y_{1}(E\Phi_{21}+E\Phi_{22}\circ\tau)+y_{1}y_{2}\Lambda^{\circ}.\label{eq:Lambda-condition}
\end{equation}
\end{claim}

\begin{proof}
Multiply the second equation of (\ref{eq:Lambda-obvious-condition})
by $\tau$ and obtain: 
\[
-\tau y_{2}\Lambda''=\tau\ci E\Phi_{21}+\tau\ci E\Phi_{22}\ci\tau.
\]
 Multiply the first by $\tau$ and the second by $\tau y_{1}$ and
identify the results to obtain: 
\[
\tau y_{2}\Lambda'=y_{1}y_{2}\tau\Lambda''+\tau y_{1}\ci\bigl(E\Phi_{21}+E\Phi_{22}\ci\tau\bigr)-\tau\ci\bigl(E\Phi_{11}+E\Phi_{12}\ci\tau\bigr).
\]
 Then: 
\begin{align*}
\Lambda' & =(y_{1}\tau-\tau y_{2})\ci\Lambda'\\
 & =y_{1}\tau\Lambda'-y_{1}y_{2}\tau\Lambda''-\tau y_{1}\ci\bigl(E\Phi_{21}+E\Phi_{22}\ci\tau\bigr)+\tau\ci\bigl(E\Phi_{11}+E\Phi_{12}\ci\tau\bigr)\\
 & =y_{1}\bigl(\tau\Lambda'-y_{2}\tau\Lambda''\bigr)-\tau y_{1}\ci\bigl(E\Phi_{21}+E\Phi_{22}\ci\tau\bigr)+\tau\ci\bigl(E\Phi_{11}+E\Phi_{12}\ci\tau\bigr).
\end{align*}
 Let $\Lambda^{\ci}=\tau\Lambda'-y_{2}\tau\Lambda''$. Then: 
\begin{align*}
\Lambda & =E\Phi_{11}+E\Phi_{12}\ci\tau+y_{1}y_{2}\Lambda^{\ci}\\
 & \qquad-y_{2}\tau y_{1}\ci\bigl(E\Phi_{21}+E\Phi_{22}\ci\tau\bigr)+y_{2}\tau\ci\bigl(E\Phi_{11}+E\Phi_{12}\ci\tau\bigr)\\
 & =\tau y_{1}\ci(E\Phi_{11}+E\Phi_{12}\circ\tau)-y_{2}\tau y_{1}\ci(E\Phi_{21}+E\Phi_{22}\circ\tau)+y_{1}y_{2}\Lambda^{\circ},
\end{align*}
 as desired.
\end{proof}
\begin{claim}
Now suppose $(\Phi,\Lambda)$ and $\Lambda^{\circ}$ are given such
that (\ref{eq:Lambda-condition}) holds. Then there are $\Lambda'$,
$\Lambda''$ such that (\ref{eq:Lambda-obvious-condition}) holds.
\end{claim}

\begin{proof}
Let 
\begin{align*}
\Lambda' & =\tau\ci\bigl(E\Phi_{11}+E\Phi_{12}\ci\tau\bigr)-\tau y_{1}\ci\bigl(E\Phi_{21}+E\Phi_{22}\ci\tau\bigr)+y_{1}\Lambda^{\ci},\\
\Lambda'' & =\tau\ci\bigl(E\Phi_{21}+E\Phi_{22}\ci\tau\bigr)+y_{1}\tau\Lambda^{\ci}.
\end{align*}
 Multiplying the first by $y_{2}$, adding $E\Phi_{11}+E\Phi_{12}\ci\tau$,
and simplifying with (\ref{eq:Lambda-condition}), we find: 
\[
y_{2}\Lambda'+E\Phi_{11}+E\Phi_{12}\ci\tau=\Lambda.
\]
 Multiplying the second by $y_{2}$ and adding $E\Phi_{21}+E\Phi_{22}\ci\tau$,
we find: 
\begin{align*}
y_{2}\Lambda''+E\Phi_{21}+E\Phi_{22}\ci\tau & =\tau y_{1}\ci\bigl(E\Phi_{21}+E\Phi_{22}\ci\tau\bigr)+\tau y_{1}y_{2}\Lambda^{\ci},
\end{align*}
 while 
\begin{align*}
\tau\Lambda & =-\tau y_{2}\tau y_{1}\ci\bigl(E\Phi_{21}+E\Phi_{22}\ci\tau\bigr)+\tau y_{1}y_{2}\Lambda^{\ci}\\
 & =\tau y_{1}\ci\bigl(E\Phi_{21}+E\Phi_{22}\ci\tau\bigr)+y_{1}y_{2}\tau\Lambda^{\ci}
\end{align*}
 using (\ref{eq:Lambda-condition}). So the pair of equations (\ref{eq:Lambda-obvious-condition})
is satisfied.
\end{proof}
The proposition follows.
\end{proof}
We will need one more description of $U$:
\begin{lem}
\label{lem:U-composition} The composition map $L_{2}\otimes_{G_{1}^{\op}}G_{2}\to U$
is an isomorphism.
\end{lem}

\begin{proof}
Consider the triangulated functor: 
\[
\mathscr{H}om_{B}(X_{2},-):K^{b}(B)\to K^{b}(G_{1}^{\op}).
\]
 By the same reasoning as in §3.3.2 of \cite{mcmillanTensor2product2representations2022},
this functor descends to the derived categories 
\[
\mathscr{H}om_{B}(X_{2},-):D^{b}(B)\to D^{b}(G_{1}^{\op}),
\]
 it is fully faithful when restricted to $\langle X_{2}\rangle_{\Delta}$,
and it is essentially surjective from $\langle X_{2}\rangle_{\Delta}$
(because the image of $X_{2}$ is quasi-isomorphic to $G_{1}^{\op}$).
The inverse is given by $X_{2}\otimes_{G_{1}^{\op}}-$. It follows
from $R\in\langle X_{2}\rangle_{\Delta}$ (Lemma 3.12 of \cite{mcmillanTensor2product2representations2022})
and 
\begin{align*}
\Hom_{K^{b}(B)}(X_{2},R) & \iso\Hom_{K^{b}(B)}(X_{2},E'X_{2})\\
 & \xrightarrow{q.i.}\mathscr{H}om_{B}(X_{2},E'X_{2})\\
 & \xrightarrow{q.i.}\mathscr{H}om_{B}(X_{2},R)
\end{align*}
 that the evaluation map is an isomorphism: 
\[
X_{2}\otimes_{G_{1}^{\op}}\Hom_{K^{b}(B)}(X_{2},R)\iso R.
\]
 This shows that the map in the lemma statement is an isomorphism: 

\begin{align*}
 & \Hom_{K^{b}(B)}(R,X_{2})\otimes_{G_{1}^{\op}}\Hom_{K^{b}(B)}(X_{2},R)\\
 & \iso\Hom_{K^{b}(B)}\bigl(R,X_{2}\otimes_{G_{1}^{\op}}\Hom_{K^{b}(B)}(X_{2},R)\bigr)\\
 & \iso\Hom_{K^{b}(B)}(R,R).
\end{align*}
\end{proof}
We will need to know the $(A[y],A[y])$-bimodule structure of the
components of $\tilde{E}$ and $\tilde{E}^{2}$ and $\tilde{F}$.
These may be read off of presentations we have given by using the
fact that $y_{i}=x_{i}-y$ is injective as an endomorphism of $E^{n}[y]$
(for any $n$). We write $y_{i}^{-1}$ for the inverse morphism defined
on the image $y_{i}E^{n}[y]$.
\begin{prop}
\label{prop:various-as-A-A-bimod} We have isomorphisms of $(A[y],A[y])$-bimodules: 
\begin{itemize}
\item $y_{1}\dots y_{n}E^{n}[y]\iso E^{n}[y]$ given by application of $(y_{1}\dots y_{n})^{-1}$.
\item $L_{1}=G_{1}\iso A[y]\oplus FE[y]$ given by $(\theta,\varphi)\mapsto(\theta,\varphi_{1})$,
where 
\[
\varphi_{1}=y_{1}^{-1}(\varphi-\theta)
\]
 is interpreted in $FE[y]$. Note that the summand $FE^{2}[y]$ is
not only a left $A[y]$-submodule of $G_{2}$, but moreover a left
$G_{1}^{\op}$-submodule of $G_{2}$.
\item $G_{2}\iso E[y]\oplus E[y]\oplus FE^{2}[y]$ given by $(e',e,\xi)\mapsto(e',e,\xi')$,
where 
\[
\xi'=(y_{1}y_{2})^{-1}\bigl(\xi-\_\otimes e-y_{2}\tau(\_\otimes(e-y_{1}e'))\bigr)
\]
 is interpreted in $FE^{2}[y]$. Note that the summand $FE^{2}[y]$
is a left $G_{1}^{\op}$-submodule of $G_{2}$.
\item $L_{2}\iso F[y]\oplus F[y]\oplus F^{2}E[y]$ given by $(f',f,\rho)\mapsto\left(f',f,\rho_{1}\right)$,
where 
\[
\rho_{1}=y_{1}^{-1}\bigl(\rho-Ef-Ef'\circ\tau\bigr)
\]
 is interpreted in $F^{2}E[y]$. Note that the summand $F^{2}E[y]$
is a left $G_{1}^{\op}$-submodule of $L_{2}$.
\item $U\iso FE[y]^{\oplus4}\oplus F^{2}E^{2}[y]$ given by 
\[
(\Phi_{11},\Phi_{21},\Phi_{12},\Phi_{22},\Lambda)\mapsto(\Phi_{11},\Phi_{21},\Phi_{12},\Phi_{22},\Lambda^{\circ}),
\]
 where 
\begin{equation}
\Lambda=\tau y_{1}(E\Phi_{11}+E\Phi_{12}\circ\tau)-y_{2}\tau y_{1}(E\Phi_{21}+E\Phi_{22}\circ\tau)+y_{1}y_{2}\Lambda^{\circ}\label{eq:Lambda^circ}
\end{equation}
determines $\Lambda^{\circ}$, which is interpreted in $F^{2}E^{2}[y]$.
Note that the summand $F^{2}E^{2}[y]$ is a left $G_{1}^{\op}$-submodule
of $U$.
\end{itemize}
\end{prop}

\begin{proof}
The first point is obvious. The second point follows from Prop.~\ref{prop:G-bar_i}
because $(\theta,\varphi_{1})$ may be chosen arbitrarily in $A^{\op}[y]\oplus\mathrm{End}_{A}(_{A}E)[y]$,
determining $\varphi$; while a given choice of $(\theta,\varphi)$
satisfying the condition determines $\varphi_{1}$ by way of $y_{1}^{-1}$.
Similar reasoning applies to $G_{2}$, $L_{2}$, and $U$, working
with Defs.~\ref{def:G_2'-bar}, \ref{def:L_2-bar}, and \ref{def:U-bar},
respectively.
\end{proof}
In what follows, we will frequently use the bimodule descriptions
on the right sides of the isomorphisms in Prop.~\ref{prop:various-as-A-A-bimod}.
Sometimes, to avoid confusion, we will use the shorthand expression
\textquoteleft submodule form\textquoteright{} to refer to the left
sides of the isomorphisms (i.e.~presentations as submodules cut out
by conditions, as in the definitions of these structures), and \textquoteleft bimodule
form\textquoteright{} to refer to the right sides of the isomorphisms.
Considering the component data of an element in one of these structures,
the components in submodule form and bimodule form differ only in
the last component: in submodule form the last component gives the
full morphism on the degree $0$ part of the top row of the complex,
and in bimodule form the last component gives the remainder term \textquoteleft $\varphi_{1}$\textquoteright ,
\textquoteleft $\xi'$\textquoteright , \textquoteleft $\chi''$\textquoteright ,
\textquoteleft $\rho_{1}$\textquoteright , or \textquoteleft $\Lambda^{\circ}$\textquoteright{}
that is produced from the conditions by inverting some $y_{i}$.

\section{Adjunction}
\begin{defn}
Let $\tilde{F}$ denote the $(C,C)$-bimodule $\phantom{}^{\vee\negmedspace}\tilde{E}$,
that is, $\Hom_{C}(_{C}\tilde{E},C)$.
\end{defn}

We know that $_{C}\tilde{E}$ is f.g.~projective. It follows that
the right adjoint functor $\Hom_{C}(_{C}\tilde{E},-)$ of $\tilde{E}\otimes_{C}-$
is canonically isomorphic to $\tilde{F}\otimes_{C}-$. We have already
defined $\tilde{x}$ and $\tilde{\tau}$. We define $\tilde{\varepsilon}:\tilde{E}\tilde{F}\to C$
and $\tilde{\eta}:C\to\tilde{F}\tilde{E}$ using the duality, and
then $\tilde{\sigma}$ and $\tilde{\rho}_{\lambda}$ using the formulas
in §\ref{subsec:Commutator-defs} with $(A,E,F,x,\tau,\eta,\veps)$
replaced by $(C,\tilde{E},\tilde{F},\tilde{x},\tilde{\tau},\tilde{\eta},\tilde{\veps})$.
Note that sometimes we view $\tilde{F}\tilde{E}$ through the canonical
isomorphism $\Hom(\tilde{E},C)\otimes_{C}\tilde{E}\iso\Hom(\tilde{E},\tilde{E})$.

Now we construct an isomorphism of $(C,C)$-bimodules 
\[
\tilde{F}\iso\Hom_{K^{b}(B)}(X_{2}\oplus R,X)
\]
 as follows: 
\begin{align*}
\tilde{F} & =\;\Hom_{C}(_{C}\tilde{E},C)\\
 & \iso\Hom_{D^{b}(C)}(\tilde{E},C)\\
 & \iso\;\Hom_{D^{b}(C)}\left(\mathscr{H}om_{B}(X,E'X),\mathscr{H}om_{B}(X,X)\right)\\
 & \iso\Hom_{D^{b}(B)}(E'X,X)\\
 & \iso\Hom_{D^{b}(B)}(X_{2}\oplus R,X_{1}\oplus X_{2})\\
 & \iso\Hom_{K^{b}(B)}(X_{2}\oplus R,X_{1}\oplus X_{2}).
\end{align*}
The second arrow comes from the quasi-isomorphisms of Lemma 3.33 and
Corollary 3.41 of \cite{mcmillanTensor2product2representations2022}.
The third arrow comes from the equivalence $\text{per }B\iso\text{per }C$.
(By a similar calculation we have $\tilde{F}^{n}\iso\mathrm{Hom}_{D^{b}(B)}(E'^{n}X,X)$.
This explains the use of the derived category for $G_{n}$ and $L_{n}$.)
The fourth arrow holds because $R\overset{q.i.}{\to}E'X_{2}$ (Lemma
\ref{lem:R-quasi-iso-E'X}), and the fifth holds because $R$ is strictly
perfect (Lemma \ref{lem:R-strictly-perfect}).

With this description of $\tilde{F}$, we can express it as a $2\times2$
matrix of $(A[y],A[y])$-bimodules whose entries we have studied:
\begin{equation}
[\tilde{F}]\iso\begin{pmatrix}\mathrm{Hom}(X_{2},X_{1}) & \mathrm{Hom}(X_{2},X_{2})\\
\mathrm{Hom}(R,X_{1}) & \mathrm{Hom}(R,X_{2})
\end{pmatrix}\iso\begin{pmatrix}F[y] & L_{1}\\
F^{2}[y] & L_{2}
\end{pmatrix}.\label{eq:Matrix-tilde-F}
\end{equation}

The top row of $[\tilde{F}]$ has been computed as the bottom row
of $[C]$. We found $\mathrm{Hom}_{K^{b}(B)}(R,X_{1})$ in Prop.~\ref{prop:F2-from-R},
and we found $\mathrm{Hom}_{K^{b}(B)}(R,X_{2})$ in Prop.~\ref{prop:L_2-from-R}.

We have $C=\End_{K^{b}(B)}(X_{1}\oplus X_{2})$, and the right action
of $C$ on $\tilde{F}$ is given by post-composition. The left action
of $C$ is by pre-composition, but one must first apply functoriality
of $E'$ and use the quasi-isomorphism from Lemma \ref{lem:R-quasi-iso-E'X},
which we write $\gamma:R\xrightarrow{q.i.}E'X_{2}$; we have: 
\begin{itemize}
\item A generator $\phi\in Z^{0}\mathscr{H}om_{B}(X_{1},X_{1})^{\op}\cong A[y]\subset C$
determines $E'\phi\in\Hom_{K^{b}(B)}(X_{2},X_{2})$ that acts on $\tilde{F}$
(on the top row) by pre-composition. An element $\phi=\theta\in A[y]$
acts in the obvious way on the left on $F[y]$ and $L_{1}$.
\item A generator $\phi\in Z^{0}\mathscr{H}om_{B}(X_{2},X_{1})\cong F[y]\subset C$
determines 
\[
E'\phi\in\Hom_{D^{b}(B)}(E'X_{2},E'X_{1})\xrightarrow[\sim]{-\ci\gamma}\Hom_{K^{b}(B)}(R,X_{2}).
\]
 So $\phi$ acts on $\tilde{F}$ by pre-composition with $E'\phi\ci\gamma:R\to X_{2}$,
taking the top row to the bottom row. Recall that we have the model
$\bar{L}_{2}$ for $\Hom_{K^{b}(B)}(R,X_{2})$. An element $\phi=f\in F[y]$
acts by pre-composition with the morphism determined by $E'\phi\ci\gamma=(0,f,0)\in\bar{L}_{2}$.
\item A generator $\phi\in Z^{0}\mathscr{H}om_{B}(X_{1},X_{2})\cong y_{1}E[y]\subset C$
determines 
\[
E'\phi\in\Hom_{K^{b}(B)}(E'X_{1},E'X_{2})\xleftarrow[\sim]{\gamma\ci-}\Hom_{K^{b}(B)}(X_{2},R).
\]
 Recall that we have the models $\bar{G}_{2}$ for $\Hom_{K^{b}(B)}(X_{2},E'X_{2})$
and $\bar{G}_{2}'$ for $\Hom_{K^{b}(B)}(X_{2},R)$, and the isomorphism
$\bar{G}_{2}\iso\bar{G}_{2}'$ given by $(e_{1},e_{2},\xi')\mapsto(y_{1}^{-1}(e_{1}-e_{2}),e_{1},\xi')$
(in bimodule forms). An element $\phi=y_{1}e\in y_{1}E[y]$ determines
$E'\phi=(y_{1}e,0,0)\in\bar{G}_{2}$, so this acts on $\tilde{F}$
by pre-composition with the morphism determined by $(e,y_{1}e,0)\in\bar{G}_{2}'$,
taking the bottom row to the top row.
\item A generator 
\[
\phi\in Z^{0}\mathscr{H}om_{B}(X_{2},X_{2})^{\op}\cong G_{1}^{\op}\subset C
\]
 determines $\phi_{R}\in\Hom_{K^{b}(B)}(R,R)$ from the right action
of $G_{1}^{\op}$ on $R$. In terms of the model $\bar{U}$, we have
$\phi_{R}=(\varphi,\varphi_{1},0,\theta,E\varphi)$ (in submodule
form), determined by $\phi=(\theta,\varphi)\in G_{1}^{\op}$. This
acts on $\tilde{F}$ (on the bottom row) by pre-composition.
\end{itemize}

\section{Isomorphisms $\tilde{\rho}_{\lambda}$}

\subsection{Some tensor products of $(C,C)$-bimodules}

In this section we compute three tensor products of bimodules over
$C$, namely $\tilde{E}\tilde{E}$, $\tilde{F}\tilde{E}$, and $\tilde{E}\tilde{F}$,
and describe the products in each case as matrices of $(A[y],A[y])$-bimodules.
These calculations are used in the remaining sections to verify that
$\tilde{\rho}_{\lambda}$ are isomorphisms. Note that the product
$\tilde{E}\tilde{E}=\tilde{E}^{2}$ is already given a description
(Eq.~\ref{eq:Matrix-tilde-E^2}) as a matrix of $(A[y],A[y])$-bimodules
using the identification with $\Hom_{K^{b}(B)}(X,E'^{2}X)$, but in
order to compute $\tilde{\sigma}$ it is also necessary to realize
$\tilde{E}^{2}$ as the tensor product over $C$ of the bimodule $\tilde{E}$
with itself.

These tensor products are computed according to the general formulation
described in §2.4 of \cite{mcmillanTensor2product2representations2022}.
First we take the tensor product over the subalgebra $\Delta:=\left(\begin{smallmatrix}A[y] & 0\\
0 & G_{1}^{\op}
\end{smallmatrix}\right)\subset C$. This product is given on components by matrix multiplication and
tensor product over $A[y]$ or $G_{1}^{\op}$. After this we must
take the quotient by $\mathrm{Im}(I_{y_{1}E[y]})+\mathrm{Im}(I_{F[y]})$,
where $I_{y_{1}E[y]}$ and $I_{F[y]}$ apply the actions of the off-diagonal
generators in $C$. This quotient may be taken separately on each
coefficient of the product over $\Delta$.

Using the language of §2.4 of \cite{mcmillanTensor2product2representations2022}
with some given $M_{R}=\left(\begin{smallmatrix}M_{1} & M_{2}\end{smallmatrix}\right)$,
$_{R}N=\left(\begin{smallmatrix}N_{1}\\
N_{2}
\end{smallmatrix}\right)$, and $R=\left(\begin{smallmatrix}A & B\\
C & D
\end{smallmatrix}\right)$, the simplest technique at our disposal for computing a quotient
by the image of (say) $I_{B}$ is to identify one of its projections
as an isomorphism. (In §2.4 of \cite{mcmillanTensor2product2representations2022},
there is, for example, a projection of $I_{B}$ to $M_{1}\otimes_{A}N_{1}$
and another projection to $M_{2}\otimes_{D}N_{2}$.) In this situation
the quotient by $\mathrm{Im}(I_{B})$ reduces to the summand of the
second projection, because every element of the first summand (in
the quotient) has a unique representative in the second summand. If
it also happens that $\mathrm{Im}(I_{C})\subset\mathrm{Im}(I_{B})$,
then the quotient by the sum $\mathrm{Im}(I_{B})+\mathrm{Im}(I_{C})$
is still isomorphic to the second summand. Many of the components
computed below are found in this way, but a few of them require more
complicated reasoning.

Let us write, in general, $I_{\beta}'$ for the projection of $I_{B}$
to the first summand, and $-I_{\beta}''$ for the projection to the
second. Similarly write $I_{\delta}'$ and $-I_{\delta}''$ for the
projections of $I_{C}$. Here \textquoteleft first\textquoteright{}
and \textquoteleft second\textquoteright{} summand and \textquoteleft $I_{B}$\textquoteright{}
and \textquoteleft $I_{C}$\textquoteright{} are understood as in
§2.4 of \cite{mcmillanTensor2product2representations2022}. In a tensor
product of $(C,C)$-bimodules, each of the four coefficients will
have its own set of maps $I_{\beta}'$, $I_{\beta}''$, $I_{\delta}'$,
$I_{\delta}''$.

The matrix forms of the bimodules $\tilde{E}$, $\tilde{E}^{2}$,
and $\tilde{F}$ are given in Eqs.~(\ref{eq:Matrix-tilde-E}), (\ref{eq:Matrix-tilde-E^2}),
and (\ref{eq:Matrix-tilde-F}). For some of our calculations it helps
to be clear about the structures of the component bimodules, so we
translate the components to the bimodule forms on the right sides
of the isomorphisms in Prop.~\ref{prop:various-as-A-A-bimod}. Note
the consequence that the formulas for multiplication within $C$ and
for the actions of elements of $C$ on components of $\tilde{E}$
or $\tilde{F}$ are more complicated; this is illustrated, for example,
in the formula for $\Gamma_{12}\mid_{E^{2}[y]_{G}G_{2}}$ in the next
section.

\subsubsection{$\tilde{E}\tilde{E}$}

For the product $\tilde{E}\tilde{E}$, we already know the structures
of the coefficients of the matrix presentation from Eq.~(\ref{eq:Matrix-tilde-E^2})
(and Prop.~(\ref{prop:various-as-A-A-bimod})). We will need to compute
the action of $\tilde{\tau}$ on elements of $\tilde{E}\tilde{E}$
in order to compute $\tilde{\sigma}$, and for this it will be necessary
and sufficient to identify the map from the tensor product over $\Delta$
to the tensor product over $C$, i.e.~to the quotient by $\mathrm{Im}(I_{\beta}'-I_{\beta}'')+\mathrm{Im}(I_{\delta}'-I_{\delta}'')$.
Write $\Gamma$ for this map. Let the subscript \textquoteleft $G$\textquoteright{}
between concatenated modules indicate the tensor product over $G_{1}^{\op}$.
(An empty subscript indicates the product over $A[y]$.) So we have:
\begin{equation}
\begin{split}\tilde{E}\otimes_{\Delta}\tilde{E} & \cong\begin{pmatrix}E[y] & E^{2}[y]\\
G_{1} & G_{2}
\end{pmatrix}\otimes_{\Delta}\begin{pmatrix}E[y] & E^{2}[y]\\
G_{1} & G_{2}
\end{pmatrix}\\
 & \cong\begin{pmatrix}EE[y]\oplus E^{2}[y]_{G}G_{1} & EE^{2}[y]\oplus E^{2}[y]_{G}G_{2}\\
G_{1}E[y]\oplus(G_{2})_{G}G_{1} & G_{1}E^{2}[y]\oplus(G_{2})_{G}G_{2}
\end{pmatrix}\\
 & \xrightarrow{\Gamma}\begin{pmatrix}E^{2}[y] & E^{3}[y]\\
G_{2} & G_{3}
\end{pmatrix}\cong[\tilde{E}^{2}],
\end{split}
\label{eq:EE-expansion}
\end{equation}
 and we wish to understand the map $\Gamma$ on each component. We
find these component maps for simple tensors using the following steps.
First we interpret a pair of elements of the left and right tensor
factors $\tilde{E}$ as morphisms in $Z^{0}\mathscr{H}om_{B}(X,E'X)$
using isomorphisms such as Prop.~3.18 in \cite{mcmillanTensor2product2representations2022}.
Then we apply $E'$ to the morphism of the right factor, and post-compose
the result with the morphism of the left, obtaining an element of
$Z^{0}\mathscr{H}om_{B}(X,E'^{2}X)$. That element is interpreted
again in $\tilde{E}^{2}$. (See Prop.~3.37 and Lemmas 3.44 and 3.45
of \cite{mcmillanTensor2product2representations2022}.)

To facilitate checking these steps, the reader is encouraged to write
out the complexes for $X$, $E'X$, and $E'^{2}X$, and to be familiar
with Prop.~\ref{prop:various-as-A-A-bimod} and the interpretations
of elements of the structures in that proposition as homomorphisms
of complexes. With this in mind, the calculations are mechanical,
if tedious. We demonstrate the first cases with detailed explanation,
and for the remaining cases we record the results.
\begin{itemize}
\item For $\Gamma_{11}$, we have: 
\begin{itemize}
\item $\Gamma_{11}\mid_{EE[y]}$ is given by $\idop_{EE[y]}$.\\
\\
To see this, let $e_{1}$ represent an element of the left factor
$E[y]$, and $e_{2}$ an element of the right factor $E[y]$. (We
are suppressing the isomorphism $EE[y]\iso E[y]E[y]$.) Viewed in
$[\tilde{E}]_{11}$ through Prop.~\ref{prop:various-as-A-A-bimod},
these correspond to $y_{1}e_{1}$ and $y_{1}e_{2}$. As a homomorphism
of complexes, $y_{1}e_{2}$ sends $\left(\begin{smallmatrix}1\\
0
\end{smallmatrix}\right)\in X_{1}$ to $\left(\begin{smallmatrix}(y_{1}e_{2},0)\\
(0,0)
\end{smallmatrix}\right)\in E'X_{1}=X_{2}$. (Use Prop.~3.31, Lemma 3.34, and Lemma 3.7 of \cite{mcmillanTensor2product2representations2022}.)
Applying functoriality (with Lemma 3.8 of \cite{mcmillanTensor2product2representations2022}
in mind to notate $E'X_{2}$), we have $E'(y_{1}e_{2}):X_{2}\to E'X_{2}$
by a map that (among other things) takes the element $\left(\begin{smallmatrix}(y_{1}e_{1},0)\\
(0,0)
\end{smallmatrix}\right)$ to 
\[
\begin{pmatrix}(y_{1}e_{1}\otimes y_{1}e_{2},\left(\begin{smallmatrix}0\\
0
\end{smallmatrix}\right),0)\\
(0,\left(\begin{smallmatrix}0\\
0
\end{smallmatrix}\right),0)
\end{pmatrix}=\begin{pmatrix}(y_{1}y_{2}(e_{1}\otimes e_{2}),\left(\begin{smallmatrix}0\\
0
\end{smallmatrix}\right),0)\\
(0,\left(\begin{smallmatrix}0\\
0
\end{smallmatrix}\right),0)
\end{pmatrix},
\]
 which is therefore the image of $\left(\begin{smallmatrix}1\\
0
\end{smallmatrix}\right)\in X_{1}$ under $E'(y_{1}e_{2})\ci y_{1}e_{1}:X_{1}\to E'X_{2}$. By Lemma
3.47 of \cite{mcmillanTensor2product2representations2022} combined
with Prop.~\ref{prop:various-as-A-A-bimod}, this image corresponds
to $e_{1}\otimes e_{2}$ in $E^{2}[y]\cong[\tilde{E}^{2}]_{11}$,
so $\Gamma_{11}$ is the identity.\\
\item $\Gamma_{11}\mid_{E^{2}[y]_{G}G_{1}}$ is given as the inverse of
$E^{2}[y]\iso E^{2}[y]_{G}G_{1}$, $ee\mapsto ee\otimes1_{G_{1}}$.\\
\\
Here $ee$ corresponds to a map $X_{1}\to E'X_{2}$. The right action
of $g\in G_{1}$ on $E^{2}[y]$ is given in terms of maps of complexes
by post-composing with the induced map $E'(g):E'X_{2}\to E'X_{2}$,
but this is also the effect of $\Gamma_{11}$ on terms in $E^{2}[y]_{G}G_{1}$.
(The apparent coincidence derives from the definition $X_{2}=E'X_{1}$
and the matrix descriptions of $C$ and $\tilde{E}$.)\\
\end{itemize}
\item For $\Gamma_{21}$, we have:
\begin{itemize}
\item $\Gamma_{21}\mid_{G_{1}E[y]}$ is given (in bimodule forms) by 
\[
(\theta,\varphi_{1})\otimes e\mapsto(\theta e,\theta y_{1}e,\varphi_{1}(-)\otimes e)\in G_{2}.
\]
\\
The map $(\theta,\varphi_{1}):X_{2}\to X_{2}$ is determined (as in
Prop.~3.18 of \cite{mcmillanTensor2product2representations2022}
except in bimodule form) by: 
\[
\begin{pmatrix}(e_{1},0)\\
(0,1)
\end{pmatrix}\overset{(\theta,\varphi_{1})}{\longmapsto}\begin{pmatrix}(e_{1}.\theta+y_{1}\varphi_{1}(e_{1}),0)\\
(0,\theta)
\end{pmatrix}.
\]
Further, $e$ corresponds to the map $e:X_{1}\to X_{2}$ given by
$\left(\begin{smallmatrix}1\\
0
\end{smallmatrix}\right)\mapsto\left(\begin{smallmatrix}(y_{1}e,0)\\
(0,0)
\end{smallmatrix}\right)$, which by functoriality induces a map $E'(e):E'X_{1}\to E'X_{2}$
that is given (similarly to $\Gamma_{11}\mid_{EE[y]}$ above) by:
\[
\begin{pmatrix}(e_{1},0)\\
(0,1)
\end{pmatrix}\overset{E'(e)}{\longmapsto}\begin{pmatrix}(e_{1}\otimes y_{1}e,\left(\begin{smallmatrix}0\\
0
\end{smallmatrix}\right),0)\\
(0,\left(\begin{smallmatrix}y_{1}e\\
0
\end{smallmatrix}\right),0)
\end{pmatrix}.
\]
 Therefore the composition is given by: 
\[
\begin{pmatrix}(e_{1},0)\\
(0,1)
\end{pmatrix}\overset{E'(e)\ci(\theta,\varphi_{1})}{\longmapsto}\begin{pmatrix}(e_{1}.\theta\otimes y_{1}e+y_{1}\varphi_{1}(e_{1})\otimes y_{1}e,\left(\begin{smallmatrix}0\\
0
\end{smallmatrix}\right),0)\\
(0,\left(\begin{smallmatrix}\theta y_{1}e\\
0
\end{smallmatrix}\right),0)
\end{pmatrix}.
\]
 This image is in $E'X_{2}$. The map corresponds to the element 
\[
(\theta y_{1}e,0,\_\otimes\theta y_{1}e+y_{1}y_{2}(\varphi_{1}(-)\otimes e)\in\bar{G}_{2}
\]
(written in submodule form), which translates (see the paragraph after
Prop.~\ref{prop:G_2'-homs}) to the element 
\[
(\theta e,\theta y_{1}e,\varphi_{1}(-)\otimes e)\in\bar{G}_{2}'
\]
 (written in bimodule form) considering Def.~\ref{def:G_2'-bar}
and Prop.~\ref{prop:G_2'-homs} (and Prop.~\ref{prop:various-as-A-A-bimod}
for the bimodule form); this is the formula we wished to establish.\\
\item $\Gamma_{21}\mid_{(G_{2})_{G}G_{1}}$ is given as the inverse of $G_{2}\iso(G_{2})_{G}G_{1}$,
$g_{2}\mapsto g_{2}\otimes1_{G_{1}}$.\\
\\
This is similar to $\Gamma_{11}\mid_{E^{2}[y]_{G}G_{1}}$ above.\\
\end{itemize}
\item For $\Gamma_{12}$, we have:
\begin{itemize}
\item $\Gamma_{12}\mid_{EE^{2}[y]}$ is given by $\idop_{E^{3}[y]}$.\\
\\
We leave this and the remaining calculations to the reader, and record
the results.\\
\item $\Gamma_{12}\mid_{E^{2}[y]_{G}G_{2}}$ is given (in bimodule forms)
by 
\[
ee\otimes(e',e,\xi')\mapsto(y_{1}y_{2}y_{3})^{-1}(E\xi)(y_{1}y_{2}ee)\in E^{3}[y].
\]
\end{itemize}
\item For $\Gamma_{22}$, we have:
\begin{itemize}
\item $\Gamma_{22}\mid_{G_{1}E^{2}[y]}$ is given by 
\[
(\theta,\varphi_{1})\otimes ee\mapsto(\theta y_{1}y_{2}ee,0,0,\varphi_{1}(-)\otimes ee)\in G_{3}.
\]
\\
A remark: in §4.2 of \cite{mcmillanTensor2product2representations2022},
the horizontal arrows of the diagrams are based on the same type of
calculations, except that submodule forms are used. For example, the
map specified just above is $E'(ee)\ci(\theta,\varphi_{1})$. In submodule
form it would be written $E'(y_{1}y_{2}ee)\ci(\theta,\varphi)$ with
$\varphi=\_.\theta+y_{1}\varphi_{1}$, and this results from the horizontal
arrow in Diagram $D_{1\mid2}(2,1,1)$. So this composite would be
written in $G_{3}$ in submodule form as 
\[
E'(y_{1}y_{2}ee)\ci(\theta,\varphi)=(\theta y_{1}y_{2}ee,0,0,\varphi\otimes ee).
\]
 \\
We encourage the reader to check §4.2 of \cite{mcmillanTensor2product2representations2022}
for occasional further hints regarding some of these calculations.\\
\item $\Gamma_{22}\mid_{(G_{2})_{G}G_{2}}$ is given (in submodule forms)
by 
\[
(e_{1},e_{2},\xi)\otimes(\bar{e}_{1},\bar{e}_{2},\bar{\xi})\mapsto(\bar{\xi}(e_{1}),e_{2}\otimes\bar{e}_{1},e_{2}\otimes\bar{e}_{2},E\bar{\xi}\circ\xi)\in G_{3}
\]
 (c.f.~Diagram $D_{1\mid2}(2,2,1)$). We will need to have this map
written for the bimodule forms. First translate the notation from
$\bar{G}_{2}$ to $\bar{G}_{2}'$ (cf.~the paragraph after Prop.~\ref{prop:G_2'-homs})
using $e'$, $e$ for the first factor (so $e'=y_{1}^{-1}(e_{1}-e_{2})$
and $e=e_{1}$) and $\bar{e}'$, $\bar{e}$ for the second factor.
Then expand $\xi$ in terms of $e'$, $e$, and $\xi'$ according
to the condition for elements of $G_{2}$ (Def.~\ref{def:G_2'-bar}),
and likewise expand $\bar{\xi}$. \label{Gamma-22-EE-bimods} Now
compute $E\bar{\xi}\circ\xi$: 
\begin{align*}
E\bar{\xi}\circ\xi & =\bigl(\_\:\_\otimes\bar{e}+y_{2}\tau(\_\:\_\otimes(\bar{e}-y_{1}\bar{e}'))+y_{1}y_{2}E\bar{\xi}'\bigr)\\
 & \circ\bigl(\_\otimes e+y_{2}\tau(\_\otimes(e-y_{1}e'))+y_{1}y_{2}\xi'\bigr)\\
 & =\_\otimes\bigl(e\bar{e}+y_{2}\tau(e\bar{e}-y_{1}e\bar{e}')+y_{1}y_{2}\bar{\xi}'(e)\bigr)\\
 & +y_{3}\circ\tau E\bigl(\_\otimes(e\bar{e}-y_{2}e'\bar{e})\bigr)\\
 & +y_{2}y_{3}\circ E\tau\circ\tau E\bigl(\_\otimes(e-y_{1}e')(\bar{e}-y_{1}\bar{e}')\bigr)\\
 & +y_{1}y_{2}y_{3}\bigl(E\bar{\xi}'\circ\tau(\_\otimes(e-y_{1}e'))+E\tau(\xi'\otimes(\bar{e}-y_{1}\bar{e}'))+\xi'\otimes\bar{e}'\bigr).
\end{align*}
 To find $\chi''$, subtract all but the last term of $\chi$ in the
condition of $\bar{G}_{3}$ in Prop.~\ref{prop:G-bar_i} and remove
$y_{3}y_{2}y_{1}$. Obtain the complete image in bimodule form: 
\begin{align*}
 & (e',e,\xi')\otimes(\bar{e}',\bar{e},\bar{\xi}')\mapsto\\
 & \bigl(e\bar{e}+y_{2}\tau(e\bar{e}-y_{1}e\bar{e}')+y_{1}y_{2}\bar{\xi}'(e),e\bar{e}-y_{2}e'\bar{e},(e-y_{1}e')(\bar{e}-y_{1}\bar{e}'),\\
 & \qquad E\bar{\xi}'\circ\tau(\_\otimes(e-y_{1}e'))+E\tau(\xi'\otimes(\bar{e}-y_{1}\bar{e}'))+\xi'\otimes\bar{e}'\bigr).
\end{align*}
\end{itemize}
\end{itemize}

\subsubsection{$\tilde{F}\tilde{E}$}

For the product $\tilde{F}\tilde{E}$, we can find the $(A[y],A[y])$-bimodule
structure of the components of its matrix presentation using the same
technique as for $\tilde{F}$ and $\tilde{E}^{2}$. We have: 
\begin{align*}
\tilde{F}\tilde{E} & =\;\Hom_{C}(_{C}\tilde{E},C)\otimes_{C}\tilde{E}\\
 & \iso\Hom_{C}(_{C}\tilde{E},\tilde{E})\iso\Hom_{D^{b}(C)}(\tilde{E},\tilde{E})\iso\Hom_{D^{b}(C)}(\mathscr{E},\mathscr{E})\\
 & \iso\Hom_{D^{b}(C)}\left(\mathscr{H}om_{B}(X,E'X),\mathscr{H}om_{B}(X,E'X)\right)\\
 & =\;\Hom_{D^{b}(B)}(E'X,E'X)\\
 & \iso\Hom_{K^{b}(B)}(X_{2}\oplus R,X_{2}\oplus R).
\end{align*}
(The last isomorphism uses the quasi-isomorphism $R\overset{q.i.}{\to}E'X_{2}$
and the fact that $E'X_{1}=X_{2}$ and $R$ are strictly perfect.)
So the matrix presentation is: 
\begin{equation}
[\tilde{F}\tilde{E}]\iso\begin{pmatrix}\mathrm{Hom}(X_{2},X_{2}) & \mathrm{Hom}(X_{2},R)\\
\mathrm{Hom}(R,X_{2}) & \mathrm{Hom}(R,R)
\end{pmatrix}\iso\begin{pmatrix}G_{1} & G_{2}\\
L_{2} & U
\end{pmatrix}.\label{eq:Matrix-tilde-FE}
\end{equation}

As we did for $\tilde{E}^{2}$, we study the map $\Gamma$ from the
components of the product over $\Delta$ to those of the product over
$C$: 
\begin{equation}
\begin{split}\tilde{F}\otimes_{\Delta}\tilde{E} & \cong\begin{pmatrix}F[y] & L_{1}\\
F^{2}[y] & L_{2}
\end{pmatrix}\otimes_{\Delta}\begin{pmatrix}E[y] & E^{2}[y]\\
G_{1} & G_{2}
\end{pmatrix}\\
 & \cong\begin{pmatrix}FE[y]\oplus(L_{1})_{G}G_{1} & FE^{2}[y]\oplus(L_{1})_{G}G_{2}\\
F^{2}E[y]\oplus(L_{2})_{G}G_{1} & F^{2}E^{2}[y]\oplus(L_{2})_{G}G_{2}
\end{pmatrix}\\
 & \cong\begin{pmatrix}FE[y]\oplus G_{1} & FE^{2}[y]\oplus G_{2}\\
F^{2}E[y]\oplus L_{2} & F^{2}E^{2}[y]\oplus(L_{2})_{G}G_{2}
\end{pmatrix}\xrightarrow{\Gamma}\begin{pmatrix}G_{1} & G_{2}\\
L_{2} & U
\end{pmatrix}.
\end{split}
\label{eq:FE-expansion}
\end{equation}

The bulleted claims below are justified in the paragraphs following
them.
\begin{itemize}
\item We have $\Gamma_{11}:FE[y]\oplus G_{1}\to G_{1}$ given by $(\iota,\idop_{G_{1}})$.

Here the map $\iota:FE[y]\hookrightarrow L_{1}=G_{1}$ is the inclusion
of the second summand as written in Prop.~\ref{prop:various-as-A-A-bimod}.\\

\begin{itemize}
\item $I_{\beta}':FE[y]_{G}G_{1}\iso FE[y]$ given as the inverse of the
isomorphism $\bigl(fe\mapsto fe\otimes1_{G_{1}}\bigr)$,
\item $I_{\beta}'':FE[y]_{G}G_{1}\xrightarrow{\iota\otimes G_{1}}(L_{1})_{G}G_{1}\cong G_{1}$,
\item $I_{\delta}':(G_{1})_{G}FE[y]\iso FE[y]$ given as the inverse of
the isomorphism $\bigl(fe\mapsto1_{G_{1}}\otimes fe\bigr)$,
\item $I_{\delta}'':(G_{1})_{G}FE[y]\xrightarrow{G_{1}\otimes\iota}(G_{1})_{G}L_{1}\cong G_{1}$.\\
\end{itemize}
Using either $I_{\beta}''\circ I_{\beta}'^{-1}$, or $I_{\delta}''\circ I_{\delta}'^{-1}$,
one associates a unique representative in $(L_{1})_{G}G_{1}\cong G_{1}$
to each element of $FE[y]$. We see that $I_{\beta}''\circ I_{\beta}'^{-1}=I_{\delta}''\circ I_{\delta}'^{-1}$,
so the two associate the same representatives. It follows that the
quotient projection $\Gamma_{11}$ is given by $(\iota,\idop_{G_{1}})$
as proposed.\\

\item We have $\Gamma_{21}:F^{2}E[y]\oplus L_{2}\to L_{2}$ given by $(\iota',\idop_{L_{2}})$.

Here the map $\iota':F^{2}E[y]\hookrightarrow L_{2}$ is the inclusion
of the third summand as written in Prop.~\ref{prop:various-as-A-A-bimod}.\\

\begin{itemize}
\item $I_{\beta}':F^{2}E[y]_{G}G_{1}\iso F^{2}E[y]$ given as the inverse
of $\bigl(ffe\mapsto ffe\otimes1_{G_{1}}\bigr)$,
\item $I_{\beta}'':F^{2}E[y]_{G}G_{1}\xrightarrow{\iota'\otimes G_{1}}(L_{2})_{G}G_{1}\cong L_{2}$,
\item $I_{\delta}':(L_{2})_{G}FE[y]\to F^{2}E[y]$ given by 
\begin{align*}
(f',f,\rho')\otimes\bar{f}\bar{e} & \mapsto(\bar{f}\circ\rho)\otimes\bar{e}\\
 & \quad=\bigl(\bar{f}\circ(Ef+Ef'\circ\tau+y_{1}\rho')\bigr)\otimes\bar{e}.
\end{align*}
 Note that here $\bar{f}$ is interpreted as a map of complexes $\bar{f}:X_{2}\to X_{1}$
which is composed with the map $(f',f,\rho'):R\to X_{2}$ to obtain
a map $R\to X_{1}$. The composite $\bar{f}\ci\rho$ is a map $R\to X_{1}$,
interpreted again in $F^{2}[y]$ according to Prop.~\ref{prop:F2-from-R}.
(So, coincidentally, the notation \textquoteleft $\bar{f}\ci\rho$\textquoteright{}
has two valid interpretations: one as a map $R\to X_{1}$, and another
as a map $E^{2}[y]\to A[y]$ represented by an element of $F^{2}[y]$.)
\item $I_{\delta}'':(L_{2})_{G}FE[y]\xrightarrow{L_{2}\otimes\iota}(L_{2})_{G}G_{1}\cong L_{2}$.\\
\end{itemize}
Consider the first two maps. We have that $I_{\beta}''\circ I_{\beta}'^{-1}=\iota'$
as maps $F^{2}E[y]\to L_{2}$. Consider the last two maps. One may
check that $\iota'\circ I_{\delta}'=I_{\delta}''$. It follows that
$\mathrm{Im}(I_{\delta}'-I_{\delta}'')\subset\mathrm{Im}(I_{\beta}'-I_{\beta}'')$,
so in the quotient every element of $F^{2}E[y]$ is associated to
a unique element of $L_{2}$, given by applying the map $\iota'$.\\

\item We have $\Gamma_{12}:FE^{2}[y]\oplus G_{2}\to G_{2}$ given by $(\iota'',\idop_{G_{2}})$.

Here the map $\iota'':FE^{2}[y]\hookrightarrow G_{2}$ is the inclusion
of the third summand as written in Prop.~\ref{prop:various-as-A-A-bimod}.\\

\begin{itemize}
\item $I_{\beta}':FE[y]_{G}G_{2}\to FE^{2}[y]$ given by 
\begin{align*}
\bar{f}\bar{e}\otimes(e',e,\xi') & \mapsto\bar{f}\otimes(y_{1}y_{2})^{-1}\xi(y_{1}\bar{e})\\
 & \quad=\bar{f}\otimes\bigl(\tau(\bar{e}\otimes e)-y_{2}\tau(\bar{e}\otimes e')+\xi'(y_{1}\bar{e})\bigr).
\end{align*}
 The map is given by considering $\bar{e}$ as a map of complexes
$X_{1}\to X_{2}$, and $(e',e,\xi')$ as a map of complexes $X_{2}\to R$,
and then composing, and translating the result to bimodule form (removing
$y_{1}y_{2}$). The final expression is computed by plugging $y_{1}\bar{e}$
into 
\[
\xi=\tau y_{1}(\_\otimes e)-y_{2}\tau y_{1}(\_\otimes e')+y_{1}y_{2}\xi'
\]
 from the condition of Def.~\ref{def:G_2'-bar}, and we obtain: 
\begin{align*}
\xi(y_{1}\bar{e}) & =\tau y_{1}(y_{1}\bar{e}\otimes e)-y_{2}\tau y_{1}(y_{1}\bar{e}\otimes e')+y_{1}y_{2}\xi'\\
 & =\tau y_{1}y_{2}(\bar{e}\otimes e)-y_{2}\tau y_{1}y_{2}(\bar{e}\otimes e')+y_{1}y_{2}\xi'\\
 & =y_{1}y_{2}\left(\tau(\bar{e}\otimes e)-y_{2}\tau(\bar{e}\otimes e')+\xi'\right).
\end{align*}
\item $I_{\beta}'':FE[y]_{G}G_{2}\xrightarrow{\iota\otimes G_{2}}(L_{1})_{G}G_{2}\cong G_{2}$,
\item $I_{\delta}':(G_{1})_{G}FE^{2}[y]\iso FE^{2}[y]$ given as the inverse
of $\bigl(fee\mapsto1_{G_{1}}\otimes fee\bigr)$,
\item $I_{\delta}'':(G_{1})_{G}FE^{2}[y]\xrightarrow{G_{1}\otimes\iota''}(L_{1})_{G}G_{2}\cong G_{2}$.\\
\end{itemize}
Consider the last two maps. We have that $I_{\delta}''\circ I_{\delta}'^{-1}=\iota''$
as maps $FE^{2}[y]\to G_{2}$. Now consider the first two maps. Observe
that $I_{\beta}''=\iota''\circ I_{\beta}'$. It follows that $\mathrm{Im}(I_{\beta}'-I_{\beta}'')\subset\mathrm{Im}(I_{\delta}'-I_{\delta}'')$,
so every element of $FE^{2}[y]$ is associated in the quotient to
a unique element of $G_{2}$ by applying the map $\iota''$.\\

\item We have $\Gamma_{22}:F^{2}E^{2}[y]\oplus(L_{2})_{G}G_{2}\to U$ given
by $(\iota''',\idop_{U})$.

Here the map $\iota''':F^{2}E^{2}[y]\to U$ is the inclusion of the
fifth summand as written in Prop.~\ref{prop:various-as-A-A-bimod}.\\

\begin{itemize}
\item $I_{\beta}':F^{2}E[y]_{G}G_{2}\to F^{2}E^{2}[y]$ given by $\overline{ff}\bar{e}\otimes(e',e,\xi')\mapsto\overline{ff}\otimes(y_{1}y_{2})^{-1}\xi(y_{1}\bar{e})$,
\item $I_{\beta}'':F^{2}E[y]_{G}G_{2}\xrightarrow{\iota'\otimes G_{2}}(L_{2})_{G}G_{2}\iso U$
(using Lemma \ref{lem:U-composition}),
\item $I_{\delta}':(L_{2})_{G}FE^{2}[y]\to F^{2}E^{2}[y]$ given by $(f',f,\rho')\otimes\bar{f}\overline{ee}\mapsto(\bar{f}\circ\rho)\otimes\overline{ee}$,
\item $I_{\delta}'':(L_{2})_{G}FE^{2}[y]\xrightarrow{L_{2}\otimes\iota''}(L_{2})_{G}G_{2}\iso U$.\\
\end{itemize}
Consider the first two maps. Observe that 
\[
I_{\beta}'\bigl(\overline{ff}\bar{e}\otimes(e,y_{1}e,\xi'=0)\bigr)=\overline{ff}\otimes(\bar{e}\otimes e)\in F^{2}E^{2}[y].
\]
 It follows that $I_{\beta}'$ is surjective. Now we show that $\iota'''\circ I_{\beta}'=I_{\beta}''$
and that $\iota'''\circ I_{\delta}'=I_{\delta}''$ using the bimodule
forms. First apply $\iota'''$ to the image under $I_{\beta}'$ of
an arbitrary simple tensor $\overline{ff}\bar{e}\otimes(e',e,\xi')\in F^{2}E[y]_{G}G_{2}$:
\[
\iota'''\bigl(\overline{ff}\otimes(y_{1}y_{2})^{-1}\xi(y_{1}\bar{e})\bigr)=\bigl(0,0,0,0,\Lambda^{\circ}=\overline{ff}\otimes(y_{1}y_{2})^{-1}\xi(y_{1}\bar{e})\bigr),
\]
 then apply $I_{\beta}''$ to the same arbitrary simple tensor, and
view the result through the isomorphism $(L_{2})_{G}G_{2}\iso U$:
\begin{align*}
I_{\beta}''\bigl(\overline{ff}\bar{e}\otimes(e',e,\xi')\bigr) & =(0,0,\overline{ff}\bar{e})\otimes_{G_{1}^{\op}}(e',e,\xi')\\
 & \mapsto\bigl(0,0,0,0,\overline{ff}\otimes(y_{1}y_{2})^{-1}\xi(y_{1}\bar{e})\bigr)\in U.
\end{align*}
 So $\iota'''\circ I_{\beta}'=I_{\beta}''$. Repeat the procedure
with the second pair of maps: 
\begin{align*}
\iota'''\bigl((\bar{f}\circ\rho)\otimes\overline{ee}\bigr) & =\bigl(0,0,0,0,(\bar{f}\circ\rho)\otimes\overline{ee}\bigr),\\
I_{\delta}''\bigl((f',f,\rho')\otimes\bar{f}\overline{ee}\bigr) & =(f',f,\rho')\otimes_{G_{1}^{\op}}(0,0,\bar{f}\overline{ee})\\
 & \mapsto\bigl(0,0,0,0,(\bar{f}\circ\rho)\otimes\overline{ee}\bigr),
\end{align*}
 so $\iota'''\circ I_{\delta}'=I_{\delta}''$. It follows that every
element of $F^{2}E^{2}[y]$ is associated in the quotient to a unique
representative in $U$ by applying $\iota'''$.
\begin{rem}
The map $\iota'''$ describes the inclusion of the morphisms of $\Hom_{K^{b}(B)}(R,R)$
that factor through $X_{1}$. The maps $I_{\beta}'$ and $I_{\delta}'$
are in fact isomorphisms: 
\begin{align*}
\Hom_{K^{b}(B)}(X_{1},X_{2})\otimes_{G_{1}^{\op}}\Hom_{K^{b}(B)}(X_{2},R) & \iso\Hom_{K^{b}(B)}(X_{1},R),\\
\Hom_{K^{b}(B)}(R,X_{2})\otimes_{G_{1}^{\op}}\Hom_{K^{b}(B)}(X_{2},X_{1}) & \iso\Hom_{K^{b}(B)}(R,X_{1}).
\end{align*}
They are produced by reasoning as in Lemma \ref{lem:U-composition},
using that $R$ is a finite direct sum of summands of $X_{2}$.
\end{rem}

\end{itemize}

\subsubsection{$\tilde{E}\tilde{F}$}

We do not have a matrix presentation of the components of the product
$\tilde{E}\tilde{F}$ from the Rickard equivalence. Instead, in this
section, we proceed by studying the quotient directly, by components,
determining the structure of the quotient itself, as well as the quotient
projection $\Gamma$ from the tensor product over $\Delta$ to the
tensor product over $C$.

As before, in each bulleted section we propose a component of $\Gamma$.
Here the arguments following a bulleted line also must justify the
structure of the codomain of the $\Gamma$ component written in that
bulleted line. The domains are known, and in each case the annihilated
submodule $\mathrm{Im}(I_{\beta}'-I_{\beta}'')+\mathrm{Im}(I_{\delta}'-I_{\delta}'')$
is defined already. Our method is to write down a map called $\Gamma_{ij}$
from the appropriate domain, show that it is surjective, and show
that its kernel is $\mathrm{Im}(I_{\beta}'-I_{\beta}'')+\mathrm{Im}(I_{\delta}'-I_{\delta}'')$.
The codomain of $\Gamma$ can be summarized in a matrix: 

\begin{equation}
\begin{split}\tilde{E}\otimes_{\Delta}\tilde{F} & \cong\begin{pmatrix}E[y] & E^{2}[y]\\
G_{1} & G_{2}
\end{pmatrix}\otimes_{\Delta}\begin{pmatrix}F[y] & L_{1}\\
F^{2}[y] & L_{2}
\end{pmatrix}\\
 & \cong\begin{pmatrix}EF[y]\oplus E^{2}[y]_{G}F^{2}[y] & E[y]G_{1}\oplus E^{2}[y]_{G}L_{2}\\
G_{1}F[y]\oplus(G_{2})_{G}F^{2}[y] & G_{1}G_{1}\oplus(G_{2})_{G}L_{2}
\end{pmatrix}\\
 & \xrightarrow{\Gamma}\begin{pmatrix}EF[y] & E[y]G_{1}\\
G_{1}F[y] & G_{1}G_{1}\oplus EF[y]
\end{pmatrix}.
\end{split}
\label{eq:EF-expansion}
\end{equation}
 
\begin{itemize}
\item We have $\Gamma_{11}:EF[y]\oplus E^{2}[y]_{G}F^{2}[y]\to EF[y]$ given
by $(\idop_{EF[y]},\omega)$.\\

Define a map $\omega:E^{2}[y]\otimes_{A[y]}F^{2}[y]\to EF[y]$ by:
\begin{align*}
e_{1}e_{2}\otimes f_{2}f_{1} & \mapsto e_{1}.f_{2}(y_{1}e_{2})\otimes f_{1}=e_{1}\otimes f_{2}(y_{1}e_{2}).f_{1}.
\end{align*}
 Let $\varphi_{1}\in FE[y]$ be given in the second summand of (the
bimodule form) $G_{1}^{\op}\cong A[y]\oplus FE[y]$. Observe that
$\bigl(e_{1}\otimes\varphi_{1}(y_{1}e_{2})\bigr)\otimes f_{2}f_{1}$
and $e_{1}e_{2}\otimes\bigl((f_{2}\circ y_{1}\varphi_{1})\otimes f_{1}\bigr)$
are both sent by $\omega$ to $e_{1}.(f_{2}\ci y_{1}\varphi_{1})(y_{1}e_{2})\otimes f_{1}$.
This means $\omega$ is middle-linear over generators in both summands
of $G_{1}^{\op}$, so it descends to a map, also called $\omega$,
from the tensor product $E^{2}[y]_{G}F^{2}[y]$ taken over $G_{1}^{\op}$.\\

\begin{itemize}
\item $I_{\beta}':EE[y]_{G}F^{2}[y]\to EF[y]$ given (using bimodule forms)
by 
\begin{align*}
e_{1}\otimes e_{2}\otimes f_{2}f_{1} & \mapsto e_{1}\otimes\bigl(f_{1}\circ Ef_{2}\circ(e_{2},y_{1}e_{2},0)\bigr)\\
 & =e_{1}\otimes\bigl(f_{1}\circ Ef_{2}\circ(\_\otimes y_{1}e_{2})\bigr)\\
 & =e_{1}\otimes f_{1}(\_.f_{2}(y_{1}e_{2}))\\
 & =e_{1}\otimes f_{2}(y_{1}e_{2}).f_{1}.
\end{align*}
In this calculation, $e_{2}$ is interpreted as a map of complexes
$X_{1}\to X_{2}$ which induces a map $X_{2}\to R$ (from $E'X_{1}\to E'X_{2}$),
precomposition with which gives the left action of $e_{2}\in E[y]\subset C$
on $F^{2}[y]\subset\tilde{F}$. The induced map corresponds to $(e_{2},y_{1}e_{2},0)$
in $G_{2}$ (strictly speaking, in $\bar{G}_{2}'$) in bimodule form.
Further, $f_{2}f_{1}\in F^{2}[y]$ is interpreted as a map $R\to X_{1}$,
which is identified by $E^{2}[y]\overset{Ef_{2}}{\to}E[y]\overset{f_{1}}{\to}A[y]$
applied to $E^{2}[y]$ in the top row of $R_{0}$. The composite $X_{2}\to X_{1}$
is identified (in the second row) by the morphism $f_{1}\circ Ef_{2}\circ(\_\otimes y_{1}e_{2}):E[y]\to A[y]$,
which is evaluated in the third and fourth lines.\\
\item $I_{\beta}'':EE[y]_{G}F^{2}[y]\xrightarrow{\idop}E^{2}[y]_{G}F^{2}[y]$,
\item $I_{\delta}':E^{2}[y]_{G}FF[y]\to EF[y]$ given (using bimodule forms)
by 
\begin{align*}
e_{1}e_{2}\otimes f_{2}\otimes f_{1} & \mapsto\bigl((0,f_{2},0)\circ e_{1}e_{2}\bigr)\otimes f_{1}\\
 & =y_{1}^{-1}(Ef_{2})(y_{1}y_{2}(e_{1}e_{2}))\otimes f_{1}\\
 & =e_{1}.f_{2}(y_{1}e_{2})\otimes f_{1}.
\end{align*}
 Here $f_{2}:X_{2}\to X_{1}$ induces $(0,f_{2},0):R\to X_{2}$ in
$L_{2}$. Further, $e_{1}e_{2}:X_{1}\to R$, and the composite map
$X_{1}\to X_{2}$ is identified by applying $Ef_{2}$ to the top row
of $R_{0}$ after putting $y_{1}y_{2}e_{1}e_{2}$ in that term, and
removing the final $y_{1}$ to obtain the bimodule form.\\
\item $I_{\delta}'':E^{2}[y]_{G}FF[y]\xrightarrow{\idop}E^{2}[y]_{G}F^{2}[y]$.\\
\end{itemize}
We see that $I_{\beta}'=\omega$ and $I_{\delta}'=\omega$ after identifying
$EE[y]\cong E^{2}[y]$ and $FF[y]\cong F^{2}[y]$. It follows that
the kernel of $\Gamma_{11}$ is the image of $I_{\beta}'-I_{\beta}''$,
which is also the image of $I_{\delta}'-I_{\delta}''$, and thus $\ker(\Gamma_{11})=\mathrm{Im}(I_{\beta}'-I_{\beta}'')+\mathrm{Im}(I_{\delta}'-I_{\delta}'')$
as desired.
\begin{rem}
The map $\omega$ corresponds on the models to the map given by composition:
\[
\Hom_{K^{b}(B)}(X_{2},R)\otimes_{G_{1}^{\op}}\Hom_{K^{b}(B)}(R,X_{2})\to\Hom_{K^{b}(B)}(X_{2},X_{2}).
\]
\end{rem}

\item We have $\Gamma_{21}:G_{1}F[y]\oplus(G_{2})_{G}F^{2}[y]\to G_{1}F[y]$
given by $(\idop_{G_{1}F[y]},\omega')$.\\

Let $\omega':(G_{2})_{G}F^{2}[y]\to G_{1}F[y]$ be defined (using
bimodule forms) by 
\begin{align*}
(e',e,\xi')\otimes f_{2}f_{1} & \mapsto\bigl((0,f_{2},0)\circ(e',e,\xi')\bigr)\otimes f_{1}\\
 & =\Bigl(f_{2}(e),y_{1}^{-1}Ef_{2}\circ\bigl(y_{2}\tau(\_\otimes(e-y_{1}e'))+y_{1}y_{2}\xi'\bigr)\Bigr)\otimes f_{1}\\
 & =\Bigl(f_{2}(e),Ef_{2}\circ\tau(\_\otimes(e-y_{1}e'))+E(f_{2}\circ y_{1})\circ\xi'\Bigr)\otimes f_{1}.
\end{align*}
 Here $f_{2}:X_{2}\to X_{1}$ again induces $(0,f_{2},0):R\to X_{2}$
in $L_{2}$. The composite 
\[
X_{2}\overset{(e',e,\xi')}{\longrightarrow}R\overset{(0,f_{2},0)}{\longrightarrow}X_{2}
\]
is a map $X_{2}\to X_{2}$ identified by the element of $G_{1}$ given
in the next line.\\

\begin{itemize}
\item $I_{\beta}':G_{1}E[y]_{G}F^{2}[y]\to G_{1}F[y]$ given (using bimodule
form) by 
\begin{align*}
ge\otimes f_{2}f_{1} & \mapsto g\otimes\bigl((f_{1}\circ Ef_{2})\circ(e,y_{1}e,0)\bigr)\\
 & =g\otimes f_{2}(y_{1}e).f_{1}.
\end{align*}
 The map $E[y]_{G}F^{2}[y]\to F[y]$ used here is the same as the
one in $I_{\beta}'$ of $\Gamma_{11}$ above.\\
\item $I_{\beta}'':G_{1}E[y]_{G}F^{2}[y]\to(G_{2})_{G}F^{2}[y]$ given (using
bimodule forms) by 
\begin{align*}
(\theta,\varphi_{1})\otimes e\otimes f_{2}f_{1} & \mapsto\bigl((e,y_{1}e,0)\circ(\theta,\varphi_{1})\bigr)\otimes f_{2}f_{1}\\
 & =(\theta e,\theta y_{1}e,\varphi_{1}(-)\otimes e)\otimes f_{2}f_{1}.
\end{align*}
The map $G_{1}E[y]\to G_{2}$ used here is the same as the one in
$\Gamma_{21}$ for $\tilde{E}\tilde{E}$.\\
\item $I_{\delta}':(G_{2})_{G}FF[y]\to G_{1}F[y]$ given by the map $\omega'$
(after identifying $FF[y]$ with $F^{2}[y]$),
\item $I_{\delta}'':(G_{2})_{G}FF[y]\xrightarrow{\idop}(G_{2})_{G}F^{2}[y]$.\\
\end{itemize}
We show that $\omega'\circ I_{\beta}''=I_{\beta}'$: 
\begin{align*}
\omega'\bigl((\theta e,\theta y_{1}e,\varphi_{1}\otimes e)\otimes f_{2}f_{1}\bigr) & =\bigl(f_{2}(\theta y_{1}e),E(f_{2}\circ y_{1})\circ(\varphi_{1}\otimes e)\bigr)\otimes f_{1}\\
 & =\bigl(\theta f_{2}(y_{1}e),\varphi_{1}.f_{2}(y_{1}e)\bigr)\otimes f_{1}\\
 & =(\theta,\varphi_{1}).f_{2}(y_{1}e)\otimes f_{1}\\
 & =I_{\beta}'\bigl((\theta,\varphi_{1})\otimes e\otimes f_{2}f_{1}\bigr).
\end{align*}
In the first line, note that $(\theta y_{1}e)-y_{1}(\theta e)=0$
so the term \textquoteleft $Ef_{2}\circ\tau(\_\otimes(e-y_{1}e'))$\textquoteright{}
in the image under $\omega'$ disappears. Then $f_{2}\ci y_{1}$ applied
to $e$ produces $f_{2}(y_{1}e)\in A[y]$, which acts on the right
on $\varphi_{1}$ for the second line. For the third line, the element
$f_{2}(y_{1}e)\in A[y]$ acts on $G_{1}$ on the right diagonally.
It follows that $I_{\beta}'-I_{\beta}''=(\omega'-\idop)I_{\beta}''$,
and therefore $\mathrm{Im}(I_{\beta}'-I_{\beta}'')\subset\mathrm{Im}(I_{\delta}'-I_{\delta}'')$.
Thus $\ker(\Gamma_{21})=\mathrm{Im}(I_{\beta}'-I_{\beta}'')+\mathrm{Im}(I_{\delta}'-I_{\delta}'')$,
as desired.\\

\item We have $\Gamma_{12}:E[y]G_{1}\oplus E^{2}[y]_{G}L_{2}\to E[y]G_{1}$
given by $(\idop_{E[y]G_{1}},\omega'')$.\\

Let $\omega'':E^{2}[y]_{G}L_{2}\to E[y]G_{1}$ be defined (using bimodule
forms) by 
\begin{align*}
e_{1}e_{2}\otimes(f',f,\rho') & \mapsto e_{1}\otimes\bigl((f',f,\rho')\circ(e_{2},y_{1}e_{2},0)\bigr)\\
 & =e_{1}\otimes\bigl(f(y_{1}e_{2})+f'(e_{2}),Ef'\circ\tau(\_\otimes e_{2})+\rho'(\_\otimes y_{1}e_{2})\bigr).
\end{align*}

\begin{itemize}
\item $I_{\beta}':EE[y]_{G}L_{2}\to E[y]G_{1}$ given by the map $\omega''$
(after identifying $EE[y]$ with $E^{2}[y]$),
\item $I_{\beta}'':EE[y]_{G}L_{2}\xrightarrow{\idop}E^{2}[y]_{G}L_{2}$,
\item $I_{\delta}':E^{2}[y]_{G}F[y]G_{1}\to E[y]G_{1}$ given (borrowing
from $I_{\delta}'$ of $\Gamma_{11}$) by 
\[
e_{1}e_{2}\otimes f_{2}\otimes g\mapsto e_{1}\otimes f_{2}(y_{1}e_{2}).g,
\]
\item $I_{\delta}'':E^{2}[y]_{G}F[y]G_{1}\to E^{2}[y]_{G}L_{2}$ given (using
bimodule forms) by 
\begin{align*}
e_{1}e_{2}\otimes f\otimes(\theta,\varphi_{1}) & \mapsto e_{1}e_{2}\otimes\bigl((\theta,\varphi_{1})\circ(0,f,0)\bigr)\\
 & =e_{1}e_{2}\otimes(0,f.\theta,f\otimes\varphi_{1}).
\end{align*}
Here $f:X_{2}\to X_{1}$ induces $(0,f,0):R\to X_{2}$, and the reader
may check the composition with $(\theta,\varphi_{1}):X_{2}\to X_{2}$.\\
\end{itemize}
We show that $\omega''\circ I_{\delta}''=I_{\delta}'$: 
\begin{align*}
\omega''\bigl(e_{1}e_{2}\otimes(0,f.\theta,\varphi_{1}\circ Ef)\bigr) & =e_{1}\otimes\bigl(f(y_{1}e_{2}).\theta,(\varphi_{1}\circ Ef)(\_\otimes y_{1}e_{2})\bigr)\\
 & =e_{1}\otimes\bigl(f(y_{1}e_{2}).\theta,\varphi_{1}(\_.f(y_{1}e_{2}))\bigr)\\
 & =e_{1}\otimes\bigl(f(y_{1}e_{2}).\theta,f(y_{1}e_{2}).\varphi_{1}\bigr)\\
 & =e_{1}\otimes f(y_{1}e_{2}).(\theta,\varphi_{1})\\
 & =I_{\delta}'\bigl(e_{1}e_{2}\otimes f\otimes(\theta,\varphi_{1})\bigr).
\end{align*}
 Thus $I_{\delta}'-I_{\delta}''=(\omega''-\idop)I_{\delta}''$, and
therefore $\mathrm{Im}(I_{\delta}'-I_{\delta}'')\subset\mathrm{Im}(I_{\beta}'-I_{\beta}'')$.
It follows that $\ker(\Gamma_{12})=\mathrm{Im}(I_{\beta}'-I_{\beta}'')+\mathrm{Im}(I_{\delta}'-I_{\delta}'')$,
as desired.\\

\item We have $\Gamma_{22}:G_{1}G_{1}\oplus(G_{2})_{G}L_{2}\to G_{1}G_{1}\oplus EF[y]$
given by $\begin{pmatrix}\idop_{G_{1}G_{1}} & \omega'''\\
0 & \kappa
\end{pmatrix}$.\\

Below we describe the maps $I_{\beta}'$, $I_{\beta}''$, $I_{\delta}'$,
$I_{\delta}''$, and define a map $\omega''':(G_{2})_{G}L_{2}\to G_{1}G_{1}$,
and we show that $\omega'''\circ I_{\beta}''=I_{\beta}'$ and $\omega'''\circ I_{\delta}''=I_{\delta}'$.
Then we describe a decomposition of $(G_{2})_{G}L_{2}$ into $(A[y],A[y])$-sub-bimodules
$(G_{2})_{G}L_{2}\cong H\oplus EF[y]$ where $H=\mathrm{Im}(I_{\beta}'')+\mathrm{Im}(I_{\delta}'')$.
The projection onto $EF[y]$ is called $\kappa$. (This copy of $EF[y]$
lies in the kernel of $\omega'''$.) From all this it follows that
$\ker(\Gamma_{22})=\mathrm{Im}(I_{\beta}'-I_{\beta}'')+\mathrm{Im}(I_{\delta}'-I_{\delta}'')$
and $\Gamma_{22}$ describes the projection to the quotient.\\

\begin{itemize}
\item $I_{\beta}':G_{1}E[y]_{G}L_{2}\to G_{1}G_{1}$ given (borrowing from
$I_{\beta}'$ of $\Gamma_{12}$) by 
\begin{flalign*}
 & g\otimes e\otimes(f',f,\rho')\mapsto\\
 & \quad g\otimes\bigl(f'(e)+f(y_{1}e),Ef'\circ\tau(\_\otimes e)+\rho'(\_\otimes y_{1}e)\bigr),
\end{flalign*}
\item $I_{\beta}'':G_{1}E[y]_{G}L_{2}\to(G_{2})_{G}L_{2}$ given (borrowing
from $I_{\beta}''$ of $\Gamma_{21}$) by 
\[
(\theta,\varphi_{1})\otimes e\otimes\ell\mapsto(\theta e,\theta y_{1}e,\varphi_{1}(-)\otimes e)\otimes\ell,
\]
\item $I_{\delta}':(G_{2})_{G}F[y]G_{1}\to G_{1}G_{1}$ given (borrowing
from $I_{\delta}'$ of $\Gamma_{21}$) by 
\begin{flalign*}
 & (e',e,\xi')\otimes f\otimes g\mapsto\\
 & \quad\bigl(f(e),Ef\circ\tau(\_\otimes(e-y_{1}e'))+E(f\circ y_{1})\circ\xi'\bigr)\otimes g,
\end{flalign*}
\item $I_{\delta}'':(G_{2})_{G}F[y]G_{1}\to(G_{2})_{G}L_{2}$ given (borrowing
from $I_{\delta}''$ of $\Gamma_{12}$) by 
\[
g\otimes f\otimes(\theta,\varphi_{1})\mapsto g\otimes(0,f.\theta,f\otimes\varphi_{1}).
\]
\\
\end{itemize}
Now we define a morphism of $(A[y],A[y])$-bimodules $\omega''':G_{2}\otimes_{A[y]}L_{2}\to G_{1}G_{1}$,
and in a subsequent lemma we show that $\omega'''$ descends to a
morphism $\omega''':G_{2}\otimes_{G_{1}^{\op}}L_{2}\to G_{1}G_{1}$
by showing that it is also middle-linear over generators of $G_{1}^{\op}$
in $FE[y]$. (Since $G_{1}\cong A[y]\oplus FE[y]$ as a bimodule,
this ensures linearity over all of $G_{1}^{\op}$.) Let $(e',e,\xi')\otimes(f',f,\rho')\in G_{2}\otimes_{A[y]}L_{2}$
be an arbitrary simple tensor. We define: 
\begin{multline*}
\omega''':(e',e,\xi')\otimes(f',f,\rho')\mapsto\Bigl(\veps(e'\otimes f')+\veps(e\otimes f),FE(\varepsilon\circ y_{1}F)(\xi'\otimes f)\\
+FE\veps(\xi'\otimes f')+\sigma(e\otimes f)-\sigma(y_{1}e'\otimes f)\Bigr)\otimes(1,0)\\
+(1,0)\otimes\Bigl(0,\veps FE(e\otimes\rho')+\sigma(e'\otimes f')\Bigr)+\sigma FE(e\otimes\rho')\\
-\sigma FE(y_{1}e'\otimes\rho')+FE\sigma(\xi'\otimes f')+FE(\veps\circ y_{1}F)FE(\xi'\otimes\rho').
\end{multline*}
 The last four terms, beginning with $\sigma FE(e\otimes\rho')$,
are elements of $FEFE[y]$. They should be interpreted in the last
summand of $G_{1}G_{1}$ that appears in the following decomposition
of bimodules: 
\begin{equation}
\begin{split}G_{1}\otimes_{A[y]}G_{1} & \iso A[y]\oplus FE[y]\oplus FE[y]\oplus FEFE[y],\\
(\theta,\varphi_{1})\otimes(\theta',\varphi_{1}') & \mapsto\bigl(\theta\theta',\theta.\varphi_{1}',\varphi_{1}.\theta',\varphi_{1}\otimes\varphi_{1}'\bigr).
\end{split}
\label{eq:G1G1-ordered-decomp}
\end{equation}
 At this point $\omega'''$ has been defined as a map $G_{2}\otimes_{A[y]}L_{2}\to G_{1}G_{1}$.
It is useful to go further and record the data of $\omega'''$ as
a matrix. We can give a decomposition of $G_{2}\otimes_{A[y]}L_{2}$
into a direct sum of $(A[y],A[y])$-bimodules: 
\begin{align*}
 & G_{2}\otimes_{A[y]}L_{2}\iso EF[y]^{\oplus4}\oplus FE^{2}F[y]^{\oplus2}\oplus EF^{2}E[y]^{\oplus2}\oplus FE^{2}F^{2}E[y],\\
 & (e',e,\xi')\otimes(f',f,\rho')\mapsto(e'\otimes f',e'\otimes f,e\otimes f',e\otimes f)\\
 & \qquad\qquad\qquad\qquad\oplus(e'\otimes\rho',e\otimes\rho')\oplus(\xi'\otimes f',\xi'\otimes f)\oplus(\xi'\otimes\rho').
\end{align*}
 Each of the terms in the formula for $\omega'''$ is a morphism of
$(A[y],A[y])$-bimodules.
\begin{defn}
Using the ordered decompositions of $G_{2}\otimes_{A[y]}L_{2}$ and
of $G_{1}G_{1}$ above, the map $\omega''':G_{2}\otimes_{A[y]}L_{2}\to G_{1}G_{1}$
is given by the following matrix: {\scriptsize{}
\[
\begin{pmatrix}\veps & 0 & 0 & \veps & 0 & 0 & 0 & 0 & 0\\
\sigma & 0 & 0 & 0 & 0 & \veps FE & 0 & 0 & 0\\
0 & -\sigma\circ y_{1}F & 0 & \sigma & 0 & 0 & FE\veps & FE(\varepsilon\circ y_{1}F) & 0\\
0 & 0 & 0 & 0 & -(\sigma\circ y_{1}F)FE & \sigma FE & FE\sigma & 0 & FE(\veps\circ y_{1}F)FE
\end{pmatrix}.
\]
}{\scriptsize\par}
\end{defn}

\begin{lem}
\label{lem:omega'''-middle-linear} The map $\omega'''$ is middle-linear
over the action of generators of the summand $FE[y]\subset G_{1}^{\op}$.
\end{lem}

\begin{proof}
We first compute the middle actions $(e',e,\xi').\varphi_{1}$ and
$\varphi_{1}.(f',f,\rho')$ for $\varphi_{1}\in FE[y]\subset G_{1}^{\op}$,
$(e',e,\xi')\in G_{2}$, and $(f',f,\rho')\in L_{2}$, both in bimodule
form. These are: 
\begin{align*}
(e',e,\xi').\varphi_{1} & =\bigl(\varphi_{1}(e),y_{1}\varphi_{1}(e),E\varphi_{1}\circ\tau(\_\otimes(e-y_{1}e'))+E(\varphi_{1}y_{1})\circ\xi'\bigr)\\
\varphi_{1}.(f',f,\rho') & =\bigl(0,f\circ y_{1}\varphi_{1}+f'\circ\varphi_{1},Ef'\circ\tau\circ E\varphi_{1}+\rho'\circ E(y_{1}\varphi_{1})\bigr).
\end{align*}
 Using the formulas above, one easily computes the images under $\omega'''$
of $(e',e,\xi').\varphi_{1}\otimes(f',f,\rho')$ and $(e',e,\xi')\otimes\varphi_{1}.(f',f,\rho')$
and checks that they agree.
\end{proof}
\begin{cor}
It follows from Lemma \ref{lem:omega'''-middle-linear} that $\omega'''$
determines a morphism of $(A[y],A[y])$-bimodules $\omega''':(G_{2})_{G}L_{2}\to G_{1}G_{1}$.
\end{cor}

We show next that $\omega'''\circ I_{\beta}''=I_{\beta}'$ and $\omega'''\circ I_{\delta}''=I_{\delta}'$.
The formula for $\omega'''$ is determined by these conditions and
may be derived from them. Evaluating the left side of the first equation:
\begin{flalign*}
 & \omega'''\circ I_{\beta}''\bigl((\theta,\varphi_{1})\otimes e\otimes(f',f,\rho')\bigr)\\
 & =\omega'''\bigl((\theta e,\theta y_{1}e,\varphi_{1}\otimes e)\otimes(f',f,\rho')\bigr)\\
 & =\bigl(f'(\theta e)+f(\theta y_{1}e),\varphi_{1}.f(y_{1}e)+\varphi_{1}.f'(e)\bigr)\otimes(1,0)\\
 & \quad+(1,0)\otimes\bigl(0,\rho'(\_\otimes\theta y_{1}e)+Ef'\circ\tau(\_\otimes\theta e)\bigr)\\
 & \quad+(0,\varphi_{1})\otimes\bigl(0,Ef'\circ\tau(\_\otimes e)+\rho'(\_\otimes y_{1}e)\bigr)\\
 & =\bigl(\theta.(f'(e)+f(y_{1}e)),\varphi_{1}.(f'(e)+f(y_{1}e))\bigr)\otimes(1,0)\\
 & \quad+(\theta,\varphi_{1})\otimes\bigl(0,Ef'\circ\tau(\_\otimes e)+\rho'(\_\otimes y_{1}e)\bigr)\\
 & =(\theta,\varphi_{1})\otimes\bigl(f'(e)+f(y_{1}e),Ef'\circ\tau(\_\otimes e)+\rho'(\_\otimes y_{1}e)\bigr)\\
 & =I_{\beta}'\bigl((\theta,\varphi_{1})\otimes e\otimes(f',f,\rho')\bigr).
\end{flalign*}
 Now evaluating the left side of the second equation: 
\begin{flalign*}
 & \omega'''\circ I_{\delta}''\bigl((e',e,\xi')\otimes f\otimes(\theta,\varphi_{1})\bigr)\\
 & =\omega'''\bigl((e',e,\xi')\otimes(0,f.\theta,f\otimes\varphi_{1})\bigr)\\
 & =\bigl(f(e).\theta,E(f.\theta\circ y_{1})\circ\xi'+E(f.\theta)\circ\tau(\_\otimes(e-y_{1}e'))\bigr)\otimes(1,0)\\
 & \quad+(1,0)\otimes\bigl(0,f(e).\varphi_{1}\bigr)+\bigl(0,Ef\circ\tau(\_\otimes(e-y_{1}e')\bigr)\otimes(0,\varphi_{1})\\
 & \quad+\bigl(0,E(f\ci y_{1})\circ\xi'\bigr)\otimes(0,\varphi_{1})\\
 & =\bigl(f(e),E(f\ci y_{1})\circ\xi'+Ef\circ\tau(\_\otimes(e-y_{1}e'))\bigr)\otimes(\theta,0)\\
 & \quad+\bigl(f(e),Ef\circ\tau(\_\otimes(e-y_{1}e'))+E(f\ci y_{1})\circ\xi'\bigr)\otimes(0,\varphi_{1})\\
 & =\bigl(f(e),Ef\circ\tau(\_\otimes(e-y_{1}e'))+E(f\ci y_{1})\circ\xi'\bigr)\otimes(\theta,\varphi_{1})\\
 & =I_{\delta}'\bigl((e',e,\xi')\otimes f\otimes(\theta,\varphi_{1})\bigr).
\end{flalign*}

Now the product $(G_{2})_{G}L_{2}$ is the quotient of the product
$(G_{2})_{A[y]}L_{2}$ by the image of $\gamma'-\gamma''$, where:
\\

\begin{itemize}
\item $\gamma':\bigl(G_{2}\otimes_{A[y]}FE[y]\bigr)\otimes_{A[y]}L_{2}\to(G_{2})_{A[y]}L_{2}$
given by 
\begin{align*}
 & (e',e,\xi')\otimes\varphi_{1}\otimes\ell\mapsto\\
 & \bigl(\varphi_{1}(e),y_{1}\varphi_{1}(e),E\varphi_{1}\circ\tau(\_\otimes(e-y_{1}e'))+E(\varphi_{1}\ci y_{1})\circ\xi'\bigr)\otimes\ell,
\end{align*}
\item $\gamma'':G_{2}\otimes_{A[y]}\bigl(FE[y]\otimes_{A[y]}L_{2}\bigr)\to(G_{2})_{A[y]}L_{2}$
given by 
\begin{align*}
 & g\otimes\varphi_{1}\otimes(f',f,\rho')\mapsto\\
 & g\otimes\bigl(0,f'\circ\varphi_{1}+f\circ y_{1}\varphi_{1},Ef'\circ\tau\circ E\varphi_{1}+\rho'\circ E(y_{1}\varphi_{1})\bigr).
\end{align*}
\end{itemize}
There is a copy of $EF[y]$ in $(G_{2})_{A[y]}L_{2}$ generated by
terms of the form $(0,e,0)\otimes(f',0,0)$. Let $\bar{H}$ be its
direct complement. The images of $\gamma'$ and $\gamma''$ lie in
$\bar{H}$, so $(G_{2})_{G}L_{2}\cong H\oplus EF[y]$, where $H$
is the quotient of $\bar{H}$ by the image of $\gamma'-\gamma''$.

The image of $I_{\beta}''$ includes every term of the form $(e,y_{1}e,\varphi_{1}\otimes e)\otimes\ell$,
and the image of $I_{\delta}''$ includes every term of the form $g\otimes(0,f,f\otimes\varphi_{1})$.
By adding appropriate linear combinations of terms of the first form,
one obtains any element $(e,y_{1}e,\xi')\otimes\ell$, and similarly
from terms of the second form one obtains any $g\otimes(0,f,\rho')$.
It follows that $\mathrm{Im}(I_{\beta}''+I_{\delta}'')=H$.
\end{itemize}

\subsection{Maps $\tilde{\rho}_{\lambda}$: formulas}

In this section we derive formulas by matrix components for the maps
$\tilde{\sigma}=\tilde{F}\tilde{E}\tilde{\veps}\ci\tilde{F}\tilde{\tau}\tilde{F}\ci\tilde{\eta}\tilde{E}\tilde{F}$,
$\tilde{\veps}\circ\tilde{x}^{i}\tilde{F}$, and $\tilde{F}\tilde{x}^{i}\circ\tilde{\eta}$
that are used to define the maps $\tilde{\rho}_{\lambda}$. We will
be using the matrix components for $\tilde{E}$, $\tilde{F}$, $\tilde{E}\tilde{E}$,
$\tilde{F}\tilde{E}$, and $\tilde{E}\tilde{F}$ that were found and
studied in previous sections. (See Eqs.~(\ref{eq:Matrix-tilde-E}),
(\ref{eq:Matrix-tilde-F}), (\ref{eq:Matrix-tilde-E^2}), (\ref{eq:Matrix-tilde-FE}),
and (\ref{eq:EF-expansion}), respectively.) The unit and counit $\tilde{\eta}$
and $\tilde{\veps}$ are given by the duality pairing and thus are
easily interpreted in terms of maps between complexes where that is
convenient. The morphisms $\tilde{x}$ and $\tilde{\tau}$ were given
on components in Eqs.~(\ref{eq:Def-tilde-x}) and (\ref{eq:Def-tilde-tau}).

\subsubsection{Map $\tilde{\sigma}:\tilde{E}\tilde{F}\to\tilde{F}\tilde{E}$}

We begin by computing the map $\tilde{\sigma}:\tilde{E}\tilde{F}\to\tilde{F}\tilde{E}$.
Recall that $\tilde{\sigma}$ is defined by $\tilde{\sigma}=\tilde{F}\tilde{E}\tilde{\veps}\ci\tilde{F}\tilde{\tau}\tilde{F}\ci\tilde{\eta}\tilde{E}\tilde{F}$,
and $\tilde{\eta}$, $\tilde{\veps}$, and $\tilde{\tau}$ are determined
already. We will need formulas for each component of $\tilde{\sigma}$
in its matrix presentation.

We use the following technique to derive the formulas. We start with
an appropriate matrix coefficient of the element $[\tilde{\eta}(1)]\in[\tilde{F}\tilde{E}]$,
together with an arbitrary generator of a component of the matrix
$[\tilde{E}\tilde{F}]$. Then we write the latter as a sum of simple
tensor products of elements of $[\tilde{E}]$ with elements of $[\tilde{F}]$.
As a point of notation, this will be said to lie in $[\tilde{E}]\cdot[\tilde{F}]$
(and similarly for other matrix products). Then we write $[\tilde{\eta}(1)]$
in $[\tilde{F}]\cdot[\tilde{E}]$, and taking another tensor product
we have an element we can write in $[\tilde{F}]\cdot[\tilde{E}]\cdot[\tilde{E}]\cdot[\tilde{F}]$.
Upon this we apply $[\tilde{F}]\cdot[\tilde{\tau}]\cdot[\tilde{F}]$
using (\ref{eq:Def-tilde-tau}). We view the result in $[\tilde{F}]\cdot[\tilde{E}]\cdot[\tilde{E}\tilde{F}]$,
apply $[\tilde{F}]\cdot[\tilde{E}]\cdot[\tilde{\veps}]$ to obtain
an element of $[\tilde{F}]\cdot[\tilde{E}]\cdot[C]$, view this in
$[\tilde{F}\tilde{E}]\cdot[C]$, and allow the coefficient in $[C]$
to act on the right on the coefficient in $[\tilde{F}\tilde{E}]$.
The result is the image under $[\tilde{\sigma}]$ of the arbitrary
generator in $[\tilde{E}\tilde{F}]$ with which we began.

The following bulleted lines state the results of this procedure,
and the procedure itself is carried out in detail in the paragraphs
below those lines.
\begin{itemize}
\item We have $[\tilde{\sigma}]_{11}:[\tilde{E}\tilde{F}]_{11}\to[\tilde{F}\tilde{E}]_{11}$
given by $\left(\begin{smallmatrix}\veps\\
\sigma
\end{smallmatrix}\right)$ using the decompositions: 
\begin{itemize}
\item $[\tilde{E}\tilde{F}]_{11}\cong EF[y]$,
\item $[\tilde{F}\tilde{E}]_{11}\cong(G_{1})_{G}G_{1}\cong G_{1}\cong A[y]\oplus FE[y]$.
\end{itemize}
We take $[\tilde{\eta}(1)]_{11}=(1,0)\otimes(1,0)\in(G_{1})_{G}G_{1}\cong[\tilde{F}\tilde{E}]_{11}$
(using bimodule form), and an arbitrary generator $e\otimes f\in EF[y]\cong[\tilde{E}\tilde{F}]_{11}$.
The product of these in $[\tilde{F}\tilde{E}]\cdot[\tilde{E}\tilde{F}]$
can be represented in $[\tilde{F}]\cdot[\tilde{E}]\cdot[\tilde{E}]\cdot[\tilde{F}]$
by: 
\begin{multline*}
\left(\begin{smallmatrix}0 & (1,0)\\
0 & 0
\end{smallmatrix}\right)\cdot\left(\begin{smallmatrix}0 & 0\\
(1,0) & 0
\end{smallmatrix}\right)\cdot\left(\begin{smallmatrix}e & 0\\
0 & 0
\end{smallmatrix}\right)\cdot\left(\begin{smallmatrix}f & 0\\
0 & 0
\end{smallmatrix}\right)\in\left(\begin{smallmatrix}F[y] & L_{1}\\
F^{2}[y] & L_{2}
\end{smallmatrix}\right)\cdot\left(\begin{smallmatrix}E[y] & E^{2}[y]\\
G_{1} & G_{2}
\end{smallmatrix}\right)\cdot\left(\begin{smallmatrix}E[y] & E^{2}[y]\\
G_{1} & G_{2}
\end{smallmatrix}\right)\cdot\left(\begin{smallmatrix}F[y] & L_{1}\\
F^{2}[y] & L_{2}
\end{smallmatrix}\right).
\end{multline*}
 The middle factors give $(1,0)\otimes e\in G_{1}\otimes_{A[y]}E[y]$.
Passing through $\Gamma_{21}$ of Eq.~(\ref{eq:EE-expansion}), this
represents $(e,y_{1}e,0)\in G_{2}\cong[\tilde{E}^{2}]_{21}$. To apply
$[\tilde{\tau}]_{21}$ from Eq.~(\ref{eq:Def-tilde-tau}) we translate
that formula from the terms of $\bar{G}_{2}$ to those of $\bar{G}_{2}'$
in bimodule form. Using $e_{1}-e_{2}=y_{1}e'$, we have: 
\begin{align*}
\xi & =\_\otimes e_{1}+y_{2}\tau(\_\otimes e_{2})+y_{1}y_{2}\xi'\\
 & =\_\otimes(e_{1}-e_{2})+\tau y_{1}(\_\otimes e_{2})+y_{1}y_{2}\xi';\\
\tau\ci\xi & =\tau(\_\otimes(e_{1}-e_{2}))+\tau^{2}y_{1}(\_\otimes e_{2})+\tau y_{1}y_{2}\xi'\\
 & =\tau y_{1}(\_\otimes e')+y_{1}y_{2}\tau\ci\xi'\\
 & =\_\otimes e'+y_{2}\tau(\_\otimes e')+y_{1}y_{2}(\tau\ci\xi').
\end{align*}
 So instead of $(e_{1},e_{2},\xi)\overset{\tilde{\tau}}{\mapsto}(e',e',\tau\ci\xi)$
in the terms of $\bar{G}_{2}$, the formula is $(e',e,\xi')\overset{\tilde{\tau}}{\mapsto}(0,e',\tau\ci\xi')$
in the terms of $\bar{G}_{2}'$. Application of $[\tilde{\tau}]_{21}$
to $(e,y_{1}e,0)$ therefore yields $(0,e,0)$, which may be represented
in $[\tilde{E}]\cdot[\tilde{E}]$ by: 
\[
\left(\begin{smallmatrix}0 & 0\\
0 & (0,e,0)
\end{smallmatrix}\right)\cdot\left(\begin{smallmatrix}0 & 0\\
(1,0) & 0
\end{smallmatrix}\right)\in[\tilde{E}]\cdot[\tilde{E}].
\]
 Then: 
\[
\left(\begin{smallmatrix}0 & 0\\
(1,0) & 0
\end{smallmatrix}\right)\cdot\left(\begin{smallmatrix}f & 0\\
0 & 0
\end{smallmatrix}\right)\overset{\tilde{\veps}}{\longmapsto}\left(\begin{smallmatrix}0 & 0\\
f & 0
\end{smallmatrix}\right)\in[C]
\]
 and 
\[
\left(\begin{smallmatrix}0 & (1,0)\\
0 & 0
\end{smallmatrix}\right)\cdot\left(\begin{smallmatrix}0 & 0\\
0 & (0,e,0)
\end{smallmatrix}\right)\overset{\Gamma_{12}}{\longmapsto}\left(\begin{smallmatrix}0 & (0,e,0)\\
0 & 0
\end{smallmatrix}\right)\in\left(\begin{smallmatrix}G_{1} & G_{2}\\
L_{2} & U
\end{smallmatrix}\right)=[\tilde{F}\tilde{E}]
\]
using $\Gamma_{12}$ in Eq.~(\ref{eq:FE-expansion}). Finally letting
$f\in C$ act on the right, we have: 
\[
\left(\begin{smallmatrix}0 & (0,e,0)\\
0 & 0
\end{smallmatrix}\right)\cdot\left(\begin{smallmatrix}0 & 0\\
f & 0
\end{smallmatrix}\right)=\left(\begin{smallmatrix}\bigl(f(e),Ef\circ\tau(\_\otimes e)\bigr) & 0\\
0 & 0
\end{smallmatrix}\right)\in[\tilde{F}\tilde{E}].
\]
 The nonzero coefficient may be interpreted as $\bigl(\veps(e\otimes f),\sigma(e\otimes f)\bigr)$.
\item We have $[\tilde{\sigma}]_{21}:[\tilde{E}\tilde{F}]_{21}\to[\tilde{F}\tilde{E}]_{21}$
given by $\left(\begin{smallmatrix}1 & 0\\
0 & F\veps\\
0 & F\sigma
\end{smallmatrix}\right)$ using the decompositions: 
\begin{itemize}
\item $[\tilde{E}\tilde{F}]_{21}\cong G_{1}F[y]\cong F[y]\oplus FEF[y]$,
\item $[\tilde{F}\tilde{E}]_{21}\cong L_{2}\cong F[y]\oplus F[y]\oplus F^{2}E[y]$.
\end{itemize}
Considering the isomorphism $FE\iso\mathrm{Hom}_{A}(_{A}E,E)$, we
can choose an expression for $\eta(1)\in FE\subset FE[y]$ corresponding
to $\idop_{E}\in\Hom_{A}(_{A}E,E)$ as a sum of simple tensors: 
\[
\eta(1)=\sum_{a\in Q}f_{a}\otimes e_{a}\in FE\subset FE[y],
\]
 where $Q$ is some finite index set. Using $f_{a}$, $e_{a}$ for
$a\in Q$, we find an expression for $[\tilde{\eta}(1)]_{22}$ in
$(L_{2})_{G}G_{2}$: 
\begin{lem}
\label{lem:sigma-21-eta} The element 
\[
\sum_{a\in Q}(f_{a},0,0)\otimes(e_{a},0,0)+\sum_{b\in Q}(0,f_{b},0)\otimes(0,e_{b},0)\in(L_{2})_{G}G_{2}
\]
 (written using bimodule forms) is sent to $\idop_{R}\in U$ under
the composition morphism $(L_{2})_{G}G_{2}\iso U$of Lemma \ref{lem:U-composition}.
We write $[\tilde{\eta}(1)]_{22}$ for this element.
\end{lem}

\begin{proof}
We first take composition of the first sum, and then of the second.
\begin{claim}
Under the map $(L_{2})_{G}G_{2}\iso U$, we have: 
\[
\sum_{a\in Q}(f_{a},0,0)\otimes(e_{a},0,0)\mapsto(0,0,0,\idop_{E[y]},0).
\]
\end{claim}

\begin{proof}[Proof of Claim]
 The matrix $[\Phi]$ giving the degree $1$ bottom row part of the
image, which is a morphism in $\Hom_{K^{b}(B)}(R,R)$ written in $U$,
is $\sum_{a\in Q}\left(\begin{smallmatrix}0 & 0\\
0 & f_{a}(\_)\otimes e_{a}
\end{smallmatrix}\right)=\left(\begin{smallmatrix}0 & 0\\
0 & 1
\end{smallmatrix}\right)$. To compute the fifth coefficient $\Lambda^{\ci}$ of the image,
we find the degree $0$ part $\Lambda$ of the map on the top row,
given by taking the composition $E^{2}[y]\to A[y]\to E^{2}[y]$: 
\begin{align*}
 & \sum_{a\in Q}-y_{2}\tau(\_\otimes y_{1}e_{a})\circ(Ef_{a}\circ\tau)\\
 & =\sum_{a\in Q,d\in P}-y_{2}\tau y_{1}\bigl(\tau(\_\:\_)_{(1d)}\otimes f_{a}(\tau(\_\:\_)_{(2d)}).e_{a}\bigr)\\
 & =-y_{2}\tau y_{1}\tau=-y_{2}\tau
\end{align*}
 (in the second line we introduce notation for a decomposition $\tau(ee)=\sum_{d\in P}\tau(ee)_{(1d)}\otimes\tau(ee)_{(2d)}$
for some choices of $\tau(ee)_{(id)}$, $i=1,2$ and finite index
set $P$, and in the third line we use that $\sum_{a\in Q}f_{a}(e^{*}).e_{a}=e^{*}$
for any $e^{*}\in E[y]$). Then $\Lambda^{\ci}=0$ is determined by
Eq.~(\ref{eq:Lambda^circ}) with this $\Lambda$ and $\Phi$.
\end{proof}
\begin{claim}
Under the map $(L_{2})_{G}G_{2}\iso U$, we have: 
\[
\sum_{b\in Q}(0,f_{b},0)\otimes(0,e_{b},0)\mapsto(\idop_{E[y]},0,0,0,0).
\]
\end{claim}

\begin{proof}[Proof of Claim]
 Computing as above, the matrix $[\Phi]$ is given by $\left(\begin{smallmatrix}1 & 0\\
0 & 0
\end{smallmatrix}\right)$, and we have: 
\begin{align*}
\sum_{b\in Q}\bigl(\_\otimes e_{b}+y_{2}\tau(\_\otimes e_{b})\bigr)\circ Ef_{b} & =\sum_{b\in Q}\tau y_{1}(\_\otimes e_{b})\circ Ef_{b}\\
 & =\tau y_{1}(\sum_{b\in Q}f_{b}(\_).e_{b})=\tau y_{1}.
\end{align*}
 Again, $\Lambda^{\ci}=0$ is determined by Eq.~(\ref{eq:Lambda^circ})
with this $\Lambda$ and $\Phi$.
\end{proof}
So $[\tilde{\eta}(1)]_{22}$ is sent to $(1,0,0,1,0)\in U$, which
indeed corresponds to $\idop_{R}$.
\end{proof}
Then we take an arbitrary generator $(\theta,\varphi_{1})\otimes f\in G_{1}F[y]\cong[\tilde{E}\tilde{F}]_{21}$.
Expressing the product $\tilde{\eta}(1)\otimes(\theta,\varphi_{1})\otimes f$
in $[\tilde{F}]\cdot[\tilde{E}]\cdot[\tilde{E}]\cdot[\tilde{F}]$,
we have: 
\begin{align*}
 & \phantom{+}\sum_{a\in Q}\left(\begin{smallmatrix}0 & 0\\
0 & (f_{a},0,0)
\end{smallmatrix}\right)\cdot\left(\begin{smallmatrix}0 & 0\\
0 & (e_{a},0,0)
\end{smallmatrix}\right)\cdot\left(\begin{smallmatrix}0 & 0\\
(\theta,\varphi_{1}) & 0
\end{smallmatrix}\right)\cdot\left(\begin{smallmatrix}f & 0\\
0 & 0
\end{smallmatrix}\right)\\
 & +\sum_{b\in Q}\left(\begin{smallmatrix}0 & 0\\
0 & (0,f_{b},0)
\end{smallmatrix}\right)\cdot\left(\begin{smallmatrix}0 & 0\\
0 & (0,e_{b},0)
\end{smallmatrix}\right)\cdot\left(\begin{smallmatrix}0 & 0\\
(\theta,\varphi_{1}) & 0
\end{smallmatrix}\right)\cdot\left(\begin{smallmatrix}f & 0\\
0 & 0
\end{smallmatrix}\right)\\
 & \in\left(\begin{smallmatrix}F[y] & L_{1}\\
F^{2}[y] & L_{2}
\end{smallmatrix}\right)\cdot\left(\begin{smallmatrix}E[y] & E^{2}[y]\\
G_{1} & G_{2}
\end{smallmatrix}\right)\cdot\left(\begin{smallmatrix}E[y] & E^{2}[y]\\
G_{1} & G_{2}
\end{smallmatrix}\right)\cdot\left(\begin{smallmatrix}F[y] & L_{1}\\
F^{2}[y] & L_{2}
\end{smallmatrix}\right).
\end{align*}

Now we interpret $\left(\begin{smallmatrix}0 & 0\\
0 & (e_{a},0,0)
\end{smallmatrix}\right)\cdot\left(\begin{smallmatrix}0 & 0\\
(\theta,\varphi_{1}) & 0
\end{smallmatrix}\right)$ and $\left(\begin{smallmatrix}0 & 0\\
0 & (0,e_{b},0)
\end{smallmatrix}\right)\cdot\left(\begin{smallmatrix}0 & 0\\
(\theta,\varphi_{1}) & 0
\end{smallmatrix}\right)$ in $[\tilde{E}^{2}]$ using $\Gamma_{21}$ from Eq.~(\ref{eq:EE-expansion});
this requires the right action of $G_{1}^{\op}$ on $R$ from Lemma
\ref{lem:G-action-R}. Then we apply $[\tilde{\tau}]$. 
\begin{align*}
\Gamma_{21}:(e_{a},0,0)\otimes(\theta,\varphi_{1}) & \mapsto(e_{a},0,0).(\theta,\varphi_{1})\\
 & =\bigl(e_{a}.\theta,0,-E\varphi_{1}\circ\tau(\_\otimes y_{1}e_{a})\bigr)\in G_{2}=[\tilde{E}^{2}]_{21}\\
 & \overset{\tilde{\tau}}{\mapsto}\bigl(0,e_{a}.\theta,-\tau\circ E\varphi_{1}\circ\tau(\_\otimes y_{1}e_{a})\bigr)\in[\tilde{E}^{2}]_{21},\\
\Gamma_{21}:(0,e_{b},0)\otimes(\theta,\varphi_{1}) & \mapsto(0,e_{b},0).(\theta,\varphi_{1})\\
 & =\bigl(\varphi_{1}(e_{b}),\varphi(e_{b}),E\varphi_{1}\circ\tau(\_\otimes e_{b})\bigr)\in[\tilde{E}^{2}]_{21}\\
 & \overset{\tilde{\tau}}{\mapsto}\bigl(0,\varphi_{1}(e_{b}),\tau\circ E\varphi_{1}\circ\tau(\_\otimes e_{b})\bigr)\in[\tilde{E}^{2}]_{21}.
\end{align*}
 We can represent these in $[\tilde{E}]\cdot[\tilde{E}]$ using the
isomorphism $G_{2}\iso(G_{2})_{G}G_{1}$, $g\mapsto g\otimes(1,0)$.
So, after applying $[\tilde{F}]\cdot[\tilde{\tau}]\cdot[\tilde{F}]$
to the middle terms, we have: 
\begin{align*}
 & \phantom{+}\sum_{a\in Q}\left(\begin{smallmatrix}0 & 0\\
0 & (f_{a},0,0)
\end{smallmatrix}\right)\cdot\left(\begin{smallmatrix}0 & 0\\
0 & \bigl(0,e_{a}.\theta,-\tau\circ E\varphi_{1}\circ\tau(\_\otimes y_{1}e_{a})\bigr)
\end{smallmatrix}\right)\cdot\left(\begin{smallmatrix}0 & 0\\
(1,0) & 0
\end{smallmatrix}\right)\cdot\left(\begin{smallmatrix}f & 0\\
0 & 0
\end{smallmatrix}\right)\\
 & +\sum_{b\in Q}\left(\begin{smallmatrix}0 & 0\\
0 & (0,f_{b},0)
\end{smallmatrix}\right)\cdot\left(\begin{smallmatrix}0 & 0\\
0 & \bigl(0,\varphi_{1}(e_{b}),\tau\circ E\varphi_{1}\circ\tau(\_\otimes e_{b})\bigr)
\end{smallmatrix}\right)\cdot\left(\begin{smallmatrix}0 & 0\\
(1,0) & 0
\end{smallmatrix}\right)\cdot\left(\begin{smallmatrix}f & 0\\
0 & 0
\end{smallmatrix}\right).
\end{align*}
 Then $\tilde{\veps}:\left(\begin{smallmatrix}0 & 0\\
(1,0) & 0
\end{smallmatrix}\right)\cdot\left(\begin{smallmatrix}f & 0\\
0 & 0
\end{smallmatrix}\right)\mapsto f\in F[y]\cong[C]_{21}$, so by applying $[\tilde{F}]\cdot[\tilde{E}]\cdot[\tilde{\veps}]$
and viewing the first two factors in $(L_{2})_{G}G_{2}\subset[\tilde{F}\tilde{E}]_{22}$
we obtain: 
\[
\begin{pmatrix}0 & 0\\
0 & \begin{smallmatrix}\phantom{+}\sum_{a\in Q}(f_{a},0,0)\otimes\bigl(0,e_{a}.\theta,-\tau\circ E\varphi_{1}\circ\tau(\_\otimes y_{1}e_{a})\bigr)\\
+\sum_{b\in Q}(0,f_{b},0)\otimes\bigl(0,\varphi_{1}(e_{b}),\tau\circ E\varphi_{1}\circ\tau(\_\otimes e_{b})\bigr)
\end{smallmatrix}
\end{pmatrix}\cdot\begin{pmatrix}0 & 0\\
f & 0
\end{pmatrix}\in[\tilde{F}\tilde{E}]\cdot[C].
\]
 Now we express this element in $L_{2}=[\tilde{F}\tilde{E}]_{21}$
by applying the composition map $(L_{2})_{G}G_{2}\iso U$ and then
evaluating the action of $f\in[C]_{21}$ on the right. The latter
may be computed by embedding $f$ in $L_{2}$ as $(0,f,0)$ and post-composing
with this element.

Passing first through the composition map $(L_{2})_{G}G_{2}\iso U$,
we have: 
\[
\sum_{a\in Q}\bigl(0,e_{a}.\theta,-\tau\circ E\varphi_{1}\circ\tau(\_\otimes y_{1}e_{a})\bigr)\circ(f_{a},0,0)\mapsto(0,0,\theta,0,-\tau\circ E\varphi_{1}\circ\tau).
\]
 In the first components of this calculation, we have used: 
\begin{align*}
\sum_{a\in Q}f_{a}(-).(e_{a}.\theta):E[y] & \to E[y]\\
e & \mapsto\sum_{a\in Q}f_{a}(e).(e_{a}.\theta)=e.\theta,
\end{align*}

and for the last component we have used: 
\begin{align*}
 & \sum_{a\in Q}\bigl(\_\otimes e_{a}.\theta+y_{2}\tau(\_\otimes e_{a}.\theta)-y_{1}y_{2}\tau\circ E\varphi_{1}\circ\tau(\_\otimes y_{1}e_{a})\bigr)\ci(Ef_{a}\ci\tau)\\
 & =\tau.\theta-y_{1}y_{2}\tau\circ E\varphi_{1}\circ\tau\\
 & =\tau y_{1}(E\theta\circ\tau)+y_{1}y_{2}(-\tau\circ E\varphi_{1}\circ\tau).
\end{align*}
 The fact that $\Lambda^{\ci}=-\tau\circ E\varphi_{1}\circ\tau$ can
be deduced by comparing with Eq.~(\ref{eq:Lambda^circ}) where $[\Phi]=\left(\begin{smallmatrix}0 & E\theta\\
0 & 0
\end{smallmatrix}\right)$. Similarly, we have: 
\[
\sum_{b\in Q}\bigl(0,\varphi_{1}(e_{b}),\tau\circ E\varphi_{1}\circ\tau(\_\otimes e_{b})\bigr)\circ(0,f_{b},0)\mapsto(\varphi_{1},0,0,0,\tau\circ E\varphi_{1}\circ\tau),
\]
 where again we have used: 
\begin{align*}
 & \sum_{b\in Q}\bigl(\_\otimes\varphi_{1}(e_{b})+y_{2}\tau(\_\otimes\varphi_{1}(e_{b}))+y_{1}y_{2}\tau\circ E\varphi_{1}\circ\tau(\_\otimes e_{b})\bigr)\ci Ef_{b}\\
 & =\tau y_{1}\circ E\varphi_{1}+y_{1}y_{2}\tau\circ E\varphi_{1}\circ\tau\\
 & =\tau y_{1}(E\varphi_{1})+y_{1}y_{2}(\tau\circ E\varphi_{1}\circ\tau),
\end{align*}
 so $\Lambda^{\ci}=\tau\circ E\varphi_{1}\circ\tau$. For the sum
of the images, we have $(\varphi_{1},0,\theta,0,0)\in U$. Next we
compute the right action of $f\in[C]_{21}$ on this element: 
\[
(0,f,0)\circ(\varphi_{1},0,\theta,0,0)=(\theta.f,f\circ\varphi_{1},Ef\circ\tau\circ E\varphi_{1}),
\]
 where we have used: 
\begin{align*}
 & Ef\circ\bigl(\tau y_{1}(E\varphi_{1}+E\theta\ci\tau)\bigr)\\
 & =Ef\circ\bigl(\tau y_{1}\ci E\varphi_{1}+\tau y_{1}\tau\ci E\theta\bigr)\\
 & =Ef\circ\bigl(\tau y_{1}\ci E\varphi_{1}+\tau\ci E\theta\bigr)\\
 & =Ef\circ\bigl(y_{2}\tau\ci E\varphi_{1}+E\varphi_{1}+E\theta\ci\tau\bigr)\\
 & =E(\theta.f)\circ\tau+E(f\circ\varphi_{1})+y_{1}\circ(Ef\circ\tau\circ E\varphi_{1}).
\end{align*}
Our final expression for the image of $\left(\begin{smallmatrix}0 & 0\\
0 & (\theta,\varphi_{1})\otimes f
\end{smallmatrix}\right)$ under $[\tilde{\sigma}]_{21}$ is therefore: 
\[
\begin{pmatrix}0 & 0\\
\bigl(\theta.f,f\circ\varphi_{1},Ef\circ\tau\circ E\varphi_{1}\bigr) & 0
\end{pmatrix}\in\begin{pmatrix}G_{1} & G_{2}\\
L_{2} & U
\end{pmatrix}=[\tilde{F}\tilde{E}].
\]
 The bulleted statement follows from the fact that $f\circ\varphi_{1}=F\veps(\varphi_{1}\otimes f)$
and $Ef\circ\tau\circ E\varphi_{1}=F\sigma(\varphi_{1}\otimes f)$.
\item We have $[\tilde{\sigma}]_{12}:[\tilde{E}\tilde{F}]_{12}\to[\tilde{F}\tilde{E}]_{12}$
given by $\left(\begin{smallmatrix}0 & \veps E\\
1 & y_{1}\circ\veps E\\
0 & \sigma E
\end{smallmatrix}\right)$ using the decompositions: 
\begin{itemize}
\item $[\tilde{E}\tilde{F}]_{12}\cong E[y]G_{1}\cong E[y]\oplus EFE[y]$,
\item $[\tilde{F}\tilde{E}]_{12}\cong G_{2}\cong E[y]\oplus E[y]\oplus FE^{2}[y]$.
\end{itemize}
We take $[\tilde{\eta}(1)]_{11}=(1,0)\otimes(1,0)\in G_{1}G_{1}\cong[\tilde{F}\tilde{E}]_{11}$,
and an arbitrary generator $e\otimes(\theta,\varphi_{1})\in E[y]G_{1}\cong[\tilde{E}\tilde{F}]_{12}$.
The product of these in $[\tilde{F}\tilde{E}]\cdot[\tilde{E}\tilde{F}]$
can be expressed in $[\tilde{F}]\cdot[\tilde{E}]\cdot[\tilde{E}]\cdot[\tilde{F}]$
by: 
\[
\left(\begin{smallmatrix}0 & (1,0)\\
0 & 0
\end{smallmatrix}\right)\cdot\left(\begin{smallmatrix}0 & 0\\
(1,0) & 0
\end{smallmatrix}\right)\cdot\left(\begin{smallmatrix}e & 0\\
0 & 0
\end{smallmatrix}\right)\cdot\left(\begin{smallmatrix}0 & (\theta,\varphi_{1})\\
0 & 0
\end{smallmatrix}\right),
\]
 and application of $[\tilde{F}]\cdot[\tilde{\tau}]\cdot[\tilde{F}]$
gives: 
\[
\left(\begin{smallmatrix}0 & (1,0)\\
0 & 0
\end{smallmatrix}\right)\cdot\left(\begin{smallmatrix}0 & 0\\
0 & (0,e,0)
\end{smallmatrix}\right)\cdot\left(\begin{smallmatrix}0 & 0\\
(1,0) & 0
\end{smallmatrix}\right)\cdot\left(\begin{smallmatrix}0 & (\theta,\varphi_{1})\\
0 & 0
\end{smallmatrix}\right).
\]
 This is sent by $[\tilde{F}]\cdot[\tilde{E}]\cdot[\tilde{\veps}]$
to 
\[
\left(\begin{smallmatrix}0 & (0,e,0)\\
0 & 0
\end{smallmatrix}\right)\cdot\left(\begin{smallmatrix}0 & 0\\
0 & (\theta,\varphi_{1})
\end{smallmatrix}\right)\in[\tilde{F}\tilde{E}]\cdot[C]
\]
 which, after computing the action using Lemma \ref{lem:G-action-R},
gives 
\[
\bigl(\varphi_{1}(e),\varphi(e),E\varphi_{1}\circ\tau(\_\otimes e)\bigr)\in G_{2}\cong[\tilde{F}\tilde{E}]_{12}.
\]
For the last term, set $\xi=\_\otimes e+y_{2}\tau(\_\otimes e)$.
Then $E\varphi\ci\xi$ is: 
\begin{align*}
E\varphi\ci\xi & =\_\otimes\varphi(e)+E\theta\ci y_{2}\tau(\_\otimes e)+E(y_{1}\varphi_{1})\ci y_{2}\tau(\_\otimes e)\\
 & =\_\otimes\varphi(e)+y_{2}\tau(\_\otimes e.\theta)+y_{1}y_{2}E\varphi_{1}\ci\tau(\_\otimes e).
\end{align*}
 Subtracting $\_\otimes\varphi(e)+y_{2}\tau\bigl(\_\otimes(\varphi(e)-y_{1}\varphi_{1}(e))\bigr)$
to isolate $y_{1}y_{2}\xi'$, we obtain $y_{1}y_{2}E\varphi_{1}\circ\tau(\_\otimes e)$
and the last component follows. For the final result observe that
$E\varphi_{1}\circ\tau(\_\otimes e)=\sigma E(e\otimes\varphi_{1}).$
\item We have $[\tilde{\sigma}]_{22}:[\tilde{E}\tilde{F}]_{22}\to[\tilde{F}\tilde{E}]_{22}$
given by: 
\[
\begin{pmatrix}0 & 0 & 1 & 0 & 0\\
0 & 0 & 0 & F\veps E & 0\\
\eta & y_{1} & 0 & 0 & \sigma\\
0 & 1 & 0 & 0 & 0\\
0 & 0 & 0 & F\sigma E & 0
\end{pmatrix}
\]
 using the ordered decompositions from Eq.~\ref{eq:G1G1-ordered-decomp}
and Prop.~\ref{prop:various-as-A-A-bimod}: 
\begin{itemize}
\item $[\tilde{E}\tilde{F}]_{22}\cong G_{1}G_{1}\oplus EF[y]\cong A[y]\oplus FE[y]\oplus FE[y]\oplus FEFE[y]\oplus EF[y]$,
\item $[\tilde{F}\tilde{E}]_{22}\cong U\cong FE[y]^{\oplus4}\oplus F^{2}E^{2}[y]$.
\end{itemize}
We compute $[\tilde{\sigma}]_{22}$ first on $G_{1}G_{1}$, and afterwards
on $EF[y]$. We can use the same presentation for $[\tilde{\eta}(1)]_{22}$
as in the calculations for $[\tilde{\sigma}]_{21}$. Let $(\theta,\varphi_{1})\otimes(\theta',\varphi_{1}')\in G_{1}G_{1}$
be an arbitrary generator. Then the presentation for the product in
$[\tilde{F}]\cdot[\tilde{E}]\cdot[\tilde{E}]\cdot[\tilde{F}]$ is:
\begin{align*}
 & \phantom{+}\sum_{a\in Q}\left(\begin{smallmatrix}0 & 0\\
0 & (f_{a},0,0)
\end{smallmatrix}\right)\cdot\left(\begin{smallmatrix}0 & 0\\
0 & (e_{a},0,0)
\end{smallmatrix}\right)\cdot\left(\begin{smallmatrix}0 & 0\\
(\theta,\varphi_{1}) & 0
\end{smallmatrix}\right)\cdot\left(\begin{smallmatrix}0 & (\theta',\varphi_{1}')\\
0 & 0
\end{smallmatrix}\right)\\
 & +\sum_{b\in Q}\left(\begin{smallmatrix}0 & 0\\
0 & (0,f_{b},0)
\end{smallmatrix}\right)\cdot\left(\begin{smallmatrix}0 & 0\\
0 & (0,e_{b},0)
\end{smallmatrix}\right)\cdot\left(\begin{smallmatrix}0 & 0\\
(\theta,\varphi_{1}) & 0
\end{smallmatrix}\right)\cdot\left(\begin{smallmatrix}0 & (\theta',\varphi_{1}')\\
0 & 0
\end{smallmatrix}\right)\\
 & \in\left(\begin{smallmatrix}F[y] & L_{1}\\
F^{2}[y] & L_{2}
\end{smallmatrix}\right)\cdot\left(\begin{smallmatrix}E[y] & E^{2}[y]\\
G_{1} & G_{2}
\end{smallmatrix}\right)\cdot\left(\begin{smallmatrix}E[y] & E^{2}[y]\\
G_{1} & G_{2}
\end{smallmatrix}\right)\cdot\left(\begin{smallmatrix}F[y] & L_{1}\\
F^{2}[y] & L_{2}
\end{smallmatrix}\right).
\end{align*}
 Using again the calculations for $[\tilde{\sigma}]_{21}$, we see
that application of $[\tilde{F}\tilde{E}\tilde{\veps}]\circ[\tilde{F}\tilde{\tau}\tilde{F}]$
yields: 
\[
\begin{pmatrix}0 & 0\\
0 & (\varphi_{1},0,\theta,0,0)
\end{pmatrix}\cdot\begin{pmatrix}0 & 0\\
0 & (\theta',\varphi_{1}')
\end{pmatrix}\in[\tilde{F}\tilde{E}]\cdot[C].
\]
 Now compute the action of $(\theta',\varphi_{1}')$ on the right
on $U$ using Lemma \ref{lem:G-action-R}. For the matrix part $[\Phi]$,
we have: 
\[
\begin{pmatrix}\varphi' & 0\\
\varphi_{1}' & \theta'
\end{pmatrix}\cdot\begin{pmatrix}\varphi_{1} & \theta\\
0 & 0
\end{pmatrix}=\begin{pmatrix}\varphi'\circ\varphi_{1} & \theta.\varphi'\\
\varphi_{1}'\circ\varphi_{1} & \theta.\varphi_{1}'
\end{pmatrix}.
\]
 The submodule form of $(\varphi_{1},0,\theta,0,0)$ is $(\varphi_{1},0,\theta,0,\tau\ci E\varphi)$
using: 
\begin{align*}
\Lambda & =\tau y_{1}(E\varphi_{1}+E\theta\ci\tau)\\
 & =\tau\ci E(y_{1}\varphi_{1})+\tau\ci E\theta\\
 & =\tau\ci E\varphi.
\end{align*}
 Then after taking the action, the last coefficient of the submodule
form is given by post-composing with $E\varphi'$ to obtain $\Lambda=E\varphi'\circ\tau\circ E\varphi$,
which we expand using $\varphi=\_.\theta+y_{1}\varphi_{1}$ and similarly
for $\varphi'$: 
\begin{align*}
\Lambda & =E\varphi'\circ\tau\circ E\varphi\\
 & =E\theta'\circ\tau\circ E\theta+E\theta'\circ\tau\circ E(y_{1}\varphi_{1})\\
 & \quad+E(y_{1}\varphi_{1}')\circ\tau\circ E\theta+E(y_{1}\varphi_{1}')\circ\tau\circ E(y_{1}\varphi_{1}).
\end{align*}
 To compute the bimodule form, evaluate Eq.~(\ref{eq:Lambda^circ})
using $[\Phi]$: 
\begin{align*}
 & \Lambda=\\
 & \tau y_{1}\bigl(E(\varphi'\circ\varphi_{1})+E(\theta.\varphi')\circ\tau\bigr)-y_{2}\tau y_{1}\bigl(E(\varphi_{1}'\circ\varphi_{1})+E(\theta.\varphi_{1}')\circ\tau\bigr)+y_{1}y_{2}\Lambda^{\circ}\\
 & =\tau y_{1}\circ E\varphi'\circ\bigl(E\varphi_{1}+E\theta\circ\tau\bigr)-y_{2}\tau\circ E(y_{1}\varphi_{1}')\circ\bigl(E\varphi_{1}+E\theta\circ\tau\bigr)+y_{1}y_{2}\Lambda^{\circ}\\
 & =\tau y_{1}\circ E\theta'\circ\bigl(E\varphi_{1}+E\theta\circ\tau\bigr)+E(y_{1}\varphi_{1}')\circ\bigl(E\varphi_{1}+E\theta\circ\tau\bigr)+y_{1}y_{2}\Lambda^{\circ}.
\end{align*}
(For the last equality we expand $\varphi'=\_.\theta'+y_{1}\varphi_{1}'$
and use the relation $\tau y_{1}-y_{2}\tau=\idop$.) By identifying
the two expressions we can solve to find $\Lambda^{\ci}=E\varphi_{1}'\circ\tau\circ E\varphi_{1}$.
So the image is given using the bimodule form of $U$ by: 
\[
\begin{pmatrix}0 & 0\\
0 & \bigl(\varphi'\circ\varphi_{1},\varphi_{1}'\circ\varphi_{1},\theta.\varphi',\theta.\varphi_{1}',E\varphi_{1}'\circ\tau\circ E\varphi_{1}\bigr)
\end{pmatrix}\in\begin{pmatrix}G_{1} & G_{2}\\
L_{2} & U
\end{pmatrix}=[\tilde{F}\tilde{E}].
\]
 Using the fact that $E\varphi_{1}'\circ\tau\circ E\varphi_{1}=F\sigma E(\varphi_{1}\otimes\varphi_{1}')$
and $\varphi_{1}'\circ\varphi_{1}=F\veps E(\varphi_{1}\otimes\varphi_{1}')$,
one recovers the first four columns of the matrix of $[\tilde{\sigma}]_{22}$.

For the fifth column of $[\tilde{\sigma}]_{22}$, we start with an
arbitrary generator $e\otimes f'\in EF[y]\subset[\tilde{E}\tilde{F}]_{22}$.
The element $(0,e,0)\otimes(f',0,0)\in(G_{2})_{G}L_{2}$ is sent by
$\Gamma_{22}$ of $[\tilde{E}\tilde{F}]$ to $e\otimes f'$. So we
consider the element: 
\begin{align*}
 & \phantom{+}\sum_{a\in Q}\left(\begin{smallmatrix}0 & 0\\
0 & (f_{a},0,0)
\end{smallmatrix}\right)\cdot\left(\begin{smallmatrix}0 & 0\\
0 & (e_{a},0,0)
\end{smallmatrix}\right)\cdot\left(\begin{smallmatrix}0 & 0\\
0 & (0,e,0)
\end{smallmatrix}\right)\cdot\left(\begin{smallmatrix}0 & 0\\
0 & (f',0,0)
\end{smallmatrix}\right)\\
 & +\sum_{b\in Q}\left(\begin{smallmatrix}0 & 0\\
0 & (0,f_{b},0)
\end{smallmatrix}\right)\cdot\left(\begin{smallmatrix}0 & 0\\
0 & (0,e_{b},0)
\end{smallmatrix}\right)\cdot\left(\begin{smallmatrix}0 & 0\\
0 & (0,e,0)
\end{smallmatrix}\right)\cdot\left(\begin{smallmatrix}0 & 0\\
0 & (f',0,0)
\end{smallmatrix}\right)\\
 & \in\left(\begin{smallmatrix}F[y] & L_{1}\\
F^{2}[y] & L_{2}
\end{smallmatrix}\right)\cdot\left(\begin{smallmatrix}E[y] & E^{2}[y]\\
G_{1} & G_{2}
\end{smallmatrix}\right)\cdot\left(\begin{smallmatrix}E[y] & E^{2}[y]\\
G_{1} & G_{2}
\end{smallmatrix}\right)\cdot\left(\begin{smallmatrix}F[y] & L_{1}\\
F^{2}[y] & L_{2}
\end{smallmatrix}\right),
\end{align*}
 and we compute its image under $\tilde{F}\tilde{E}\tilde{\veps}\ci\tilde{F}\tilde{\tau}\tilde{F}$.
First apply $\Gamma_{22}$ of $[\tilde{E}\tilde{E}]$ to $(e_{a},0,0)\otimes(0,e,0)$
and $(0,e_{b},0)\otimes(0,e,0)$, using the rule for bimodule forms
on p.~\pageref{Gamma-22-EE-bimods}: 
\begin{align*}
 & (e_{a},0,0)\otimes(0,e,0)\overset{\Gamma_{22}}{\longmapsto}(0,-y_{2}(e_{a}\otimes e),-y_{2}(e_{a}\otimes e),0)\in G_{3},\\
 & (0,e_{b},0)\otimes(0,e,0)\overset{\Gamma_{22}}{\longmapsto}(\tau y_{1}(e_{b}\otimes e),e_{b}\otimes e,e_{b}\otimes e,0)\in G_{3}.
\end{align*}
 Next we apply $[\tilde{\tau}]_{22}$ to these elements: 
\begin{align*}
 & (0,-y_{2}(e_{a}\otimes e),-y_{2}(e_{a}\otimes e),0)\overset{[\tilde{\tau}]_{22}}{\longmapsto}\bigl(e_{a}\otimes e,e_{a}\otimes e,-\tau y_{2}(e_{a}\otimes e),0\bigr),\\
 & (\tau y_{1}(e_{b}\otimes e),e_{b}\otimes e,e_{b}\otimes e,0)\overset{[\tilde{\tau}]_{22}}{\longmapsto}\bigl(\tau(e_{b}\otimes e),\tau(e_{b}\otimes e),\tau(e_{b}\otimes e),0\bigr).
\end{align*}
 Note that formula (\ref{eq:Def-tilde-tau}) is given for the submodule
form of $G_{3}$. Using Prop.~3.21 of \cite{mcmillanTensor2product2representations2022},
one defines a bimodule form in the usual way, where the last coefficient
is $\chi''$ instead of $\chi$. By studying the proof of Lemma 4.3
of \cite{mcmillanTensor2product2representations2022}, one observes
that the action of $\tilde{\tau}$ on the last coefficient in this
bimodule form is (also) given by post-composition with $\tau E$,
whence the final zeros above.

The next step is to express $\bigl(e_{a}e,e_{a}e,-\tau y_{2}(e_{a}e),0\bigr)$
and $\bigl(\tau(e_{b}e),\tau(e_{b}e),\tau(e_{b}e),0\bigr)$ back in
$(G_{2})_{G}G_{2}$ (i.e.~find a preimage under $\Gamma_{22}|_{(G_{2})_{G}G_{2}}$)
in order to view them in $[\tilde{E}]\cdot[\tilde{E}]$. We will need
the notation $\tau(ee)=\sum_{d\in P}\tau(ee)_{(1d)}\otimes\tau(ee)_{(2d)}$
introduced to compute $[\tilde{\sigma}_{21}]$ above.
\begin{claim*}
We have: 
\begin{align*}
\sum_{d\in P}\substack{\bigl(0,\tau(e_{a}e)_{(1d)},0\bigr)\otimes\bigl(\tau(e_{a}e)_{(2d)},y_{1}\tau(e_{a}e)_{(2d)},0\bigr)\\
-\bigl(0,\tau y_{1}(e_{a}e)_{(1d)},0\bigr)\otimes\bigl(0,\tau y_{1}(e_{a}e)_{(2d)},0\bigr)
}
 & \overset{\Gamma_{22}}{\longmapsto}\bigl(e_{a}e,e_{a}e,-\tau y_{2}(e_{a}e),0\bigr),\\
\sum_{d\in P}\bigl(0,\tau(e_{b}e)_{(1d)},0\bigr)\otimes\bigl(0,\tau(e_{b}e)_{(2d)},0\bigr) & \overset{\Gamma_{22}}{\longmapsto}\bigl(\tau(e_{b}e),\tau(e_{b}e),\tau(e_{b}e),0\bigr).
\end{align*}
\end{claim*}
\begin{proof}
The proof is a direct calculation using the bimodules formulation
of $\Gamma_{22}|_{(G_{2})_{G}G_{2}}$ on p.~\pageref{Gamma-22-EE-bimods}. 
\end{proof}
Thus, after applying $\tilde{F}\tilde{\tau}\tilde{F}$, we have the
element: 
\begin{align*}
 & \phantom{+}\sum_{a\in Q,d\in P}\left(\begin{smallmatrix}0 & 0\\
0 & (f_{a},0,0)
\end{smallmatrix}\right)\cdot\left(\begin{smallmatrix}0 & 0\\
0 & \left(0,\tau(e_{a}e)_{(1d)},0\right)
\end{smallmatrix}\right)\cdot\left(\begin{smallmatrix}0 & 0\\
0 & \left(\tau(e_{a}e)_{(2d)},y_{1}\tau(e_{a}e)_{(2d)},0\right)
\end{smallmatrix}\right)\cdot\left(\begin{smallmatrix}0 & 0\\
0 & (f',0,0)
\end{smallmatrix}\right)\\
 & +\sum_{a\in Q,d\in P}\left(\begin{smallmatrix}0 & 0\\
0 & (f_{a},0,0)
\end{smallmatrix}\right)\cdot\left(\begin{smallmatrix}0 & 0\\
0 & -\left(0,\tau y_{1}(e_{a}e)_{(1d)},0\right)
\end{smallmatrix}\right)\cdot\left(\begin{smallmatrix}0 & 0\\
0 & \left(0,\tau y_{1}(e_{a}e)_{(2d)},0\right)
\end{smallmatrix}\right)\cdot\left(\begin{smallmatrix}0 & 0\\
0 & (f',0,0)
\end{smallmatrix}\right)\\
 & +\sum_{b\in Q,d\in P}\left(\begin{smallmatrix}0 & 0\\
0 & (0,f_{b},0)
\end{smallmatrix}\right)\cdot\left(\begin{smallmatrix}0 & 0\\
0 & \left(0,\tau(e_{b}e)_{(1d)},0\right)
\end{smallmatrix}\right)\cdot\left(\begin{smallmatrix}0 & 0\\
0 & \left(0,\tau(e_{b}e)_{(2d)},0\right)
\end{smallmatrix}\right)\cdot\left(\begin{smallmatrix}0 & 0\\
0 & (f',0,0)
\end{smallmatrix}\right),
\end{align*}
 and we need to apply $[\tilde{F}\tilde{E}]\cdot\tilde{\veps}$ and
then realize the result in $[\tilde{F}\tilde{E}]$. Observe that:
\begin{align*}
\bigl(0,\tau y_{1}(e_{a}e)_{(2d)},0\bigr)\otimes(f',0,0) & \overset{\tilde{\veps}}{\mapsto}0,\\
\bigl(0,\tau(e_{b}e)_{(2d)},0\bigr)\otimes(f',0,0) & \overset{\tilde{\veps}}{\mapsto}0.
\end{align*}
 Therefore only the top row will remain. We have in submodule form:
\begin{multline*}
\bigl(\tau(e_{a}e)_{(2d)},y_{1}\tau(e_{a}e)_{(2d)},0\bigr)\otimes(f',0,0)\\
\overset{[\tilde{\veps}]_{22}}{\longmapsto}\bigl(f'(\tau(e_{a}e)_{(2d)}),Ef'\circ\tau\circ\bigl(\_\otimes y_{1}\tau(e_{a}e)_{(2d)}\bigr)\bigr)\in G_{1}.
\end{multline*}
 We convert to bimodule form and give this a name: 
\[
(\theta,\varphi_{1})_{a,d}:=\bigl(f'(\tau(e_{a}e)_{(2d)}),Ef'\circ\tau(\_\otimes\tau(e_{a}e)_{(2d)})\bigr)\in G_{1}.
\]
 Observe that under the composition isomorphism $(L_{2})_{G}G_{2}\iso U$
we have: 
\[
(f_{a},0,0)\otimes\bigl(0,\tau(e_{a}e)_{(1d)},0\bigr)\mapsto\bigl(0,0,f_{a}(\_).\tau(e_{a}e){}_{(1d)},0,0\bigr)\in U.
\]
 We are therefore left with: 
\[
\sum_{a\in Q,d\in P}\begin{pmatrix}0 & 0\\
0 & \bigl(0,0,f_{a}(\_).\tau(e_{a}e){}_{(1d)},0,0\bigr)
\end{pmatrix}\cdot\begin{pmatrix}0 & 0\\
0 & (\theta,\varphi_{1})_{a,d}
\end{pmatrix}\in[\tilde{F}\tilde{E}]\cdot[C].
\]

It remains to use the right action of $G_{1}^{\op}$ on $R$ (Lemma
\ref{lem:G-action-R}) to compute the action of $(\theta,\varphi_{1})_{a,d}$.
The new matrix is given for each term of the sum by: 
\begin{gather*}
\begin{pmatrix}Ef'\circ\tau y_{1}\bigl(\_\otimes\tau(e_{a}e)_{(2d)}\bigr) & 0\\
Ef'\circ\tau(\_\otimes\tau(e_{a}e)_{(2d)}) & f'(\tau(e_{a}e)_{(2d)})
\end{pmatrix}\cdot\begin{pmatrix}0 & f_{a}(\_).\tau(e_{a}e){}_{(1d)}\\
0 & 0
\end{pmatrix}\\
=\begin{pmatrix}0 & Ef'\circ\tau y_{1}\Bigl(f_{a}(\_).\tau(e_{a}e){}_{(1d)}\otimes\tau(e_{a}e)_{(2d)}\Bigr)\\
0 & Ef'\circ\tau\Bigl(f_{a}(\_).\tau(e_{a}e){}_{(1d)}\otimes\tau(e_{a}e)_{(2d)}\Bigr)
\end{pmatrix}.
\end{gather*}
 After summing over $a$ and $d$ this becomes: 
\[
\overset{\sum_{a,d}}{\longsquiggly}\begin{pmatrix}0 & Ef'\circ\tau\circ y_{1}\tau(\_\otimes e)\\
0 & Ef'\circ\tau\bigl(\tau(\_\otimes e)\bigr)
\end{pmatrix}=\begin{pmatrix}0 & Ef'\circ\tau(\_\otimes e)\\
0 & 0
\end{pmatrix}.
\]
 This matrix gives the first four components of the final element
of $U$. To find the fifth in submodule form, we compute the submodule
form of $\bigl(0,0,f_{a}(\_).\tau(e_{a}e){}_{(1d)},0,0\bigr)$ and
post-compose with $E\varphi$: 
\begin{align*}
 & E\varphi\ci\Bigl(\tau y_{1}\ci E\bigl(f_{a}(\_).\tau(e_{a}e){}_{(1d)}\bigr)\ci\tau\Bigr)\\
 & =E\varphi\ci\Bigl(\tau y_{1}\bigl(Ef_{a}\circ\tau(\_\:\_)\otimes\tau(e_{a}e)_{(1d)}\bigr)\Bigr)\\
 & =\Bigl(E^{2}f'\circ E\tau\ci y_{1}\bigl(\_\:\_\otimes\tau(e_{a}e)_{(2d)}\bigr)\Bigr)\ci\Bigl(\tau y_{1}\bigl(Ef_{a}\circ\tau(\_\:\_)\otimes\tau(e_{a}e)_{(1d)}\bigr)\Bigr)\\
 & =E^{2}f'\circ E\tau\ci y_{1}\bigl(\tau y_{1}\bigl(Ef_{a}\circ\tau(\_\:\_)\otimes\tau(e_{a}e)_{(1d)}\bigr)\otimes\tau(e_{a}e)_{(2d)}\bigr)\\
 & =E^{2}f'\circ E\tau\circ\tau E\circ y_{2}y_{1}\bigl(Ef_{a}\circ\tau(\_\:\_)\otimes\tau(e_{a}e)_{(1d)}\otimes\tau(e_{a}e)_{(2d)}\bigr).
\end{align*}
 (The last equality is, schematically, $y_{1}\bigl((\tau y_{1}AA)\otimes B\bigr)=\tau E\ci y_{2}y_{1}(AA\otimes B)$.)
Summing over $d$ and $a$ we obtain: 
\[
\overset{\sum_{a,d}}{\longsquiggly}E^{2}f'\circ E\tau\circ\tau E\circ y_{2}y_{1}\bigl(E\tau\circ\tau E(\_\:\_\otimes e)\bigr).
\]
 Now observe the following calculation in the nil affine Hecke algebra:
\begin{align*}
\tau_{1}(\tau_{2}y_{2})y_{1}\tau_{1}\tau_{2} & =\tau_{1}(y_{3}\tau_{2})y_{1}\tau_{1}\tau_{2}+\tau_{1}y_{1}\tau_{1}\tau_{2}\\
 & =(\tau_{1}y_{3})(\tau_{2}y_{1})\tau_{1}\tau_{2}+(\tau_{1}y_{1})\tau_{1}\tau_{2}\\
 & =(y_{3}\tau_{1})(y_{1}\tau_{2})\tau_{1}\tau_{2}+(y_{2}\tau_{1})\tau_{1}\tau_{2}+\tau_{1}\tau_{2}\\
 & =y_{3}(\tau_{1}y_{1})\tau_{2}\tau_{1}\tau_{2}+0+\tau_{1}\tau_{2}\\
 & =y_{3}(y_{2}\tau_{1})\tau_{2}\tau_{1}\tau_{2}+y_{3}\tau_{2}\tau_{1}\tau_{2}+\tau_{1}\tau_{2}\\
 & =0+y_{3}\tau_{2}\tau_{1}\tau_{2}+\tau_{1}\tau_{2}.
\end{align*}
(Here $\tau_{i}=E^{n-i-1}\tau E^{i-1}$ for whatever $n$.) Therefore
we have: 
\begin{align*}
 & E^{2}f'\circ E\tau\circ\tau E\circ y_{2}y_{1}\bigl(E\tau\circ\tau E(\_\:\_\otimes e)\bigr)\\
 & =y_{2}E^{2}f'\circ\tau E\circ E\tau\circ\tau E(\_\:\_\otimes e)+E^{2}f'\circ E\tau\circ\tau E(\_\:\_\otimes e).
\end{align*}
 Now to find the bimodule form of the fifth component we consider:
\begin{align*}
 & \tau y_{1}\circ\bigl(E^{2}f'\circ E\tau(\_\:\_\otimes e)\circ\tau\bigr)\\
 & =\tau y_{1}\circ E^{2}f'\circ E\tau\circ\tau E(\_\:\_\otimes e)\\
 & =y_{2}E^{2}f'\circ\tau E\circ E\tau\circ\tau E(\_\:\_\otimes e)+E^{2}f'\circ E\tau\circ\tau E(\_\:\_\otimes e),
\end{align*}
 and since this agrees with the expression before it, Eq.~(\ref{eq:Lambda^circ})
implies that the fifth component in bimodule form is zero. The final
expression is $\bigl(0,0,Ef'\circ\tau\bigl(\_\otimes e\bigr),0,0\bigr)\in U\cong[\tilde{F}\tilde{E}]_{22}.$
Observe that $Ef'\circ\tau(\_\otimes e)=\sigma(e\otimes f')$. This
gives the fifth column of the matrix of $[\tilde{\sigma}]_{22}$,
and we have now justified all components of that matrix.
\end{itemize}

\subsubsection{Maps $\tilde{\protect\veps}\circ\tilde{x}^{i}\tilde{F}$ and $\tilde{F}\tilde{x}^{i}\circ\tilde{\eta}$}

We continue by computing the maps $\tilde{\veps}\circ\tilde{x}^{i}\tilde{F}$
and $\tilde{F}\tilde{x}^{i}\circ\tilde{\eta}$ on the various components
of the matrices $[\tilde{E}\tilde{F}]$, $[\tilde{F}\tilde{E}]$,
and $[C]$. As before, we propose these maps in the bulleted lines
and justify them in the paragraphs following.
\begin{itemize}
\item We have $[\tilde{\veps}\circ\tilde{x}^{i}\tilde{F}]_{11}:[\tilde{E}\tilde{F}]_{11}\to[C]_{11}$
given by $\veps\circ x^{i}y_{1}F$ using the decompositions: 
\begin{itemize}
\item $[\tilde{E}\tilde{F}]_{11}\cong EF[y]$,
\item $[C]_{11}\cong A[y]$.
\end{itemize}
The endomorphism $\tilde{x}\in\End(\tilde{E})$ (see (\ref{eq:Def-tilde-x}))
determines an endomorphism of $[\tilde{E}\tilde{F}]_{11}$ given by
$xF$ on $EF[y]$. The morphism $\tilde{\veps}$ composes elements
of $\tilde{E}$ with those of $\tilde{F}$ when they are interpreted
in $\Hom_{D^{b}(B)}(X,E'X)$ and $\Hom_{D^{b}(B)}(E'X,X)$. In particular,
$e\in E[y]\cong[\tilde{E}]_{11}$ represents the morphism $X_{1}\to E'X_{1}$
given by $1\mapsto y_{1}e$ in degree $0$ of the top row, and $f\in F[y]\cong[\tilde{F}]_{11}$
represents the morphism given by $e\mapsto f(e)$ in degree $0$ of
the top row.
\item We have $[\tilde{F}\tilde{x}^{i}\circ\tilde{\eta}]_{11}:[C]_{11}\to[\tilde{F}\tilde{E}]_{11}$
given by $\left(\begin{smallmatrix}y^{i}\\
Fh_{i-1}(x,y)\circ\eta
\end{smallmatrix}\right)$ using the decompositions: 
\begin{itemize}
\item $[C]_{11}\cong A[y]$,
\item $[\tilde{F}\tilde{E}]_{11}\cong G_{1}\cong A[y]\oplus FE[y]$.
\end{itemize}
Here $h_{i}(z_{1},\dots,z_{n})$ is the complete homogeneous symmetric
polynomial of degree $i$ in the variables $z_{1},\dots,z_{n}$. Note
the small case interpretations: 
\[
\begin{cases}
h_{i-1}(x,y)=0 & i=0\\
h_{i-1}(x,y)=1 & i=1\\
h_{i-1}(x,y)=x+y & i=2\\
\dots & \dots
\end{cases}
\]

Observe that $[\tilde{\eta}]_{11}$ is given by $1\mapsto\idop_{X_{2}}\in G_{1}^{\op}\cong\End_{K^{b}(B)}(X_{2})$,
and $\idop_{X_{2}}=(1,0)$ (in bimodule form). More generally $\theta\mapsto\_.\theta\in\Hom_{A}(_{A}E,E)[y]\cong FE[y]\subset G_{1}^{\op}$.
From (\ref{eq:Def-tilde-x}) we have the action of $[\tilde{x}]_{11}$
on $G_{1}^{\op}$ in submodule form: $\tilde{x}^{i}.(\theta,\varphi)=(y^{i}\theta,x^{i}\circ\varphi)$.
Now convert this expression to bimodule form: 
\begin{align*}
x^{i}\circ\varphi & =x^{i}\ci\_.\theta+x^{i}y_{1}\varphi_{1}\\
 & =y^{i}\theta+(x^{i}-y^{i})\ci\_.\theta+y_{1}x^{i}\varphi_{1}\\
 & =y^{i}\theta+y_{1}\bigl(h_{i-1}(x,y)\ci\_.\theta+x^{i}\circ\varphi_{1}\bigr),
\end{align*}
 so $\tilde{x}^{i}.(\theta,\varphi_{1})=(y^{i}\theta,h_{i-1}(x,y)\ci\_.\theta+x^{i}\circ\varphi_{1})$.
In particular, $\tilde{x}^{i}.(1,0)=(y^{i},h_{i-1}(x,y))$, which
gives the proposed formula by viewing $x,y$ as endofunctors of $E$
instead of as elements of $FE[y]$.
\item We have $[\tilde{\veps}\circ\tilde{x}^{i}\tilde{F}]_{21}:[\tilde{E}\tilde{F}]_{21}\to[C]_{21}$
given by $\bigl(x^{i},F(\veps\circ x^{i}y_{1}F)\bigr)$ using the
decompositions: 
\begin{itemize}
\item $[\tilde{E}\tilde{F}]_{21}\cong G_{1}F[y]\cong F[y]\oplus FEF[y]$,
\item $[C]_{21}\cong F[y]$.
\end{itemize}
(Here $x\in\End(F)[y]$ is given by $x(f)=f\circ x$.) The map $[\tilde{\veps}]_{21}:G_{1}F[y]\to F[y]$
is given (using submodule form) by $(\theta,\varphi)\otimes f\mapsto f\ci\varphi$.
The endomorphism $[\tilde{x}]_{21}$ acts on $G_{1}$ as described
under the previous bullet: $\tilde{x}^{i}.(\theta,\varphi_{1})=(y^{i}\theta,h_{i-1}(x,y)\ci\_.\theta+x^{i}\circ\varphi_{1})$.
Then $[\tilde{\veps}\circ\tilde{x}^{i}\tilde{F}]_{21}:G_{1}F[y]\to F[y]$
is given using bimodule form by: 
\begin{align*}
\tilde{x}^{i}.(\theta,\varphi_{1})\otimes f & \mapsto f\ci x^{i}\ci\varphi\\
 & \quad=f\circ x^{i}\ci\_.\theta+f\circ x^{i}y_{1}\varphi_{1},
\end{align*}
 and the component data follows from this formula.
\item We have $[\tilde{F}\tilde{x}^{i}\circ\tilde{\eta}]_{21}:[C]_{21}\to[\tilde{F}\tilde{E}]_{21}$
given by $\left(\begin{smallmatrix}0\\
y^{i}\\
F(Fh_{i-1}(x,y)\circ\eta)
\end{smallmatrix}\right)$ using the decompositions: 
\begin{itemize}
\item $[C]_{21}\cong F[y]$,
\item $[\tilde{F}\tilde{E}]_{21}\cong L_{2}\cong F[y]\oplus F[y]\oplus F^{2}E[y]$.
\end{itemize}
Let $\left(\begin{smallmatrix}0 & 0\\
f & 0
\end{smallmatrix}\right)\in\left(\begin{smallmatrix}A[y] & E[y]\\
F[y] & G_{1}^{\op}
\end{smallmatrix}\right)=[C]$, and observe that: 
\begin{alignat*}{2}
\tilde{\eta}\left(\left(\begin{smallmatrix}0 & 0\\
f & 0
\end{smallmatrix}\right)\right) & =\tilde{\eta}\left(\left(\begin{smallmatrix}0 & 0\\
f & 0
\end{smallmatrix}\right).\left(\begin{smallmatrix}1 & 0\\
0 & 0
\end{smallmatrix}\right)\right) &  & =\left(\begin{smallmatrix}0 & 0\\
f & 0
\end{smallmatrix}\right).\tilde{\eta}\left(\left(\begin{smallmatrix}1 & 0\\
0 & 0
\end{smallmatrix}\right)\right)\\
 & =\left(\begin{smallmatrix}0 & 0\\
f & 0
\end{smallmatrix}\right).\left(\begin{smallmatrix}(1,0) & 0\\
0 & 0
\end{smallmatrix}\right) &  & =\left(\begin{smallmatrix}0 & 0\\
(0,f,0) & 0
\end{smallmatrix}\right)\in\left(\begin{smallmatrix}G_{1} & G_{2}\\
L_{2} & U
\end{smallmatrix}\right)=[\tilde{F}\tilde{E}].
\end{alignat*}
 Here $(0,f,0)$ is written in the bimodule form of $L_{2}$. (The
action of $f\in F[y]\subset[C]_{21}$ on generators in $G_{1}\subset[\tilde{F}\tilde{E}]$
is given by $F[y]G_{1}\to L_{2}$, $f\otimes(\theta,\varphi)\mapsto(0,f\ci\_.\theta,\varphi\ci Ef)$
(written in submodule form), and this image is $(0,f\ci\_.\theta,\varphi_{1}\ci Ef)$
in bimodule form.)

Now we apply $[\tilde{F}\tilde{x}^{i}]_{21}$. From Eq.~(\ref{eq:FE-expansion}):
\[
[\tilde{F}]\cdot[\tilde{E}]\supset(L_{2})_{G}G_{1}\ni(0,f,0)\otimes(1,0)\overset{\Gamma_{21}}{\longmapsto}(0,f,0)\in L_{2}\subset[\tilde{F}\tilde{E}].
\]
 We have already seen that $\tilde{x}^{i}.(1,0)=(y^{i},h_{i-1}(x,y))\in G_{1}$,
so we have: 
\[
(0,f,0)\otimes(1,0)\overset{[\tilde{F}\tilde{x}^{i}]_{21}}{\longmapsto}(0,f,0)\otimes(y^{i},h_{i-1}(x,y)).
\]
Then 
\[
\Gamma_{21}:(0,f,0)\otimes(y^{i},h_{i-1}(x,y))\mapsto\bigl(0,y^{i}f,x^{i}\circ Ef\bigr)
\]
 written in submodule form. In bimodule form the image is: 
\[
\bigl(0,y^{i}f,h_{i-1}(x,y)\circ Ef\bigr),
\]
 which we compute using: 
\begin{align*}
x^{i}\circ Ef & =\bigl(y^{i}+y_{1}h_{i-1}(x,y)\bigr)\circ Ef\\
 & =E(y^{i}f)+y_{1}\bigl(h_{i-1}(x,y)\circ Ef\bigr).
\end{align*}
 Note that $F^{2}E[y]\ni h_{i-1}(x,y)\circ Ef=F\bigl(h_{i-1}(x,y)\circ\eta\bigr)(f)$.
\item We have $[\tilde{\veps}\circ\tilde{x}^{i}\tilde{F}]_{12}:[\tilde{E}\tilde{F}]_{12}\to[C]_{12}$
given by $\bigl(x^{i},(\veps\circ x^{i}y_{1}F)E\bigr)$ using the
decompositions: 
\begin{itemize}
\item $[\tilde{E}\tilde{F}]_{12}\cong E[y]G_{1}\cong E[y]\oplus EFE[y]$,
\item $[C]_{12}\cong E[y]$.
\end{itemize}
The endomorphism $[\tilde{x}]_{12}$ acts as $x$ on $E[y]=[\tilde{E}]_{11}$,
and thus as $xG_{1}$ on $E[y]G_{1}=[\tilde{E}\tilde{F}]_{12}$. The
map $[\tilde{\veps}]_{12}:E[y]G_{1}\to E[y]$ is given (using submodule
form) by $e\otimes(\theta,\varphi)\mapsto y_{1}^{-1}\varphi(y_{1}e)$.
(Recall that $e\in E[y]$ indicates the map $X_{1}\to X_{2}$ given
on the top row by $A[y]\to E[y]$, $1\mapsto y_{1}e$.) So we have:
\[
x^{i}(e)\otimes(\theta,\varphi_{1})\overset{[\tilde{\veps}]_{12}}{\longmapsto}y_{1}^{-1}\varphi(x^{i}y_{1}e)=x^{i}(e).\theta+\varphi_{1}(x^{i}y_{1}e),
\]
 and the component data follows from this formula.
\item We have $[\tilde{F}\tilde{x}^{i}\circ\tilde{\eta}]_{12}:[C]_{12}\to[\tilde{F}\tilde{E}]_{12}$
given by $\left(\begin{smallmatrix}y^{i}\\
y^{i}y_{1}\\
(Fh_{i-1}(x,y)\circ\eta)E
\end{smallmatrix}\right)$ using the decompositions: 
\begin{itemize}
\item $[C]_{12}\cong E[y]$,
\item $[\tilde{F}\tilde{E}]_{12}\cong G_{2}\cong E[y]\oplus E[y]\oplus FE^{2}[y]$.
\end{itemize}
By reasoning as in the $[\tilde{F}\tilde{x}^{i}\circ\tilde{\eta}]_{21}$
case, we find: 
\[
[C]\ni\left(\begin{smallmatrix}0 & e\\
0 & 0
\end{smallmatrix}\right)\overset{[\tilde{\eta}]}{\longmapsto}\left(\begin{smallmatrix}0 & (e,y_{1}e,0)\\
0 & 0
\end{smallmatrix}\right)\in\left(\begin{smallmatrix}G_{1} & G_{2}\\
L_{2} & U
\end{smallmatrix}\right)=[\tilde{F}\tilde{E}],
\]
 using the bimodule form of $G_{2}$. Now we apply $[\tilde{F}\tilde{x}^{i}]_{12}$.
From Eq.~(\ref{eq:FE-expansion}): 
\[
[\tilde{F}]\cdot[\tilde{E}]\supset(L_{1})_{G}G_{2}\ni(1,0)\otimes(e,y_{1}e,0)\overset{\Gamma_{12}}{\longmapsto}(e,y_{1}e,0)\in G_{2}\subset[\tilde{F}\tilde{E}].
\]
 In (\ref{eq:Def-tilde-x}) we have a formula for the action of $[\tilde{x}^{i}]_{22}$
on $G_{2}\subset[\tilde{E}]$ written in terms of the data $e_{1}$,
$e_{2}$, $\xi$. The data $(e,y_{1}e,0)$ corresponds to $e_{1}=y_{1}e$,
$e_{2}=0$, $\xi=\_\otimes y_{1}e$ (see the paragraph after Prop.~\ref{prop:G_2'-homs}).
Applying $[\tilde{x}^{i}]_{22}$ gives $e_{1}=y^{i}y_{1}e$, $e_{2}=0$,
$\xi=\_\otimes y^{i}y_{1}e+y_{1}y_{2}h_{i-1}(x_{2},y)(\_\otimes e)$,
where to compute $\xi$ we have used: 
\begin{align*}
x_{2}^{i}\circ(\_\otimes y_{1}e) & =\bigl(y^{i}+y_{2}h_{i-1}(x_{2},y)\bigr)\circ(\_\otimes y_{1}e)\\
 & =\_\otimes y^{i}y_{1}e+y_{1}y_{2}h_{i-1}(x_{2},y)(\_\otimes e).
\end{align*}
 This corresponds to the data $\bigl(y^{i}e,y^{i}y_{1}e,h_{i-1}(x_{2},y)(\_\otimes e)\bigr)\in G_{2}$
in the bimodule form. So we have: 
\begin{align*}
(1,0)\otimes(e,y_{1}e,0) & \overset{[\tilde{F}\tilde{x}^{i}]_{12}}{\longmapsto}(1,0)\otimes\bigl(y^{i}e,y^{i}y_{1}e,h_{i-1}(x_{2},y)(\_\otimes e)\bigr)\\
 & \overset{\Gamma_{12}}{\longmapsto}\bigl(y^{i}e,y^{i}y_{1}e,h_{i-1}(x_{2},y)(\_\otimes e)\bigr)\in G_{2}\subset[\tilde{F}\tilde{E}].
\end{align*}

Note that $FE^{2}[y]\ni h_{i-1}(x_{2},y)(\_\otimes e)=\bigl((Fh_{i-1}(x,y)\circ\eta)E\bigr)(e)$.
\item We have $[\tilde{\veps}\circ\tilde{x}^{i}\tilde{F}]_{22}:[\tilde{E}\tilde{F}]_{22}\to[C]_{22}$
given by: 
\[
\begin{pmatrix}y^{i} & 0 & 0 & 0 & -\veps\circ h_{i-1}(x,y)F\\
h_{i-1}(x,y)\circ\eta & x^{i}E & Fx^{i} & F(\veps\circ x^{i}y_{1}F)E & \substack{-FE\veps\circ F(\tau\circ h_{i-1}(x_{1},x_{2}))F\circ\eta EF\\
-FE\veps\circ F(h_{i-2}(x_{1},x_{2},y))F\circ\eta EF
}
\end{pmatrix}
\]
 using the ordered decompositions (recall Eq.~(\ref{eq:G1G1-ordered-decomp})): 
\begin{itemize}
\item $[\tilde{E}\tilde{F}]_{22}\cong G_{1}G_{1}\oplus EF[y]\cong A[y]\oplus FE[y]\oplus FE[y]\oplus FEFE[y]\oplus EF[y]$,
\item $[C]_{22}\cong G_{1}\cong A[y]\oplus FE[y]$.
\end{itemize}
Consider the first four columns first, i.e.~the restriction of the
map to $G_{1}G_{1}$. Take an arbitrary generator $(\theta,\varphi_{1})\otimes(\theta',\varphi_{1}')$.
Borrowing a calculation from the case ${[\tilde{\veps}\ci\tilde{x}^{i}\tilde{F}]_{21}}$
we find: 
\[
(\theta,\varphi_{1})\otimes(\theta',\varphi_{1}')\overset{[\tilde{x}^{i}\tilde{F}]_{22}}{\longmapsto}(y^{i}\theta,h_{i-1}(x,y)\ci\_.\theta+x^{i}\circ\varphi_{1})\otimes(\theta',\varphi_{1}').
\]
 Now $[\tilde{\veps}]_{22}:G_{1}G_{1}\to G_{1}$ is given by composition,
so we have: 
\begin{multline*}
(y^{i}\theta,h_{i-1}(x,y)\ci\_.\theta+x^{i}\circ\varphi_{1})\otimes(\theta',\varphi_{1}')\\
\overset{[\tilde{\veps}]_{22}}{\longmapsto}\bigl(y^{i}\theta\theta',\_.\theta'\circ h_{i-1}(x,y)\ci\_.\theta+(\_.\theta')\circ x^{i}\circ\varphi_{1}+\varphi_{1}'\circ(\_.y^{i}\theta)\\
\quad+\varphi_{1}'\circ(x^{i}-y^{i})\ci\_.\theta+\varphi_{1}'\circ y_{1}x^{i}\circ\varphi_{1}\bigr)\\
=\bigl(y^{i}\theta\theta',h_{i-1}(x,y)\ci\_.\theta\theta'+\varphi_{1}'\circ x^{i}\ci\_.\theta+x^{i}\ci\_.\theta'\circ\varphi_{1}+\varphi_{1}'\circ y_{1}x^{i}\circ\varphi_{1}\bigr).
\end{multline*}
 The first four columns of the matrix of $[\tilde{\veps}\circ\tilde{x}^{i}\tilde{F}]_{22}$
can be read off this formula.

The last column gives the restriction of $[\tilde{\veps}\circ\tilde{x}^{i}\tilde{F}]_{22}$
to a map $EF[y]\to A[y]\oplus FE[y]$. Its computation is more involved.
We start with a generator $e\otimes f$, and note that: 
\[
[\tilde{E}]\cdot[\tilde{F}]\supset(G_{2})_{G}L_{2}\ni(0,e,0)\otimes(f,0,0)\overset{\Gamma_{22}}{\longmapsto}e\otimes f\in EF[y]\subset[\tilde{E}\tilde{F}]_{22}
\]
using $\Gamma_{22}\mid_{(G_{2})_{G}L_{2}}=\kappa$ from Eq.~(\ref{eq:EF-expansion}).
Now we must apply $[\tilde{x}]_{22}$ to the first factor, and then
compose the factors, thereby applying $[\tilde{\veps}]_{22}$ and
giving an element of $G_{2}\cong\End_{K^{b}(B)}(X_{2})$.

The data $(0,e,0)$ corresponds to $e_{1}=e_{2}=e$, $\xi=\tau y_{1}(\_\otimes e)$
(see the paragraph after Prop.~\ref{prop:G_2'-homs}). The action
of $[\tilde{x}^{i}]_{22}$ on $G_{2}\subset[\tilde{E}]$ then gives
$e_{1}=y^{i}e$, $e_{2}=x^{i}e$, $\xi=x_{2}^{i}\ci\tau y_{1}(\_\otimes e)$.
We can compute the composite with $(f,0,0)$ directly using this information.
It is given in submodule form by: 
\begin{align*}
 & \bigl(f\ci y_{1}^{-1}(y^{i}e-x^{i}e),Ef\circ\tau\circ x_{2}^{i}\circ\tau y_{1}(\_\otimes e)\bigr)\\
 & =\bigl(f(-h_{i-1}(x,y)e),Ef\circ\tau\circ x_{2}^{i}\circ\tau y_{1}(\_\otimes e)\bigr)\in G_{1}.
\end{align*}

It remains to convert this to bimodule form. In the calculation we
will use three facts, easily checked by the reader: \\

\begin{itemize}
\item $x_{2}^{i}\circ\tau=\tau\circ x_{1}^{i}-h_{i-1}(x_{1},x_{2}),$
\item $x_{2}^{j}=y^{j}+y_{2}h_{i-1}(x_{2},y),$
\item $\sum_{j+k=i-1}x_{1}^{j}h_{k-1}(x_{2},y)=h_{i-2}(x_{1},x_{2},y).$\\
\end{itemize}
Then we have for the main calculation: 
\begin{align*}
 & Ef\circ\tau\circ x_{2}^{i}\circ\tau y_{1}(\_\otimes e)\\
 & =-Ef\circ\tau y_{1}\circ h_{i-1}(x_{1},x_{2})(\_\otimes e)\\
 & =-Ef\circ h_{i-1}(x_{1},x_{2})(\_\otimes e)-y_{1}Ef\circ\tau\circ h_{i-1}(x_{1},x_{2})(\_\otimes e)\\
 & =-Ef\circ\sum_{j+k=i-1}x_{1}^{j}\bigl(y^{k}+y_{2}h_{k-1}(x_{2},y)\bigr)(\_\otimes e)-y_{1}\circ Ef\circ\tau\circ h_{i-1}(x_{1},x_{2})(\_\otimes e)\\
 & =-Ef\circ h_{i-1}(x_{1},y)(\_\otimes e)-y_{1}Ef\ci\bigl(h_{i-2}(x_{1},x_{2},y)(\_\otimes e)+\tau\circ h_{i-1}(x_{1},x_{2})(\_\otimes e)\bigr).
\end{align*}
 Then observe that: 
\begin{align*}
-Ef\circ h_{i-1}(x_{1},y)(\_\otimes e) & =\_\otimes f(-h_{i-1}(x,y)e)\\
 & =(-\veps\circ h_{i-1}(x,y)F)(e\otimes f),
\end{align*}
 and that: 
\begin{align*}
 & -Ef\ci\Bigl(h_{i-2}(x_{1},x_{2},y)(\_\otimes e)+\tau\circ h_{i-1}(x_{1},x_{2})(\_\otimes e)\Bigr)\\
 & =\Bigl(-Ef\ci F\bigl(\tau\circ h_{i-1}(x_{1},x_{2})+h_{i-2}(x_{1},x_{2},y)\bigr)\circ\eta E\Bigr)(e)\\
 & =\Bigl(-FE\veps\circ F\bigl(\tau\circ h_{i-1}(x_{1},x_{2})+h_{i-2}(x_{1},x_{2},y)\bigr)F\circ\eta EF\Bigr)(e\otimes f).
\end{align*}
(We are using that $\_\otimes e$ considered in $\mathrm{Hom}_{A}(_{A}E,E^{2})[y]$
corresponds to $(\eta E)(e)$ in $FE^{2}[y]$; also note that $\veps(e\otimes f)=f(e)\in A[y]$
induces $Ef(\_\otimes e)=\_.f(e)$ considered in $\mathrm{Hom}_{A}(_{A}E,E)[y]$.)
The formulas in the last column of $[\tilde{\veps}\circ\tilde{x}^{i}\tilde{F}]_{22}$
follow.
\item We have $[\tilde{F}\tilde{x}^{i}\circ\tilde{\eta}]_{22}:[C]_{22}\to[\tilde{F}\tilde{E}]_{22}$
given by: 
\[
\begin{pmatrix}Fy^{i}\circ\eta & y^{i}y_{1}\\
-Fh_{i-1}(x,y)\circ\eta & y^{i}\\
0 & 0\\
Fx^{i}\circ\eta & 0\\
F^{2}\bigl(h_{i-1}(x_{1},x_{2})\ci\tau-h_{i-2}(x_{1},x_{2},y)\bigr)\circ\eta^{2} & F^{2}h_{i-1}(x_{2},y)\circ F\eta E
\end{pmatrix}
\]
 using the ordered decompositions: 
\begin{itemize}
\item $[C]_{22}\cong G_{1}\cong A[y]\oplus FE[y]$,
\item $[\tilde{F}\tilde{E}]_{22}\cong U\cong FE[y]^{\oplus4}\oplus F^{2}E^{2}[y]$.
\end{itemize}
Observe first that $[\tilde{\eta}]_{22}:G_{1}\to U$ is determined
by $(1,0)\mapsto\idop_{R}=(1,0,0,1,0)\in U$ (using bimodule forms).
Recall (Lemma \ref{lem:sigma-21-eta} used for $[\tilde{\sigma}]_{21}$)
that: 
\begin{multline*}
(L_{2})_{G}G_{2}\ni[\tilde{\eta}(1)]=\sum_{a\in Q}(f_{a},0,0)\otimes(e_{a},0,0)+\sum_{b\in Q}(0,f_{b},0)\otimes(0,e_{b},0)\\
\overset{\Gamma_{21}}{\longmapsto}(1,0,0,1,0)\in U.
\end{multline*}
 The map $\Gamma_{21}|_{(L_{2})_{G}G_{2}}$ of Eq.~(\ref{eq:FE-expansion})
is given by composition and hence right $G_{1}^{\op}$-equivariant,
so we can compute any $[\tilde{\eta}]_{22}\bigl((\theta,\varphi_{1})\bigr)$
as $[\tilde{\eta}(1)].(\theta,\varphi_{1})\in(L_{2})_{G}G_{2}$. The
action of $[\tilde{F}\tilde{x}^{i}]$ is applied to elements of $(L_{2})_{G}G_{2}$,
and after that we pass through $\Gamma_{21}$ again to obtain the
final image in $U$.

We treat the first column of $[\tilde{F}\tilde{x}^{i}\circ\tilde{\eta}]_{22}$
first, and consider the second column afterwards. For the first column
it is enough to consider the case $(\theta,\varphi_{1})=(1,0)$. Starting
with the first term, the data $(e_{a},0,0)$ corresponds to $e_{1}=0$,
$e_{2}=-y_{1}e_{a}$, and $\xi=y_{2}\tau(\_\otimes(-y_{1}e_{a}))$.
Application of the formula for $[\tilde{x}^{i}]_{22}$ gives $e_{1}=0$,
$e_{2}=-x^{i}y_{1}e_{a}$, and $\xi=x_{2}^{i}\ci y_{2}\tau(\_\otimes(-y_{1}e_{a}))$.
Then we convert this to bimodule form, using: 
\begin{align*}
 & x_{2}^{i}\circ y_{2}\tau(\_\otimes(-y_{1}e_{a}))\\
 & =y_{2}\ci x_{2}^{i}\tau(\_\otimes(-y_{1}e_{a}))\\
 & =y_{2}\ci\tau x_{1}^{i}(\_\otimes(-y_{1}e_{a}))+y_{1}y_{2}h_{i-1}(x_{1},x_{2})(\_\otimes e_{a})\\
 & =y_{2}\tau(\_\otimes(-y_{1}x^{i}e_{a}))+y_{1}y_{2}h_{i-1}(x_{1},x_{2})(\_\otimes e_{a}),
\end{align*}
 where in the third line we have used the first fact given under the
previous bullet. So in bimodule form we have: 
\[
[\tilde{x}^{i}]_{22}:(e_{a},0,0)\mapsto\bigl(x^{i}e_{a},0,h_{i-1}(x_{1},x_{2})(\_\otimes e_{a})\bigr).
\]
 Now applying $\Gamma_{21}$ we obtain: 
\[
\sum_{a\in Q}\bigl(x^{i}e_{a},0,h_{i-1}(x_{1},x_{2})(\_\otimes e_{a})\bigr)\circ(f_{a},0,0)=\bigl(0,0,0,x^{i},h_{i-1}(x_{1},x_{2})\circ\tau\bigr)\in U,
\]
 where the last component is computed using: 
\begin{align*}
 & \Bigl(y_{2}\tau(\_\otimes(-y_{1}x^{i}e_{a}))+y_{1}y_{2}h_{i-1}(x_{1},x_{2})(\_\otimes e_{a})\Bigr)\circ Ef_{a}\circ\tau\\
 & =-y_{2}\tau y_{1}x_{1}^{i}\tau+y_{1}y_{2}h_{i-1}(x_{1},x_{2})\tau,
\end{align*}
 together with the facts that $\Phi_{11}=\Phi_{12}=\Phi_{21}=0$ and
$\Phi_{22}=x^{i}$ so: 
\begin{align*}
\Lambda & =\tau y_{1}(0+0)-y_{2}\tau y_{1}\circ(0+E\Phi_{22}\circ\tau)+y_{1}y_{2}\Lambda^{\circ}\\
 & =-y_{2}\tau y_{1}\circ x_{1}^{i}\circ\tau+y_{1}y_{2}\Lambda^{\circ}.
\end{align*}

Continuing with the second term, the data $(0,e_{b},0)$ corresponds
to $e_{1}=e_{b}$, $e_{2}=e_{b}$, and $\xi=\tau y_{1}(\_\otimes e_{b})$.
Application of the formula for $[\tilde{x}^{i}]_{22}$ gives $e_{1}=y^{i}e_{b}$,
$e_{2}=x^{i}e_{b}$, and $\xi=x_{2}^{i}\ci\tau y_{1}(\_\otimes e_{b})$.
Then we convert this to bimodule form, using: 
\begin{align*}
 & x_{2}^{i}\circ\tau y_{1}(\_\otimes e_{b})\\
 & =\tau y_{1}(\_\otimes x^{i}e_{b})-y_{1}h_{i-1}(x_{1},x_{2})(\_\otimes e_{b})\\
 & =\tau y_{1}(\_\otimes x^{i}e_{b})-y_{1}h_{i-1}(x_{1},y)(\_\otimes e_{b})\\
 & \qquad-y_{1}y_{2}h_{i-2}(x_{1},x_{2},y)(\_\otimes e_{b})\\
 & =\_\otimes x^{i}e_{b}+y_{2}\tau(\_\otimes x^{i}e_{b})-\_\otimes(x^{i}-y^{i})e_{b}\\
 & \qquad-y_{1}y_{2}h_{i-2}(x_{1},x_{2},y)(\_\otimes e_{b})\\
 & =\_\otimes y^{i}e_{b}+y_{2}\tau(\_\otimes x^{i}e_{b})-y_{1}y_{2}h_{i-2}(x_{1},x_{2},y)(\_\otimes e_{b}),
\end{align*}
 where we have made use of the fact, easily checked by the reader,
that: 
\begin{itemize}
\item $y_{2}h_{i-2}(x_{1},x_{2},y)=h_{i-1}(x_{1},x_{2})-h_{i-1}(x_{1},y).$
\end{itemize}
So in bimodule form we have: 
\[
[\tilde{x}^{i}]_{22}:(0,e_{b},0)\mapsto\Bigl(-h_{i-1}(x_{1},y)e_{b},y^{i}e_{b},-h_{i-2}(x_{1},x_{2},y)(\_\otimes e_{b})\Bigr).
\]
 Now applying $\Gamma_{21}$ we obtain: 
\begin{align*}
 & \sum_{b\in Q}\Bigl(-h_{i-1}(x_{1},y)e_{b},y^{i}e_{b},-h_{i-2}(x_{1},x_{2},y)(\_\otimes e_{b})\Bigr)\circ(0,f_{b},0)\\
 & =\bigl(y^{i},-h_{i-1}(x_{1},y),0,0,-h_{i-2}(x_{1},x_{2},y)\bigr)\in U,
\end{align*}
where the last component is computed using: 
\[
x_{2}^{i}\circ\tau y_{1}(\_\otimes e_{b})\circ Ef_{b}=x_{2}^{i}\tau y_{1}=\tau x_{1}^{i}y_{1}-y_{1}h_{i-1}(x_{1},x_{2})
\]
 together with the facts that $\Phi_{11}=y^{i}$, $\Phi_{21}=-h_{i-1}(x_{1},y)$,
$\Phi_{12}=\Phi_{22}=0$, so: 
\begin{align*}
\Lambda & =\tau y_{1}\circ\bigl(y^{i}+0\ci\tau\bigr)-y_{2}\tau y_{1}\circ\bigl(-h_{i-1}(x_{1},y)+0\ci\tau\bigr)+y_{1}y_{2}\Lambda^{\circ}\\
 & =\tau y_{1}y^{i}+y_{2}\tau y_{1}h_{i-1}(x_{1},y)+y_{1}y_{2}\Lambda^{\circ}\\
 & =\tau y_{1}y^{i}+y_{2}\tau(x_{1}^{i}-y^{i})+y_{1}y_{2}\Lambda^{\circ}\\
 & =y^{i}+y_{2}\tau x_{1}^{i}+y_{1}y_{2}\Lambda^{\circ}\\
 & =\tau x_{1}^{i}y_{1}-y_{1}h_{i-1}(x_{1},y)+y_{1}y_{2}\Lambda^{\circ},
\end{align*}
 so, using again the fact above: 
\begin{align*}
y_{2}\Lambda^{\circ} & =-h_{i-1}(x_{1},x_{2})+h_{i-1}(x_{1},y),\\
\Lambda^{\circ} & =-h_{i-2}(x_{1},x_{2},y).
\end{align*}
 Finally taking the sum of the two terms, we conclude that $[\tilde{F}\tilde{x}^{i}\circ\tilde{\eta}]_{22}:A[y]\to U$
is determined by: 
\[
1\mapsto\bigl(y^{i},-h_{i-1}(x_{1},y),0,x^{i},h_{i-1}(x_{1},x_{2})\circ\tau-h_{i-2}(x_{1},x_{2},y)\bigr).
\]
 By describing these coefficients in $FE[y]$ and $F^{2}E^{2}[y]$
instead of in $\End(E[y])$ and $\End(E^{2}[y])$, we obtain the formulas
in the first column of the matrix of $[\tilde{F}\tilde{x}^{i}\circ\tilde{\eta}]_{22}$.

Now we consider the second column of $[\tilde{F}\tilde{x}^{i}\circ\tilde{\eta}]_{22}$,
a map $FE[y]\to U$. It is found using the same method but with $(\theta,\varphi_{1})=(0,\varphi_{1})$
for a generator $\varphi_{1}\in FE[y]$. We have in bimodule form:
\begin{align*}
(e_{a},0,0).(0,\varphi_{1}) & =\bigl(0,0,E\varphi_{1}\circ\tau(\_\otimes-y_{1}e_{a})\bigr)\\
(0,e_{b},0).(0,\varphi_{1}) & =\bigl(\varphi_{1}(e_{b}),y_{1}\varphi_{1}(e_{b}),E\varphi_{1}\circ\tau(\_\otimes e_{b})\bigr),
\end{align*}
 where we have used the calculations: 
\[
E(y_{1}\varphi_{1})\circ y_{2}\tau(\_\otimes-y_{1}e_{a})=y_{1}y_{2}E\varphi_{1}\circ\tau(\_\otimes-y_{1}e_{a})
\]
 and 
\[
E(y_{1}\varphi_{1})\circ(\_\otimes e_{b}+y_{2}\tau(\_\otimes e_{b}))=\_\otimes y_{1}\varphi_{1}(e_{b})+y_{1}y_{2}E\varphi_{1}\circ\tau(\_\otimes e_{b}).
\]
 Starting with the first term, the data $\bigl(0,0,E\varphi_{1}\circ\tau(\_\otimes-y_{1}e_{a})\bigr)$
corresponds to $e_{1}=e_{2}=0$ and $\xi=y_{1}y_{2}E\varphi_{1}\circ\tau(\_\otimes-y_{1}e_{a})$.
Application of the formula for $[\tilde{x}^{i}]_{22}$ gives $e_{1}=e_{2}=0$
and $\xi=x_{2}^{i}\ci y_{1}y_{2}E\varphi_{1}\circ\tau(\_\otimes-y_{1}e_{a})$.
Converting this data to bimodule form is trivial. So we have: 
\[
[\tilde{x}^{i}]_{22}:\bigl(0,0,E\varphi_{1}\circ\tau(\_\otimes-y_{1}e_{a})\bigr)\mapsto\bigl(0,0,x_{2}^{i}\circ E\varphi_{1}\circ\tau(\_\otimes-y_{1}e_{a})\bigr).
\]
 Now applying $\Gamma_{21}$ we obtain: 
\[
\sum_{a\in Q}\bigl(0,0,x_{2}^{i}\circ E\varphi_{1}\circ\tau(\_\otimes-y_{1}e_{a})\bigr)\circ(f_{a},0,0)=\bigl(0,0,0,0,-E\varphi_{1}\circ x_{2}^{i}\tau\bigr)\in U,
\]
 where the last component is computed using: 
\begin{align*}
 & y_{1}y_{2}x_{2}^{i}\ci E\varphi_{1}\ci\tau(\_\otimes-y_{1}e_{a})\ci Ef_{a}\ci\tau\\
 & =-x_{2}^{i}\ci y_{1}y_{2}E\varphi_{1}\ci\tau y_{1}\ci\tau\\
 & =-y_{1}y_{2}E\varphi_{1}\ci x_{2}^{i}\tau.
\end{align*}
 Continuing with the second term, the data $\bigl(\varphi_{1}(e_{b}),y_{1}\varphi_{1}(e_{b}),E\varphi_{1}\circ\tau(\_\otimes e_{b})\bigr)$
corresponds to $e_{1}=y_{1}\varphi_{1}(e_{b})$, $e_{2}=0$, and $\xi=\_\otimes y_{1}\varphi_{1}(e_{b})+y_{1}y_{2}E\varphi_{1}\circ\tau(\_\otimes e_{b})$.
Application of the formula for $[\tilde{x}^{i}]_{22}$ gives $e_{1}=y_{1}y^{i}\varphi_{1}(e_{b})$,
$e_{2}=0$, and $\xi=x_{2}^{i}\ci\bigl(\_\otimes y_{1}\varphi_{1}(e_{b})+y_{1}y_{2}E\varphi_{1}\circ\tau(\_\otimes e_{b})\bigr)$.
Then we convert this to bimodule form, using: 
\begin{align*}
 & x_{2}^{i}\ci\bigl(\_\otimes y_{1}\varphi_{1}(e_{b})+y_{1}y_{2}E\varphi_{1}\circ\tau(\_\otimes e_{b})\bigr)\\
 & =\_\otimes y^{i}y_{1}\varphi_{1}(e_{b})+y_{2}h_{i-1}(x_{2},y)(\_\otimes y_{1}\varphi_{1}(e_{b}))+y_{1}y_{2}E\varphi_{1}\ci x_{2}^{i}\tau(\_\otimes e_{b})\\
 & =\_\otimes y^{i}y_{1}\varphi_{1}(e_{b})+y_{1}y_{2}\Bigl(h_{i-1}(x_{2},y)\ci E\varphi_{1}(\_\otimes e_{b})+E\varphi_{1}\ci x_{2}^{i}\tau(\_\otimes e_{b})\Bigr)\\
 & =\_\otimes y^{i}y_{1}\varphi_{1}(e_{b})+y_{1}y_{2}E\varphi_{1}\ci\bigl(x_{2}^{i}\tau+h_{i-1}(x_{2},y)\bigr)(\_\otimes e_{b})\\
 & =\_\otimes y_{1}y^{i}\varphi_{1}(e_{b})+y_{1}y_{2}\Bigl(-E\varphi_{1}\ci y_{1}h_{i-2}(x_{1},x_{2},y)(\_\otimes e_{b})+E\varphi_{1}\ci\tau\ci x_{1}^{i}(\_\otimes e_{b})\Bigr).
\end{align*}
 So in bimodule form we have: 
\begin{align*}
 & [\tilde{x}^{i}]_{22}:\bigl(\varphi_{1}(e_{b}),y_{1}\varphi_{1}(e_{b}),E\varphi_{1}\circ\tau(\_\otimes e_{b})\bigr)\mapsto\\
 & \Bigl(y^{i}\varphi_{1}(e_{b}),y_{1}y^{i}\varphi_{1}(e_{b}),E\varphi_{1}\ci\bigl(x_{2}^{i}\tau+h_{i-1}(x_{2},y)\bigr)(\_\otimes e_{b})\Bigr).
\end{align*}
 Now applying $\Gamma_{21}$ we obtain: 
\begin{align*}
 & \sum_{b\in Q}\Bigl(y^{i}\varphi_{1}(e_{b}),y_{1}y^{i}\varphi_{1}(e_{b}),E\varphi_{1}\ci\bigl(x_{2}^{i}\tau+h_{i-1}(x_{2},y)\bigr)(\_\otimes e_{b})\Bigr)\circ(0,f_{b},0)\\
 & =\Bigl(y^{i}y_{1}\varphi_{1},y^{i}\varphi_{1},0,0,E\varphi_{1}\ci\bigl(x_{2}^{i}\tau+h_{i-1}(x_{2},y)\bigr)\Bigr)\in U,
\end{align*}
where the last component is computed using: 
\begin{align*}
 & \Bigl(\_\otimes y^{i}y_{1}\varphi_{1}(e_{b})+y_{1}y_{2}E\varphi_{1}\ci\bigl(x_{2}^{i}\tau+h_{i-1}(x_{2},y)\bigr)(\_\otimes e_{b})\Bigr)\ci Ef_{b}\\
 & =y^{i}y_{1}E\varphi_{1}+y_{1}y_{2}E\varphi_{1}\ci\bigl(x_{2}^{i}\tau+h_{i-1}(x_{2},y)\bigr)(\_\otimes e_{b}),
\end{align*}
 together with the facts that $\Phi_{11}=y^{i}y_{1}\varphi_{1}$,
$\Phi_{21}=y^{i}\varphi_{1}$, $\Phi_{12}=\Phi_{22}=0$, so: 
\begin{align*}
\Lambda & =\tau y_{1}(y^{i}y_{1}E\varphi_{1}+0\ci\tau)-y_{2}\tau y_{1}(y^{i}E\varphi_{1}+0\ci\tau)+y_{1}y_{2}\Lambda^{\circ}\\
 & =y^{i}y_{1}E\varphi_{1}+y_{1}y_{2}\Lambda^{\ci}.
\end{align*}
 Taking the sum of the two terms, we conclude that $[\tilde{F}\tilde{x}^{i}\circ\tilde{\eta}]_{22}:FE[y]\to U$
is given by: 
\[
\varphi_{1}\mapsto\bigl(y^{i}y_{1}\varphi_{1},y^{i}\varphi_{1},0,0,E\varphi_{1}\ci h_{i-1}(x_{2},y)\bigr).
\]
 The last component, an element of $\End_{A}(_{A}E^{2})[y]$, is the
same as $\bigl(F^{2}h_{i-1}(x_{2},y)\circ F\eta E\bigr)(\varphi_{1})$.
This gives the formulas in the second column of the matrix of $[\tilde{F}\tilde{x}^{i}\circ\tilde{\eta}]_{22}$.
\end{itemize}

\subsection{Maps $\tilde{\rho}_{\lambda}$: isomorphisms}

Now we have formulas by components for the maps $\tilde{\sigma}$,
$\tilde{\veps}\circ\tilde{x}^{i}\tilde{F}$, and $\tilde{F}\tilde{x}^{i}\circ\tilde{\eta}$
that are used to define the maps $\tilde{\rho}_{\lambda}$. It remains
to make use of the isomorphisms $\rho_{\lambda}$ determined by $\sigma$,
$\veps\circ x^{i}F$, and $Fx^{i}\circ\eta$, together with these
formulas, to show that $\tilde{\rho}_{\lambda}$ are isomorphisms.
Note that $\tilde{\rho}_{\lambda}$ are already known to give morphisms
of $(C,C)$-bimodules, so it suffices to show that $\tilde{\rho}_{\lambda}$
are isomorphisms of sets. We will work again by components and show
that $[\tilde{\rho}_{\lambda}]_{ij}$ is an isomorphism of $(A[y],A[y])$-bimodules
for $i,j\in\{1,2\}$.

We remind the reader of our notational convention that $E_{\lambda}=Ee_{\lambda}$
for the idempotents $e_{\lambda}\in A_{\lambda}$ of a weight decomposition.
Recall that the bimodule $E$ satisfies $e_{j}Ee_{i}={\delta_{i+2,j}\cdot e_{i+2}Ee_{i}}$,
and similarly for $F$ but with $i-2$ instead of $i+2$. Finally,
recall Prop.~4.26 of \cite{mcmillanTensor2product2representations2022}
that gives the weight idempotents for the algebra $C$.
\begin{itemize}
\item We have for $[\tilde{\rho}_{\lambda}]_{11}$, $\lambda\geq0$: 
\[
[\tilde{\rho}_{\lambda}]_{11}:EF_{\lambda+1}[y]\to A_{\lambda+1}[y]\oplus FE_{\lambda+1}[y]\oplus A_{\lambda+1}[y]^{\oplus\lambda}
\]
 given by: 
\[
[\tilde{\rho}_{\lambda}]_{11}=\veps\oplus\sigma\oplus\bigoplus_{i=0}^{\lambda-1}\veps\ci x^{i}y_{1}F.
\]
\item We have for $[\tilde{\rho}_{\lambda}]_{11}$, $\lambda\leq0$: 
\[
[\tilde{\rho}_{\lambda}]_{11}:EF_{\lambda+1}[y]\oplus A_{\lambda+1}[y]^{\oplus-\lambda}\to A_{\lambda+1}[y]\oplus FE_{\lambda+1}[y]
\]
 given by: 
\[
[\tilde{\rho}_{\lambda}]_{11}=\left(\left(\begin{smallmatrix}\veps\\
\sigma
\end{smallmatrix}\right),\sum_{i=0}^{-\lambda-1}\left(\begin{smallmatrix}y^{i}\\
Fh_{i-1}(x,y)\ci\eta
\end{smallmatrix}\right)\right).
\]

\begin{prop}
The morphism of $(A[y],A[y])$-bimodules $[\tilde{\rho}_{\lambda}]_{11}$
is an isomorphism for all $\lambda$.
\end{prop}

\begin{proof}
When $\lambda\geq0$ and therefore $\lambda+1\geq0$, the map: 
\[
\sigma\oplus\bigoplus_{i=0}^{\lambda}\veps\ci x^{i}F:EF_{\lambda+1}[y]\iso FE_{\lambda+1}[y]\oplus A_{\lambda+1}[y]^{\oplus\lambda+1}
\]
 is just $\rho_{\lambda+1}\otimes_{k}k[y]$. It is an isomorphism
because $\rho_{\lambda+1}$ is an isomorphism.
\begin{claim}
\label{claim:veps-x^iy_1} When $\lambda\geq0$, the map 
\[
\sigma\oplus\veps\oplus\bigoplus_{i=0}^{\lambda-1}\veps\ci x^{i}y_{1}F:EF_{\lambda+1}[y]\to FE_{\lambda+1}[y]\oplus A_{\lambda+1}[y]^{\oplus\lambda+1}
\]
 is also an isomorphism.
\end{claim}

\begin{proof}
Let $M_{-y}\in\End_{A_{\lambda+1}[y]}\left(A_{\lambda+1}[y]^{\oplus\lambda+1}\right)$
be the endomorphism with matrix coefficients $[M_{-y}]\in\mathrm{Mat}_{(\lambda+1)\times(\lambda+1)}\left(A_{\lambda+1}[y]^{\op}\right)$
given by $1$ on the diagonal and $-y$ on the subdiagonal, and $0$
elsewhere. This matrix is invertible, and $M_{-y}$ is an isomorphism.
Observe that: 
\[
\veps\ci(-x^{i-1}yF)=-y\cdot\veps\ci x^{i-1}F.
\]
 Using this we write the map in question as a composition of isomorphisms:
\[
\sigma\oplus\veps\oplus\bigoplus_{i=0}^{\lambda-1}\veps\ci x^{i}y_{1}F=\begin{pmatrix}1 & 0\\
0 & M_{-y}
\end{pmatrix}\ci\left(\sigma\oplus\bigoplus_{i=0}^{\lambda}\veps\ci x^{i}F\right).
\]
 By reordering the first two summands in the codomain, we obtain the
map $[\tilde{\rho}_{\lambda}]_{11}$.
\end{proof}
When $\lambda=0$, the two formulas for $[\tilde{\rho}_{\lambda}]_{11}$
agree. Now assume $\lambda<0$, so $\lambda+1\leq0$ and the map:
\begin{equation}
\left(\sigma,\sum_{i=0}^{-(\lambda+1)-1}Fx^{i}\ci\eta\right):EF_{\lambda+1}\oplus A_{\lambda+1}[y]^{\oplus-(\lambda+1)}\iso FE_{\lambda+1}[y]\label{eq:sigma-eta-iso}
\end{equation}
 is $\rho_{\lambda+1}\otimes_{k}k[y]$, an isomorphism.
\begin{claim}
\label{claim:eta-h(x,y)} When $\lambda<0$, the map: 
\[
\left(\sigma,\sum_{i=1}^{-\lambda-1}Fh_{i-1}(x,y)\ci\eta\right):EF_{\lambda+1}[y]\oplus A_{\lambda+1}[y]^{\oplus-(\lambda+1)}\to FE_{\lambda+1}[y]
\]
 is also an isomorphism.
\end{claim}

\begin{proof}
This time we define an isomorphism $M_{h}\in\End_{A_{\lambda+1}[y]}\left(A_{\lambda+1}[y]^{\oplus-(\lambda+1)}\right)$
with components $[M_{h}]_{ii}=1$, $[M_{h}]_{ij}=y^{j-i}$ for $j>i$,
and $[M_{h}]_{ij}=0$ for $j<i$. This is an upper-triangular invertible
matrix: 
\[
[M_{h}]=\left(\begin{smallmatrix}1 & y & y^{2} & \dots & y^{-(\lambda+1)-1}\\
0 & 1 & y & \dots & y^{-(\lambda+1)-2}\\
0 & 0 & 1 & \dots & y^{-(\lambda+1)-3}\\
\dots & \dots & \dots & \dots & \dots\\
0 & 0 & 0 & \dots & 1
\end{smallmatrix}\right).
\]
 Now observe that $Fx^{i}y^{j}\ci\eta=(Fx^{i}\ci\eta)\cdot y^{j}$.
We use this and write: 
\begin{align*}
\sum_{i=1}^{-\lambda-1}Fh_{i-1}(x,y)\ci\eta & =\sum_{i=0}^{-\lambda-2}\sum_{j+k=i}Fx^{j}\ci\eta\cdot y^{k}\\
 & =\left(\sum_{i=0}^{-(\lambda+1)-1}Fx^{i}\ci\eta\right)\ci M_{h},
\end{align*}
 and it follows from this and the isomorphism above the claim that
the map of the claim is an isomorphism.
\end{proof}
By writing out terms, we have: 
\[
\left(\left(\begin{smallmatrix}\veps\\
\sigma
\end{smallmatrix}\right),\sum_{i=0}^{-\lambda-1}\left(\begin{smallmatrix}y^{i}\\
Fh_{i-1}(x,y)\ci\eta
\end{smallmatrix}\right)\right)=\begin{pmatrix}\veps & 1 & y & \dots & y^{-\lambda-1}\\
\sigma & 0 & \eta & \dots & Fh_{-\lambda-2}(x,y)\ci\eta
\end{pmatrix}.
\]
 Interchanging the first two summands of the domain, we obtain the
form: 
\[
\begin{pmatrix}1 & \bigl(\veps,y,y^{2},\dots,y^{-\lambda-1}\bigr)\\
0 & \left(\sigma,\sum\limits _{i=1}^{-\lambda-1}Fh_{i-1}(x,y)\ci\eta\right)
\end{pmatrix},
\]
 which (by the claim) is manifestly an isomorphism.
\end{proof}
\item We have for $[\tilde{\rho}_{\lambda}]_{21}$, $\lambda\geq0$: 
\[
[\tilde{\rho}_{\lambda}]_{21}:F_{\lambda+1}[y]\oplus FEF_{\lambda+1}[y]\to F_{\lambda+1}[y]\oplus F_{\lambda+1}[y]\oplus F^{2}E_{\lambda+1}[y]\oplus F_{\lambda+1}[y]^{\oplus\lambda}
\]
 given by: 
\[
[\tilde{\rho}_{\lambda}]_{21}=\begin{pmatrix}1 & 0\\
0 & F\veps\\
0 & F\sigma\\
\bigoplus\limits _{i=0}^{\lambda-1}x^{i} & \bigoplus\limits _{i=0}^{\lambda-1}F(\veps\circ x^{i}y_{1}F)
\end{pmatrix}.
\]
\item We have for $[\tilde{\rho}_{\lambda}]_{21}$, $\lambda\leq0$: 
\[
[\tilde{\rho}_{\lambda}]_{21}:F_{\lambda+1}[y]\oplus FEF_{\lambda+1}[y]\oplus F_{\lambda+1}[y]^{\oplus-\lambda}\to F_{\lambda+1}[y]\oplus F_{\lambda+1}[y]\oplus F^{2}E_{\lambda+1}[y]
\]
 given by: 
\[
[\tilde{\rho}_{\lambda}]_{21}=\begin{pmatrix}1 & 0 & 0\\
0 & F\veps & \sum\limits _{i=0}^{-\lambda-1}y^{i}\\
0 & F\sigma & \sum\limits _{i=0}^{-\lambda-1}F(Fh_{i-1}(x,y)\ci\eta)
\end{pmatrix}.
\]

\begin{prop}
The morphism of $(A[y],A[y])$-bimodules $[\tilde{\rho}_{\lambda}]_{21}$
is an isomorphism for all $\lambda$.
\end{prop}

\begin{proof}
When $\lambda\geq0$, we have that 
\begin{multline*}
F\veps\oplus F\sigma\oplus\bigoplus_{i=0}^{\lambda-1}F(\veps\circ x^{i}y_{1}F):FEF_{\lambda+1}[y]\\
\to F_{\lambda+1}[y]\oplus F^{2}E_{\lambda+1}[y]\oplus F_{\lambda+1}[y]^{\oplus\lambda}
\end{multline*}
 is an isomorphism, using Claim \ref{claim:veps-x^iy_1} and the fact
that (horizontal) composition of the identity functor on $F$ with
an isomorphism gives an isomorphism. Then $[\tilde{\rho}_{\lambda}]_{21}$
may be compressed to a lower-triangular $2\times2$ matrix with an
isomorphism in position $(2,2)$, so it is an isomorphism.

When $\lambda=0$, the two formulas for $[\tilde{\rho}_{\lambda}]_{21}$
agree. Assume now that $\lambda<0$, so the map 
\[
\left(F\sigma,\sum_{i=1}^{-\lambda-1}F\bigl(Fh_{i-1}(x,y)\ci\eta\bigr)\right):FEF_{\lambda+1}[y]\oplus F_{\lambda+1}[y]^{\oplus-(\lambda+1)}\to F^{2}E_{\lambda+1}[y]
\]
 is an isomorphism using Claim \ref{claim:eta-h(x,y)}. Now expand
the notation of the map $[\tilde{\rho}_{\lambda}]_{21}$ in the third
column: 
\[
\begin{pmatrix}1 & 0 & 0 & 0\\
0 & F\veps & 1 & \sum\limits _{i=1}^{-\lambda-1}y^{i}\\
0 & F\sigma & 0 & \sum\limits _{i=1}^{-\lambda-1}F\bigl(Fh_{i-1}(x,y)\ci\eta\bigr)
\end{pmatrix}.
\]
 After switching the second and third summands of the domain, we obtain
an upper-triangular matrix with isomorphisms on the diagonal, so $[\tilde{\rho}_{\lambda}]_{21}$
is an isomorphism.
\end{proof}
\item We have for $[\tilde{\rho}_{\lambda}]_{12}$, $\lambda\geq0$: 
\[
[\tilde{\rho}_{\lambda}]_{12}:E_{\lambda-1}[y]\oplus EFE_{\lambda-1}[y]\to E_{\lambda-1}[y]\oplus E_{\lambda-1}[y]\oplus FE_{\lambda-1}^{2}[y]\oplus E_{\lambda-1}[y]^{\oplus\lambda}
\]
 given by: 
\[
[\tilde{\rho}_{\lambda}]_{12}=\begin{pmatrix}0 & \veps E\\
1 & y_{1}\ci\veps E\\
0 & \sigma E\\
\bigoplus\limits _{i=0}^{\lambda-1}x^{i} & \bigoplus\limits _{i=0}^{\lambda-1}(\veps\circ x^{i}y_{1}F)E
\end{pmatrix}.
\]
\item We have for $[\tilde{\rho}_{\lambda}]_{12}$, $\lambda\leq0$: 
\[
[\tilde{\rho}_{\lambda}]_{12}:E_{\lambda-1}[y]\oplus EFE_{\lambda-1}[y]\oplus E_{\lambda-1}[y]^{\oplus-\lambda}\to E_{\lambda-1}[y]\oplus E_{\lambda-1}[y]\oplus FE_{\lambda-1}^{2}[y]
\]
 given by: 
\[
[\tilde{\rho}_{\lambda}]_{12}=\begin{pmatrix}0 & \veps E & \sum\limits _{i=0}^{-\lambda-1}y^{i}\\
1 & y_{1}\ci\veps E & \sum\limits _{i=0}^{-\lambda-1}y^{i}y_{1}\\
0 & \sigma E & \sum\limits _{i=0}^{-\lambda-1}(Fh_{i-1}(x,y)\ci\eta)E
\end{pmatrix}.
\]

\begin{prop}
The morphism of $(A[y],A[y])$-bimodules $[\tilde{\rho}_{\lambda}]_{12}$
is an isomorphism for all $\lambda$.
\end{prop}

\begin{proof}
When $\lambda\geq0$, we have that 
\[
\veps E\oplus\sigma E\oplus\bigoplus_{i=0}^{\lambda-1}(\veps\circ x^{i}y_{1}F)E:EFE_{\lambda-1}[y]\to E_{\lambda-1}[y]\oplus FE_{\lambda-1}^{2}[y]\oplus E_{\lambda-1}[y]^{\oplus\lambda}
\]
 is an isomorphism, using Claim \ref{claim:veps-x^iy_1} with $E$
applied on the right. Note that $E$ applied on the right here is
equivalent to $_{\lambda+1}E_{\lambda-1}$ applied on the right, and
this raises the weight by $2$, so we still invoke the isomorphism
$\rho_{\lambda+1}$ for weight $\lambda+1$.

We perform some row operations on the matrix of $[\tilde{\rho}_{\lambda}]_{12}$.
Subtract $y_{1}$ times the first row from the second to eliminate
the coefficient $y_{1}\ci\veps E$. Then exchange the first and second
rows, then exchange the second and third rows, then collapse the second
and third into the notation of the fourth. Obtain: 
\[
\begin{pmatrix}1 & 0\\
0\oplus0\oplus\bigoplus\limits _{i=0}^{\lambda-1}x^{i} & \sigma E\oplus\veps E\oplus\bigoplus\limits _{i=0}^{\lambda-1}(\veps\circ x^{i}y_{1}F)E
\end{pmatrix},
\]
which is upper-triangular with isomorphisms on the diagonal, so the
original matrix for $[\tilde{\rho}_{\lambda}]_{12}$ is an isomorphism.

When $\lambda=0$, the two formulas for $[\tilde{\rho}_{\lambda}]_{12}$
agree. Assume now that $\lambda<0$, so the map 
\[
\left(\sigma E,\sum_{i=1}^{-\lambda-1}\bigl(Fh_{i-1}(x,y)\ci\eta\bigr)E\right):EFE_{\lambda-1}[y]\oplus E_{\lambda-1}[y]^{\oplus-(\lambda+1)}\to FE_{\lambda-1}^{2}[y]
\]
 is an isomorphism using Claim \ref{claim:eta-h(x,y)}. Now expand
the notation of the map $[\tilde{\rho}_{\lambda}]_{12}$ in the third
column: 
\[
[\tilde{\rho}_{\lambda}]_{12}=\begin{pmatrix}0 & \veps E & 1 & \sum\limits _{i=1}^{-\lambda-1}y^{i}\\
1 & y_{1}\ci\veps E & y_{1} & \sum\limits _{i=1}^{-\lambda-1}y^{i}y_{1}\\
0 & \sigma E & 0 & \sum\limits _{i=1}^{-\lambda-1}\bigl(Fh_{i-1}(x,y)\ci\eta\bigr)E
\end{pmatrix}.
\]
 Exchange the first and second rows, then the second and third columns,
then collapse the third and fourth columns into the notation of the
third, and obtain: 
\[
\begin{pmatrix}1 & y_{1} & \left(y_{1}\ci\veps E,\sum\limits _{i=1}^{-\lambda-1}y^{i}y_{1}\right)\\
0 & 1 & \left(\veps E,\sum\limits _{i=1}^{-\lambda-1}y^{i}\right)\\
0 & 0 & \left(\sigma E,\sum\limits _{i=1}^{-\lambda-1}\bigl(Fh_{i-1}(x,y)\ci\eta\bigr)E\right)
\end{pmatrix}.
\]
 Since this is upper-triangular with isomorphisms on the diagonal,
the original matrix $[\tilde{\rho}_{\lambda}]_{12}$ is an isomorphism.
\end{proof}
\item We have for $[\tilde{\rho}_{\lambda}]_{22}$, $\lambda\geq0$: 
\begin{multline}
[\tilde{\rho}_{\lambda}]_{22}:A_{\lambda-1}[y]\oplus FE_{\lambda-1}[y]^{\oplus2}\oplus FEFE_{\lambda-1}[y]\oplus EF_{\lambda-1}[y]\\
\to FE_{\lambda-1}[y]^{\oplus4}\oplus F^{2}E_{\lambda-1}^{2}[y]\oplus A_{\lambda-1}[y]^{\oplus\lambda}\oplus FE_{\lambda-1}[y]^{\oplus\lambda}\label{eq:rho-22-pos}
\end{multline}
 given by: $[\tilde{\rho}_{\lambda}]_{22}=${\footnotesize{}
\[
\begin{pmatrix}0 & 0 & 1 & 0 & 0\\
0 & 0 & 0 & F\veps E & 0\\
\eta & y_{1} & 0 & 0 & \sigma\\
0 & 1 & 0 & 0 & 0\\
0 & 0 & 0 & F\sigma E & 0\\
\bigoplus\limits _{i=0}^{\lambda-1}y^{i} & 0 & 0 & 0 & \bigoplus\limits _{i=0}^{\lambda-1}-\veps\circ h_{i-1}(x,y)F\\
\bigoplus\limits _{i=0}^{\lambda-1}h_{i-1}(x,y)\circ\eta & \bigoplus\limits _{i=0}^{\lambda-1}x^{i}E & \bigoplus\limits _{i=0}^{\lambda-1}Fx^{i} & \bigoplus\limits _{i=0}^{\lambda-1}F(\veps\circ x^{i}y_{1}F)E & \Theta
\end{pmatrix},
\]
}{\scriptsize{} }where 
\[
\Theta=\bigoplus\limits _{i=0}^{\lambda-1}-FE\veps\circ F\bigl(\tau\circ h_{i-1}(x_{1},x_{2})-h_{i-2}(x_{1},x_{2},y)\bigr)F\circ\eta EF.
\]
 
\item We have for $[\tilde{\rho}_{\lambda}]_{22}$, $\lambda\leq0$: 
\begin{multline}
[\tilde{\rho}_{\lambda}]_{22}:A_{\lambda-1}[y]\oplus FE_{\lambda-1}[y]^{\oplus2}\oplus FEFE_{\lambda-1}[y]\oplus EF_{\lambda-1}[y]\\
\oplus A_{\lambda-1}[y]^{\oplus-\lambda}\oplus FE_{\lambda-1}[y]^{\oplus-\lambda}\to FE_{\lambda-1}[y]^{\oplus4}\oplus F^{2}E_{\lambda-1}^{2}[y]\label{eq:rho-22-neg}
\end{multline}
 given by: $[\tilde{\rho}_{\lambda}]_{22}=${\small{}
\[
\begin{pmatrix}0 & 0 & 1 & 0 & 0 & \sum\limits _{i=0}^{-\lambda-1}Fy^{i}\circ\eta & \sum\limits _{i=0}^{-\lambda-1}y^{i}y_{1}\\
0 & 0 & 0 & F\veps E & 0 & \sum\limits _{i=0}^{-\lambda-1}-Fh_{i-1}(x,y)\circ\eta & \sum\limits _{i=0}^{-\lambda-1}y^{i}\\
\eta & y_{1} & 0 & 0 & \sigma & 0 & 0\\
0 & 1 & 0 & 0 & 0 & \sum\limits _{i=0}^{-\lambda-1}Fx^{i}\circ\eta & 0\\
0 & 0 & 0 & F\sigma E & 0 & \Theta' & \sum\limits _{i=0}^{-\lambda-1}F^{2}\bigl(h_{i-1}(x_{2},y)\bigr)\circ F\eta E
\end{pmatrix},
\]
} where 
\[
\Theta'=\sum\limits _{i=0}^{-\lambda-1}F^{2}\bigl(h_{i-1}(x_{1},x_{2})\ci\tau-h_{i-2}(x_{1},x_{2},y)\bigr)\circ\eta^{2}.
\]
 
\begin{prop}
The morphism of $(A[y],A[y])$-bimodules $[\tilde{\rho}_{\lambda}]_{22}$
is an isomorphism for all $\lambda$.
\end{prop}

\begin{proof}
When $\lambda>0$ and therefore $\lambda-1\geq0$, the map 
\[
\sigma\oplus\bigoplus_{i=0}^{\lambda-2}-\veps\circ x^{i}F:EF_{\lambda-1}[y]\to FE_{\lambda-1}[y]\oplus A_{\lambda-1}[y]^{\oplus\lambda-1}
\]
 is an isomorphism. (The minus sign does not interfere.)
\begin{claim}
When $\lambda>0$, the map 
\[
\sigma\oplus\bigoplus_{i=1}^{\lambda-1}-\veps\circ h_{i-1}(x,y)F:EF_{\lambda-1}[y]\to FE_{\lambda-1}[y]\oplus A_{\lambda-1}[y]^{\oplus\lambda-1}
\]
 is an isomorphism.
\end{claim}

\begin{proof}
Define an isomorphism $M_{h}'\in\End_{A_{\lambda-1}[y]}\left(A_{\lambda-1}[y]^{\oplus\lambda-1}\right)$
with components $[M_{h}']_{ii}=1$, $[M_{h}']_{ij}=y^{i-j}$ for $i>j$,
and $[M_{h}']_{ij}=0$ for $i<j$. This is a lower-triangular invertible
matrix: 
\[
[M_{h}']=\left(\begin{smallmatrix}1 & 0 & 0 & \dots & 0\\
y & 1 & 0 & \dots & 0\\
y^{2} & y & 1 & \dots & 0\\
\dots & \dots & \dots & \dots & \dots\\
y^{\lambda-2} & y^{\lambda-3} & y^{\lambda-4} & \dots & 1
\end{smallmatrix}\right).
\]
 Now observe that $\veps\circ x^{i}y^{j}F=y^{j}\cdot\veps\circ x^{i}F$.
Using this, we can write: 
\begin{align*}
\bigoplus_{i=1}^{\lambda-1}-\veps\circ h_{i-1}(x,y)F & =\bigoplus_{i=0}^{\lambda-2}\sum_{j+k=i}y^{k}\cdot(-\veps\circ x^{j}F)\\
 & =M_{h}'\ci\left(\bigoplus_{i=0}^{\lambda-2}-\veps\circ x^{j}F\right),
\end{align*}
 and it follows from this and the isomorphism above the claim that
the map of the claim is an isomorphism.
\end{proof}
Now assume $\lambda>0$ and reorder the summands of the domain and
codomain to permute the rows and columns of the matrix of $[\tilde{\rho}_{\lambda}]_{22}$.
Let the domain be given in the order: 
\[
FE_{\lambda-1}[y]^{\oplus2}\oplus A_{\lambda-1}[y]\oplus EF_{\lambda-1}[y]\oplus FEFE_{\lambda-1}[y],
\]
 where the first two identical summands appear in the same order as
before. Let the codomain be given in the order: 
\begin{multline*}
FE_{\lambda-1}[y]^{\oplus2}\oplus A_{\lambda-1}[y]\oplus FE_{\lambda-1}[y]\oplus A_{\lambda-1}[y]^{\oplus\lambda-1}\\
\oplus F^{2}E_{\lambda-1}^{2}[y]\oplus FE_{\lambda-1}[y]\oplus FE_{\lambda-1}[y]^{\oplus\lambda-1},
\end{multline*}
 where the new summand number (numbered left to right) and corresponding
old summand number are given precisely in the following chart: 
\[
\begin{array}{ccccccccccccccc}
\text{new:} & 1 & 2 & 3 & 4 & 5 & 6 & 7 & \dots & \lambda+3 & \lambda+4 & \lambda+5 & \lambda+6 & \dots & 2\lambda+5\\
\text{old:} & 4 & 1 & 6 & 3 & 7 & 8 & 9 & \dots & \lambda+5 & 2 & 5 & \lambda+6 & \dots & 2\lambda+5.
\end{array}
\]
 Writing the matrix of $[\tilde{\rho}_{\lambda}]_{22}$ for $\lambda>0$,
with columns and rows changed by the above permutations, we obtain:
{\footnotesize{}
\[
\begin{pmatrix}1 & 0 & 0 & 0 & 0\\
0 & 1 & 0 & 0 & 0\\
0 & 0 & 1 & 0 & 0\\
y_{1} & 0 & \eta & \sigma & 0\\
0 & 0 & \bigoplus\limits _{i=1}^{\lambda-1}y^{i} & \bigoplus\limits _{i=1}^{\lambda-1}-\veps\circ h_{i-1}(x,y)F & 0\\
0 & 0 & 0 & 0 & F\veps E\\
0 & 0 & 0 & 0 & F\sigma E\\
\bigoplus\limits _{i=0}^{\lambda-1}x^{i}E & \bigoplus\limits _{i=0}^{\lambda-1}Fx^{i} & \bigoplus\limits _{i=0}^{\lambda-1}h_{i-1}(x,y)\circ\eta & \Theta & \bigoplus\limits _{i=0}^{\lambda-1}F(\veps\circ x^{i}y_{1}F)E
\end{pmatrix}.
\]
} After compressing the notation of rows 4-5 and 6-8 of this matrix,
we obtain a lower-triangular matrix. The last two diagonal entries
are: 
\[
\begin{pmatrix}\sigma\\
\bigoplus\limits _{i=1}^{\lambda-1}-\veps\circ h_{i-1}(x,y)F
\end{pmatrix},
\]
 which is an isomorphism by the claim, and: 
\[
\left(\begin{smallmatrix}F\veps E\\
F\sigma E\\
\bigoplus\limits _{i=0}^{\lambda-1}F(\veps\circ x^{i}y_{1}F)E
\end{smallmatrix}\right):FEFE_{\lambda-1}[y]\to FE_{\lambda-1}[y]\oplus F^{2}E_{\lambda-1}^{2}[y]\oplus FE_{\lambda-1}[y]^{\oplus\lambda},
\]
 which is an isomorphism for $\lambda>0$, and therefore for $\lambda+1\geq0$,
using Claim \ref{claim:veps-x^iy_1} with $F$ applied on the left
and $E$ on the right.

When $\lambda=0$ the matrix of $[\tilde{\rho}_{\lambda}]_{22}$ is
given by removing rows 3, 5--$(\lambda+3)$, and $(\lambda+6)$--$(2\lambda+5)$:
\[
\begin{pmatrix}1 & 0 & 0 & 0 & 0\\
0 & 1 & 0 & 0 & 0\\
y_{1} & 0 & \eta & \sigma & 0\\
0 & 0 & 0 & 0 & F\veps E\\
0 & 0 & 0 & 0 & F\sigma E
\end{pmatrix}.
\]
 When $\lambda=0$ we also have isomorphisms: 
\[
(\eta,\sigma):EF_{\lambda-1}[y]\oplus A[y]_{\lambda-1}\iso FE_{\lambda-1}[y]
\]
 and 
\[
\left(\begin{smallmatrix}F\veps E\\
F\sigma E
\end{smallmatrix}\right):F\bigl(FE_{\lambda+1}\bigr)E[y]\iso F_{\lambda+1}E[y]\oplus F\bigl(EF_{\lambda+1}\bigr)E[y],
\]
so we see that again the matrix can be written as a lower-triangular
matrix with invertible diagonal entries.

Finally, assume $\lambda<0$. We have an isomorphism: 
\[
\left(\sigma,\sum_{i=0}^{-\lambda}Fx^{i}\ci\eta\right):EF_{\lambda-1}[y]\oplus A_{\lambda-1}[y]^{\oplus-(\lambda-1)}\iso FE_{\lambda-1}[y],
\]
 which is the isomorphism $\rho_{\lambda-1}\otimes_{k}k[y]$. There
is a final claim to check: 
\begin{claim}
When $\lambda<0$, the map 
\[
\left(\sigma,\eta,\sum_{i=0}^{-\lambda-1}-Fx^{i}y_{1}\ci\eta\right):EF_{\lambda-1}[y]\oplus A_{\lambda-1}[y]^{\oplus-(\lambda-1)}\to FE_{\lambda-1}[y]
\]
 is an isomorphism.
\end{claim}

\begin{proof}
Define an isomorphism $M_{-y}'\in\End_{A_{\lambda-1}[y]}\left(A_{\lambda-1}[y]^{\oplus-(\lambda-1)}\right)$
with components $[M_{h}]_{ij}$ given by $1$ along the diagonal and
$-y$ along the subdiagonal. This is a lower-triangular invertible
matrix. We write the map in question as a composition of isomorphisms:
\begin{align*}
 & \left(\sigma,\eta,\sum_{i=0}^{-\lambda-1}-Fx^{i}y_{1}\ci\eta\right)=\left(\sigma,\eta,\sum_{i=1}^{-\lambda}Fx^{i}\ci\eta\right)\\
 & \ci\begin{pmatrix}\idop_{EF_{\lambda-1}[y]} & 0 & 0\\
0 & \idop_{A_{\lambda-1}[y]} & 0\\
0 & 0 & -\idop_{A_{\lambda-1}[y]^{\oplus-\lambda}}
\end{pmatrix}\ci\begin{pmatrix}\idop_{EF_{\lambda-1}[y]} & 0\\
0 & M_{-y}'
\end{pmatrix}.
\end{align*}
\end{proof}
Now let $W$ be the endomorphism of the codomain of $[\tilde{\rho}_{\lambda}]_{22}$
given by the invertible matrix: 
\[
[W]=\left(\begin{smallmatrix}1 & 0 & 0 & 0 & 0\\
0 & 1 & 0 & 0 & 0\\
0 & 0 & 1 & -y_{1} & 0\\
0 & 0 & 0 & 1 & 0\\
0 & 0 & 0 & 0 & 1
\end{smallmatrix}\right).
\]
 We show that $[W]\cdot[\tilde{\rho}_{\lambda}]_{22}$ is equivalent
to a lower-triangular matrix after giving a suitable permutation of
the domain and codomain summands. Let the domain be given in the order:
\[
EF_{\lambda-1}[y]\oplus A_{\lambda-1}[y]^{-(\lambda-1)}\oplus FE_{\lambda-1}[y]\oplus FEFE_{\lambda-1}[y]\oplus FE_{\lambda-1}[y]^{-(\lambda+1)}\oplus FE_{\lambda-1}[y]^{\oplus2},
\]
 where the change of summand numbers is given by the following chart:
\[
\begin{array}{ccccccccc}
\text{new:} & 1 & 2 & 3 & 4 & \dots & -\lambda+2 & -\lambda+3 & -\lambda+4\\
\text{old:} & 5 & 1 & 6 & 7 & \dots & -\lambda+5 & 2 & 4\\
\\
\text{new:} & -\lambda+5 & -\lambda+6 & \dots & -2\lambda+4 & -2\lambda+5 & -2\lambda+6\\
\text{old:} & -\lambda+7 & -\lambda+8 & \dots & -2\lambda+5 & -\lambda+6 & 3.
\end{array}
\]
 Let the codomain be given in the order: 
\[
FE_{\lambda-1}[y]^{\oplus4}\oplus F^{2}E_{\lambda-1}^{2}[y],
\]
 where the change of summand numbers is given by the following chart:
\[
\begin{array}{cccccc}
\text{new:} & 1 & 2 & 3 & 4 & 5\\
\text{old:} & 3 & 4 & 5 & 2 & 1.
\end{array}
\]
 The matrix of $[W]\cdot[\tilde{\rho}_{\lambda}]_{22}$ for $\lambda<0$
agrees with that for $[\tilde{\rho}_{\lambda}]_{22}$ except in the
third row, where it is: 
\[
\begin{pmatrix}\eta & 0 & 0 & 0 & \sigma & \sum\limits _{i=0}^{-\lambda-1}-Fx^{i}y_{1}\circ\eta & 0 & 0\end{pmatrix}.
\]
 Writing now the matrix of $[W]\cdot[\tilde{\rho}_{\lambda}]_{22}$
with columns and rows changed by the above permutations, and compressing
the notation for some columns, we obtain: {\footnotesize{}
\[
\begin{pmatrix}\left(\sigma,\eta,\sum\limits _{i=0}^{-\lambda-1}-Fx^{i}y_{1}\circ\eta\right) & 0 & (0,0) & 0 & 0\\
\left(0,0,\sum\limits _{i=0}^{-\lambda-1}Fx^{i}\circ\eta\right) & 1 & (0,0) & 0 & 0\\
\left(0,0,\Theta'\right) & 0 & \left(F\sigma E,\sum\limits _{i=1}^{-\lambda-1}F^{2}h_{i-1}(x_{2},y)\circ F\eta E\right) & 0 & 0\\
\left(0,0,\sum\limits _{i=0}^{-\lambda-1}-Fh_{i-1}(x,y)\circ\eta\right) & 0 & \left(F\veps E,\sum\limits _{i=1}^{-\lambda-1}y^{i}\right) & 1 & 0\\
\left(0,0,\sum\limits _{i=0}^{-\lambda-1}Fy^{i}\circ\eta\right) & 0 & \left(0,\sum\limits _{i=1}^{-\lambda-1}y^{i}y_{1}\right) & y_{1} & 1
\end{pmatrix}.
\]
} The upper left map is an isomorphism by the Claim proved above.
The middle diagonal map is an isomorphism because it is the isomorphism
of Claim \ref{claim:eta-h(x,y)} with $F$ applied on the left and
$E$ on the right. So the matrix is lower-triangular with isomorphisms
along the diagonal.
\end{proof}
\end{itemize}
\bibliographystyle{amsalpha}
\providecommand{\bysame}{\leavevmode\hbox to3em{\hrulefill}\thinspace}
\providecommand{\MR}{\relax\ifhmode\unskip\space\fi MR }
\providecommand{\MRhref}[2]{%
	\href{http://www.ams.org/mathscinet-getitem?mr=#1}{#2}
}
\providecommand{\href}[2]{#2}

\end{document}